\magnification=\magstep1
\input amstex
\documentstyle{amsppt}
\catcode`\@=11 \loadmathfont{rsfs}
\def\mycal{\mathfont@\rsfs}
\csname rsfs \endcsname \catcode`\@=\active

\vsize=7.5in

\topmatter

\title
$W^*$-representations of subfactors \\ and restrictions on the Jones index \\ 
$\text{\it to Vaughan Jones, in memoriam}$
\endtitle

\author Sorin Popa\endauthor

\rightheadtext{Restrictions on the index}

\affil University of California Los Angeles \endaffil

\email popa\@math.ucla.edu \endemail

\thanks Supported in part by NSF Grant DMS-1955812  and the Takesaki Endowed Chair at UCLA\endthanks



\abstract  A {\it $W^*$-representation} of a  
II$_1$ subfactor $N\subset M$ with finite Jones index, $[M:N]<\infty$, 
is a non-degenerate commuting square embedding of $N\subset M$ into an inclusion of atomic von Neumann algebras 
$\oplus_{i\in I} \Cal B(\Cal K_i)=\Cal N \subset^{\Cal E} \Cal M=\oplus_{j\in J} \Cal B(\Cal H_j)$.  
We undertake here a systematic study of this notion, first introduced in [P92a], giving examples and considering invariants such as the (bipartite) 
{\it inclusion graph} $\Lambda_{\Cal N \subset \Cal M}$, the {\it coupling vector} $(\text{\rm dim}(_M\Cal H_j))_j$ 
and the {\it RC-algebra} (relative commutant) $M'\cap \Cal N$, for which we establish some basic properties. 
We then prove that if $N\subset M$ admits a $W^*$-representation $\Cal N\subset^{\Cal E}\Cal M$,   
with the expectation $\Cal E$ preserving a semifinite trace on $\Cal M$, such that there exists a norm one projection 
of $\Cal M$ onto $M$ commuting with $\Cal E$, a property of $N\subset M$ that we call {\it weak injectivity/amenability},   
then $[M:N]$ equals the square norm of the inclusion graph $\Lambda_{\Cal N \subset \Cal M}$.   
\endabstract 

\endtopmatter 

\document

\heading 1.   Introduction   \endheading

A  tracial von Neumann algebra $(M, \tau)$ is naturally represented as the algebra of left multiplication operators 
on the Hilbert space $L^2M$, obtained by completing $M$ in the norm $\|\cdot \|_2$ given by the trace $\tau$. This is called the {\it standard representation} of  $M$.   
In fact, $M$ acts on $L^2M$ on the right as well, giving  a representation of $M^{op}$, the opposite of the algebra $M$. The left-right multiplication algebras $M, M^{op}$
commute, one being the centralizer of the other, $M^{op}=M'$. 
Any other representation $M \subset \Cal B(\Cal H)$ of
$M$ as a von Neumann algebra, or {\it left Hilbert} $M$-{\it module}
$_M\Cal H$, is of the form $\Cal H \simeq \oplus_k L^2(M)p_k$, for
some projections $\{p_k\}_{k\in K} \subset \Cal P(M)$, with the action of $M$ by left
multiplication. If $M$ is a factor, 
then dim$(_M\Cal H)\overset{\text{\rm def}}\to{=}\sum_k \tau(p_k)\in [0,\infty]$ characterizes the
isomorphism class of $_M\Cal H$. The $M$-module $_M\Cal H$ can be alternatively
described as the (left) $M$-module $e_{00}L^2(M^\infty, \text{\rm Tr})p$, 
where $M^\infty$ is the II$_\infty$ factor 
$M\overline{\otimes} \Cal B(\ell^2K)$, $\text{\rm Tr}$ is its normal semifinite 
faithful (n.s.f.) trace $\tau\otimes \text{\rm Tr}_{\Cal B(\ell^2K)}$, 
$e_{00}\in \Cal B(\ell^2K)$ is a rank one projection and $p\in M^\infty$ 
is any projection satisfying $\text{\rm Tr}(p)=\text{\rm dim}(_M\Cal H)$. 
Thus, the commutant $M'$ of $M$ in $\Cal B(\Cal H)$ is a II$_1$
factor iff dim$(_M\Cal H) < \infty$ and if this is the case then 
$M'\simeq (M^t)^{op}$, where $t=\text{\rm dim}(_M\Cal H)=\text{\rm Tr}(p)$.  

This summarizes Murray-von Neumann famous theory of continuous dimension for type II factors. The  number $\text{\rm dim}(_M\Cal H)\in [0, \infty]$ is called 
the {\it dimension} of the $M$-module $_M\Cal H$. When viewed as the amplifying number $t$, measuring the ratio between the size of $M$ and $M'= (M^t)^{op}$ 
in $\Cal B(\Cal H)$, it is  called the {\it coupling constant}.  

We consider in this paper the analogue for  a subfactor $N\subset M$ of finite Jones index, $[M:N]<\infty$, of the representations    
of a single II$_1$ factor. This concept was introduced in (Section 2 of [P92a]), but we undertake here 
a systematic study of this notion and its applications. Thus, a $W^*$-{\it representation} of $N\subset M$ is a nondegenerate embedding of $N\subset M$ into an 
inclusion of  atomic von Neumann algebras with expectation,  
$\oplus_{i\in I} \Cal B(\Cal K_i)=\Cal N\subset^{\Cal E} \Cal M=\oplus_{j\in J} \Cal B(\Cal H_j)$.  
This means $M$ is embedded as a von Neumann algebra into $\Cal M$ such that: $N\subset \Cal N$, with the restriction of $\Cal E$ to $M$ equal to $E_N$, the $\tau$-preserving 
expectation of $M$ onto $N$ ({\it commuting square} condition); the span of $M\Cal N$ equals $\Cal M$ ({\it non-degeneracy} condition).

These conditions imply that any orthonormal basis of $M$ over $N$ is an orthonormal basis for $\Cal E$, so the index Ind$(\Cal E)$ of this expectation equals $[M:N]$ and the 
{\it inclusion} (bipartite) {\it graph} $\Lambda=\Lambda_{\Cal N\subset \Cal M}=(b_{ij})_{i\in I, j\in J}$ of $\Cal N\subset \Cal M$ satisfies $\|\Lambda\|^2\leq [M:N]$, 
where $b_{ij}$ gives the multiplicity of $\Cal B(\Cal K_i)$ in $\Cal B(\Cal H_j)$ and $\Lambda$ is  
viewed as an $I\times J$ matrix (or element in $\Cal B(\ell^2J, \ell^2I)$). 

The role of the Murray-von Neumann dimension, or coupling constant, 
is played here by the {\it dimension/coupling vectors} $\vec{d}_M=(\text{\rm dim}(_M\Cal H_j))_j$, 
$\vec{d}_N=(\text{\rm dim}(_N\Cal K_i))_i$. For these two vectors to have all entries finite it is sufficient that one of the entries is finite, and if this is the case then 
$\Lambda^t\Lambda(\vec{d}_M)=[M:N]\vec{d}_M$,   $\vec{d}_N=\Lambda(\vec{d}_M)$. 

Thus, while for a single factor all representations are stably equivalent, $W^*$-representations of subfactors  
appear, from the outset, as a far more complex notion. 

Indeed, a subfactor $N\subset M$ admits a large variety of $W^*$-representations, with a central role played by  
the {\it standard representation}, $\Cal N^{st}\subset^{\Cal E^{st}}\Cal M^{st}$, whose inclusion graph $\Lambda_{\Cal N^{st}\subset \Cal M^{st}}$ 
coincides with the transpose of the standard graph $\Gamma_{N\subset M}$ of $N\subset M$, with the coupling vectors given by the standard weights of $\Gamma_{N\subset M}$.  

This leads right away to natural concepts of injectivity and amenability, by analogy with the single II$_1$ factor case. Thus, an inclusion of II$_1$ factors $N\subset M$ 
is {\it injective} if there exists an expectation of $(\Cal N^{st}\subset^{\Cal E^{st}}\Cal M^{st})$ onto $(N\subset M)$, i.e. 
a norm one projection $\Phi:\Cal M^{st}\rightarrow M$ that's $\Cal E^{st}$-invariant. It is {\it amenable} if there exists  
a $(N\subset M)$-{\it hypertrace} on $(\Cal N^{st}\subset^{\Cal E^{st}}\Cal M^{st})$, i.e., a state $\varphi$ on $\Cal M^{st}$ that's $\Cal E^{st}$-invariant and has $M$ in its centralizer. 
These two concepts are easily seen to be equivalent and they imply $N, M$ are amenable/injective as 
single factors, thus being isomorphic to the hyperfinite II$_1$ factor $R$ by Connes theorem ([C76]). They also imply that $N\subset M$ is the range of a norm-one projection in {\it any} 
of its $W^*$-representations $(N\subset M)\subset (\Cal N\subset^{\Cal E}\Cal M)$, once some natural compatibility of higher relative commutants is satisfied ({\it smoothness}). 

Amenability of $N\subset M$ was introduced in ([P92a], [93a], [P97a]) and shown there to be equivalent to the condition $N, M \simeq R$ 
and $\Gamma_{N\subset M}$ amenable, i.e. $\|\Gamma_{N\subset M}\|^2=[M:N]$. It was also shown equivalent to the fact that $N\subset M$ can be exhausted by 
higher relative commutants of a ``$(N\subset M)$-compatible tunnel'' of subfactors, $M\supset N \supset P_1 \supset P_2 \supset ...$, obtained 
by iterative choices of the downward basic construction for induction/reduction by projections $p\in P_n'\cap M_{k(n)}$, at each step $n$, thus being completely classified 
by its standard invariant, $\Cal G_{N\subset M}$. We revisit these results in Section 4 of the paper (see Theorem 4.5). 

The complexity of $W^*$-representation theory for a subfactor naturally leads to several weaker amenability properties as well. 
Thus, a subfactor of finite Jones index $N\subset M$ 
is {\it weakly injective} (resp. {\it weakly amenable}) if it admits a 
$W^*$-representation $(N\subset M)\subset (\Cal N\subset^{\Cal E}\Cal M)$ with an $\Cal E$-invariant normal semifinite faithful ({\it n.s.f.}) trace,  
such that $(N\subset M)$  is the range of a norm-one projection from $(\Cal N\subset^\Cal E \Cal M)$  (resp. 
if $(\Cal N\subset^{\Cal E}\Cal M)$ has a $(N\subset M)$-hypertrace).  
One should note that the standard representation 
does satisfy this ``Traciality'' property (more on this below), so amenability/injectivity does imply weak amenablity/injectivity. 

Again, these two notions are equivalent and they imply $N, M\simeq R$ (by [C76]).  
They are also hereditary, i.e., if $(Q\subset P)\subset (N\subset M)$ is a non-degenerate commuting square of subfactors and $N\subset M$ is weakly amenable, then so is $Q\subset P$. 
A main result in this paper is the following: 

\proclaim{1.1. Theorem} If an extremal $($e.g. irreducible$)$ inclusion of $\text{\rm II}_1$ factors with finite index $N\subset M$ is weakly amenable, then $[M:N]$ 
is the square norm of a bipartite graph. More precisely, if $(N\subset M)\subset (\Cal N\subset^{\Cal E}\Cal M)$ is a Tracial $W^*$-representation    
so that $N\subset M$ is the range of a norm-one projection, then $[M:N]=\|\Lambda_{\Cal N \subset \Cal M}\|^2$. 
\endproclaim  

The above result indicates that  $W^*$-representations may help detect (and explain!) 
restrictions on the set $\mycal C(M)$ of indices of irreducible subfactors of a given II$_1$ factor $M$,  especially  for $M=R$

Recall in this respect Jones fundamental result in [J82], showing that the index of any subfactor lies in the spectrum $\{4\cos^2 (\pi/n)\mid n\geq 3\}\cup [4, \infty)$. 
One of his proofs of this result amounts to showing that if $[M:N]<4$ then it must equal the square norm of the standard graph $\Gamma_{N\subset M}$, and using 
the fact that the set $\Bbb E^2$ of square norms of bipartite graphs has only the values $4\cos^2(\pi/n)$ when $<4$ ([J86]). 

The set $\Bbb E^2$  contains the half line $[2+\sqrt{5}, \infty)$, but $\Bbb E^2\cap [1, 2+\sqrt{5}]$ is a closed countable set, 
consisting of an increasing sequence of accumulation points converging to $2+\sqrt{5}$, the first of which being $4=\lim_n 4\cos^2(\pi/n)$ (cf [CDG82]; see also [GHJ88]). 
Yet it is known that any number $>4$ can occur as index of an irreducible subfactor ([P90], [P94]), and that $\mycal C(L\Bbb F_\infty)=\{4\cos^2(\pi/n)\mid n\geq 3\}\cup [4, \infty)$ ([PS01]). 
On the other hand, if $M$ is constructed out of a free ergodic probability measure preserving action of a non-elementary hyperbolic group, 
then $\mycal C(M)=\{1, 2, 3, ...\}$ ([PV13]). So $\mycal C(M)$ appears to depend in very subtle ways on the nature of the factor $M$. 

But the most important question along these lines, of calculating $\mycal C(R)$, remained open. Our work in this paper  
attempts to provide some tools for approaching the ``restrictions'' part of this   problem, more specifically for showing that $\mycal C(R)\subset \Bbb E^2$. 

By [H93], if an irreducible subfactor $N\subset M$ satisfies $4<[M:N]\leq 2+\sqrt{5}$, 
then its standard graph $\Gamma_{N\subset M}$ equals $A_\infty$. Equivalently, the higher relative commutants 
$N'\cap M_n$ in the Jones tower $N\subset M \subset_{e_0} M_1 \subset_{e_1} M_2...$ are generated by the Jones projections $e_0, e_1, ...$.  
So in order to show $\mycal C(R)\subset \Bbb E^2$, it is sufficient to prove that any $A_\infty$-subfactor of $R$ has index equal to the square norm of a 
bipartite graph. Our belief is that in fact any hyperfinite $A_\infty$-subfactor is weakly amenable. If true, then Theorem 1.1 above would imply that $\mycal C(R)\subset \Bbb E^2$.  

In order to make this speculation more specific, we need to fix some terminology and explain ways of producing $W^*$-representations.  
 
Thus, a $W^*$-representation $\Cal N\subset \Cal M$ is {\it irreducible} if $\Cal Z(\Cal M)\cap \Cal Z(\Cal N)=\Bbb C$. This is equivalent to $\Lambda_{\Cal N\subset \Cal M}$ 
being connected as a graph (irreducible as a matrix).

The $W^*$-representation is {\it Tracial} if it admits a n.s.f. trace Tr that's $\Cal E$-invariant. If 
$(N\subset M)\subset (\Cal N \subset^{\Cal E} \Cal M)$ is Tracial and has finite couplings, with the n.s.f. trace Tr given by $\vec{d}_M$ being $\Cal E$-preserving, 
then the $W^*$-representation is {\it canonically Tracial}. 

The simplest example of a $W^*$-representation occurs from a {\it graphage} of $N\subset M$, i.e., a non-degenerate commuting square 
$(Q\subset P)\subset (N\subset M)$, with $Q, P$ finite dimensional, by taking the basic construction inclusion $(N\subset M)\subset (\Cal N =\langle N, e^M_P\rangle \subset 
\langle M, e^M_P\rangle=\Cal M)$ (cf. 2.1 in [P92a]; see 3.2.2 below). It is easy to see that such $W^*$-representations have finite couplings and are canonically Tracial. 

Another isomorphism invariant for a $W^*$-representation $\Cal N\subset^{\Cal E} \Cal M$, 
besides the inclusion graph $\Lambda_{\Cal N\subset \Cal M}$ and the coupling vectors, is the isomorphism class of the {\it relative commutant} ({\it RC}) algebra 
$M'\cap \Cal N$. If $\Cal N\subset \Cal M$ has finite couplings, then $M'\cap \Cal N$ identifies naturally to a von Neumann subalgebra of $(M^t)^{op}$, so it is finite. 

A $W^*$-representation $(N\subset M) \subset (\Cal N\subset^{\Cal E}\Cal M)$ is {\it exact}, if $M\vee (M'\cap \Cal N)=\Cal M$. 
Such a representation is irreducible if and only if its RC-algebra $M'\cap \Cal N$ is a factor. 

We are  particularly interested in irreducible exact $W^*$-representation with finite couplings. They are all 
of the form $\Cal N_P:= N\vee P^{op} \subset^{\Cal E_P} M\vee P^{op}=:\Cal M_P$, where $P\subset M$ is an irreducible subfactor and  
$M$ acts here by left multiplication on $L^2(M_\infty)$,  $P$ acts on the right, and $\Cal E_P$ is the unique expectation extending $E_N\otimes id_{P^{op}}$, 
$M_\infty$ denoting the enveloping algebra of the Jones tower for $N\subset M$. 

The standard $W^*$-representation $(N\subset M)\subset (\Cal N^{st} \subset^{\Cal E^{st}} \Cal M^{st})$, corresponds to the case $P=M$ of this construction, 
i.e. to $(N\subset M)\subset (N\vee M^{op}\subset M\vee M^{op})$. 

Taking the direct sum $\oplus_P (\Cal N_P \subset^{\Cal E_P} \Cal M_P)$, over all isomorphism classes of irreducible $P\subset M$, 
one obtains the universal exact $W^*$-representation with finite couplings $(N\subset M)\subset (\Cal N^{u,fc}\subset^{\Cal E^{u,fc}}\Cal M^{u,fc})$. 
We say that $N\subset M$ is {\it ufc-amenable} (resp. {\it ufc-injective}) if this representation has a $(N\subset M)$-hypertrace (resp. norm-one projection). 
Our conjecture is that any hyperfinite $A_\infty$-subfactor $N\subset M$ is ufc-amenable,  and thus $[M:N]=\|\Lambda_{\Cal N^{u,fc}\subset \Cal M^{u,fc}}\|^2$. 

The paper is organized as follows. Section 2 recalls some basic facts in Jones theory of subfactors, and its generalization for arbitrary  $W^*$-inclusions. 
In Section 3 we develop the concept of $W^*$-representation of subfactors, give formal definitions, prove general results and provide examples. In Section 4 we 
recall from ([P92a], [P97a]) the definition of amenability and injectivity for a subfactor $N\subset M$, as well as a result establishing the equivalence of amenability/injectivity with a 
series of other properties of $N\subset M$ (see Theorem 4.5).  In Section 5 we introduce the concepts of weak-amenability and weak-injectivity and prove Theorem 1.1 (see 5.4). 
We also define here ufc-amenability and ufc-injectivity of a subfactor $N\subset M$ (see 5.7) 
and state a result establishing the equivalence of these two properties with several 
other structural properties of $N\subset M$ (Theorem 5.8), which we will prove in a follow up to this paper. Section 6 contains many comments and open problems.

We mention that this paper is an outgrowth of our unpublished 1997 note entitled ``Biduals associated to subfactors, hypertraces and restrictions for the index''.

\heading 1. Preliminaries \endheading

We recall in this section some basic facts about Jones' index theory for inclusions of II$_1$ factors and, more generally, for inclusions of von Neumann algebras. 
We typically use the notations $M, N, P, Q, B$ for tracial von Neumann algebras,  
with the generic notation $\tau$ for the corresponding (faithful normal) trace state on it, sometimes with an index specifying the algebra on which it is defined 
(e.g. $\tau_M$ for the trace on $M$). The generic notation for arbitrary von Neumann algebras will be $\Cal M, \Cal N, \Cal P, \Cal Q$, etc. 

For basics on II$_1$ factors and tracial von Neumann algebras we refer to [AP17] and for general von Neumann algebras 
to [T03]. 

\vskip.05in

\noindent
{\bf 2.1. Jones index for subfactors and  the basic construction} ([J82]). If $B\subset M$ 
is a von Neumann subalgebra of the tracial von Neumann algebra $M$, then $E_B=E^M_B$ denotes the unique 
trace preserving conditional expectation of $M$ onto $B$. It extends to a projection $e_B$ of $L^2M$ onto $L^2B$, which is selfadjoint and positive on 
$L^2M$, when viewed as the space of square summable operators affiliated to $M$. 

Identifying $M$ 
with its standard representation $M \hookrightarrow \Cal B(L^2M)$ (as left multiplication operators on $L^2M$),   
one has: 

$(i)$ $e_Bxe_B=E_B(x)e_B$, $\forall x\in M$;  

$(ii)$ $B=\{x\in M\mid [x, e_B]=0\}$;  $B\ni b \mapsto be_B$ is an onto isomorphism; 

$(iii)$ $\vee_{u\in \Cal U(B)} ue_Bu^* =1.$

\noindent 
The von Neumann algebra $\langle M, e_B\rangle$ (or $\langle M, B\rangle$) generated in $\Cal B(L^2M)$  by $M$ and $e_B$ is called the {\it extension of $M$ by $B$}. 
It coincides with $(B^{op})'=(J_MBJ_M)'\cap \Cal B(L^2M)$ (the commutant of the operators of right multiplication by elements in $B$). 
The span of $\{xe_B y \mid x, y \in M\}$ is a wo-dense $^*$-subalgebra of $\langle M, e_B \rangle$ whose wo-closure contains $M$. Thus, this wo-closure is equal 
to $\langle M, e_B\rangle$. 
This construction of the new inclusion $M \subset \langle M, e_B\rangle$ from an initial inclusion $B\subset M$, is called {\it Jones basic construction}. 

Given an inclusion of II$_1$ factors $N\subset M$, the 
 {\it Jones index of $N$ in $M$}, denoted $[M:N]$, is defined as  the Murray-von Neumann coupling constant of $N$ when acting on 
 the Hilbert space $L^2M$ by left multiplication 
operators (as a subalgebra of $M$), or equivalently the Hilbert-dimension of $L^2M$ as a (left) $N$-Hilbert module, $[M:N]=\text{\rm dim}_NL^2M$. 

Thus, $[M:N]<\infty$ iff 
$N'\cap \Cal B(L^2M)$ is a II$_1$ factor, while $[M:N]=\infty$ iff $N'$ is of type II$_\infty$. If the index is finite, then $M\subset  M_1:=\langle M, e_N \rangle$ 
is an inclusion of II$_1$ factors with index $[M_1:M]=[M:N]$ and the trace state on $M_1$ satisfies $\tau_{M_1}(xe_Ny)=\lambda \tau_M(xy)$, 
$\forall x, y\in M$, where $\lambda = [M:N]^{-1}$. 

The index of subfactors is multiplicative, in the sense that if $P\subset N \subset M$ are subfactors, then $[M:P]=[M:N][N:P]$. 

Letting $(M_{-1} \subset M_0)=(N\subset M)$ and $e_0=e_N$, this allows constructing iteratively a whole {\it tower} of inclusions of II$_1$ factors, 
$M_{-1} \subset M_0\subset_{e_0} M_1 \subset_{e_1} \subset M_2 \subset ....$, with each $M_{i+1}$, $i\geq 0$, being generated by $M_i$ and a projection $e_i$ of trace $\lambda=[M:N]^{-1}$, 
having index $[M_{i+1}:M_{i}]=[M:N]$ and satisfying the properties: 

$(a)$ $e_ixe_i=E^{M_i}_{M_{i-1}}(x)e_i$, $\forall x\in M_i$; 

$(b)$ $\{e_i\}'\cap M_i=M_{i-1}$;  

$(c)$ $\tau(xe_i)=\lambda \tau(x)$, $\forall x\in M_i$. 

\noindent
In particular, the {\it $\lambda$-sequence of Jones projections} $\{e_i\}_{i\geq 0}$ with the trace $\tau$ satisfy the  conditions: 

$(a')$ $e_ie_{i\pm 1} e_i=\lambda e_i$; 

$(b')$ $[e_i, e_j]=0$, $\forall j> i+1$; 

$(c')$ $\tau(xe_{i+1})=\lambda \tau(x)$, $\forall x\in Alg(\{e_0, e_1, ..., e_i\})$. 

\vskip.05in

Jones celebrated theorem shows that the conditions $(a')-(c')$ imply the restrictions on the index $[M:N]=\lambda^{-1} \in \{4\cos^2(\pi/n) \mid n\geq 3\} \cup [4, \infty)$, 
with this latter set called the {\it Jones spectrum}.  The key ingredient in his proof of these restrictions on the index is the fact that 
axioms $(a')-(c')$ above imply that $A_n=C^*(1, e_0, ..., e_n)$ is finite dimensional, with $e_0 \vee ... \vee e_n$ a central projection in it, with the trace of its 
complement equal to $P_{n+1}(\lambda)=P_{n}(\lambda)-\lambda P_{n-1}(\lambda)$ whenever $P_{n}(\lambda), P_{n-1}(\lambda) >0$, where the polynomials 
$P_n(t)$, $n \geq -1$, are defined recursively by the formulas $P_{-1}=1, P_0=1$, $P_{n+1}(t)=P_{n}(t) - t P_{n-1}(t)$, $n\geq 0$.

\vskip.05in
\noindent
{\bf 2.2. Extremal subfactors}. It is shown in [J82] that if $N\subset M$ is a subfactor of finite index and $p$ is a projection in $N'\cap M$, then $[pMp:Np]=\tau_M(p)\tau_{N'}(p)[M:N]$ 
({\it Jones local index formula}). 
Thus, if $p_1, ..., p_n \in N'\cap M$ is a partition of 1 with projections, then $[M:N]=\sum_i \frac{[p_iMp_i:Np_i]}{\tau_M(p_i)}$. Note that in particular, this implies 
$\text{\rm dim}(N'\cap M)\leq [M:N]<\infty$, and more generally $ \text{\rm dim}(M_{i}'\cap M_j)\leq [M_j:M_i]=[M:N]^{j-i}<\infty$. 

Following (1.2.5 in [P92a]), an inclusion of II$_1$ factors $N\subset M$ is called {\it extremal} if $\tau_M(p)=\tau_{N'}(p)$ for all $p\in \Cal P(N'\cap M)$, or equivalently 
$[pMp:Np]=[M:N]\tau(p)^2$, $\forall p\in \Cal P(N'\cap M)$ non-zero. Recall from  (Corollary 4.5 in [PP84]) that this condition is equivalent 
to $E_{N'\cap M}(e)=\lambda 1$ for any Jones projection $e\in M$ (i.e., a projection whose expectation on $N$ is equal to $\lambda 1$). By ([P97b] or the Appendix in [P97a]) 
this is also equivalent to the fact that the norm closure of the convex hull of $\{ueu^* \mid u\in \Cal U(N)\}$ contains $\lambda 1$. 

Irreducible subfactors are automatically extremal. As pointed out in ([PP84]), 
if $N\subset M$ is extremal and $4< [M:N] < 3+2\sqrt{2}$, then $N\subset M$ is irreducible. 
So for  subfactors with small index $>4$, extremality is same as irreducibility.  

\vskip.05in
\noindent
{\bf 2.3. The standard invariant and graph of a subfactor}. The higher relative commutants in the Jones tower $\{M_i'\cap M_j\}_{j\geq i\geq -1}$ 
form a lattice of finite dimensional von Neumann algebras, with a trace $\tau$ inherited from $\cup_i M_i$, with 
inclusions $M_i'\cap M_j \subset M_l'\cap M_k$ whenever $-1\leq l \leq i \leq j \leq k$. Moreover, the Jones projection $e_j$ lies in $M_{j-1}'\cap M_{j+1}$ 
and implements the $\tau$-preserving expectation of $M_{i}'\cap M_j$ onto $M_{i}'\cap M_{j-1}$, $\forall j>i$. Also, any two such expectations commute, 
$E_{M_i'\cap M_j}E_{M_k'\cap M_l}=E_{M_k'\cap M_l}E_{M_i'\cap M_j}=E_{M_k'\cap M_j}$. 
In other words $\{M_i'\cap M_j\}_{j\geq i\geq -1}$ is a lattice of {\it commuting square} inclusions. 

The lattice of higher relative commutants, with its trace and Jones projections, is clearly an isomorphism invariant 
for $N\subset M$. It is called the {\it standard invariant of} $N\subset M$ and denoted $\Cal G_{N\subset M}$. 
See [P94] for more on this object and of a way to axiomatize it as an abstract object, called {\it standard $\lambda$-lattice}. 
It is pointed out in [P94] that, due to the  {\it duality} result in [PP84], which shows that there is a natural identification between the tower   
$(M_{i+1}\subset M_{i+2} \subset ...)$ and the $[M:N]$-amplification of $(M_{i-1} \subset M_i \subset ...)$, 
all informations is in fact contained in the first two rows of the lattice, i.e.,  $\Cal G_{N\subset M}=(\{M_i'\cap M_j\}_{j\geq i=0,-1}, \{e_i\}_{i\geq 0}, \tau)$

The 1st row of consecutive inclusions $\Bbb C=M_0'\cap M_0 \subset M_0'\cap M_1 \subset ....$ in $\Cal G_{N\subset M}$ 
is determined by a pointed bi-partite connected graph, called the {\it standard} (or {\it principal}) {\it graph}  
of $N\subset M$ and denoted $\Gamma_{N\subset M}=(a_{kl})_{k\in K, l\in L}$, with the {\it even vertices} being indexed by a set $K$ containing the ``initial'' vertex $*$, 
the {\it odd vertices} indexed by a set $L$, with $a_{kl}\geq 0$ denoting the number of edges between $k$ and $l$ (see e.g. Section 1.3.5 in [P92a]). 

To explain this in details, let $K_0=\{*\}$, $L_i=\{l\in L \mid \exists a_{kl} \neq 0, k\in K_{i-1}\}$, 
$K_i=\{k\in K\mid \exists a_{kl}\neq 0, l\in L_i\}$, $i\geq 1$, and note that $K=\cup_{i\geq 0}K_i$, $L=\cup_{i\geq 1}L_i$ 
(because $\Gamma_{N\subset M}$ is connected). The set of irreducible components of $M'\cap M_{2i}$, $i\geq 0$ (resp. $M'\cap M_{2i-1}$, $i\geq 1$) is identified with 
$K_i$ (resp. $L_i$), with the embeddings $K_i \subset K_{i+1}$, $i\geq 0$ (resp. $L_i \subset L_{i+1}$, $i\geq 1$)  being implemented by the $1$ to $1$ 
map from $\Cal Z(M'\cap M_{2i})$ into $\Cal Z(M'\cap M_{2i+2})$ (resp. $\Cal Z(M'\cap M_{2i-1})$ into $\Cal Z(M'\cap M_{2i+1})$): 
$$
\Cal Z(M'\cap M_{2i})\ni z \mapsto \text{\rm unique} \  z'\in \Cal Z(M'\cap M_{2i+2}) \  \text{\rm with} \ ze_{2i+1} = z'e_{2i+1}, 
$$
(similarly for $L_i \subset L_{i+1}$). Each index $k\in K$ (resp. $l\in L$) can also be viewed as labeling the 
irreducible subfactor $Mp\subset pM_{2i}p$ (resp. $Mq\subset qM_{2i-1}q$) with $p$ a minimal 
projection in the $k$'th direct summand of $M'\cap M_{2i}$ (resp. $l$'th summand of $M'\cap M_{2i-1}$) 
for any $i\geq 0$ with $K_i \ni k$ (resp. any $i\geq 1$ with $L_i\ni l$). 

With these conventions, the bipartite graph (diagram) for the embedding $M'\cap M_{2i}\subset M'\cap M_{2i+1}$, $i\geq 0$  
(resp. $M'\cap M_{2i-1}\subset M'\cap M_{2i}$, $i\geq 1$) is given by $\Gamma_{|K_i}=(a_{kl})_{k\in K_i, l\in L_{i+1}}$ (resp. $\Gamma^t_{|L_i}=(b_{lk})_{l\in L_i, k\in K_i}$, 
where $b_{lk}=a_{kl}$). 

Let us now assume $N\subset M$ is extremal. In this case,  the trace $\tau$ on the finite dimensional algebras $M_i'\cap M_j$ is uniquely determined by the {\it standard vectors} 
$\vec{v}=(v_k)_{k\in K}$, $\vec{u}=(u_l)_{l\in L}$, given by square roots of 
indices of irreducible inclusions appearing in the Jones tower, as follows. 

First note that Jones local index formula combined with the duality imply that if $k\in K$ (resp $l\in L$) 
then the index $[pM_{2i}p: Mp]$ (resp $[pM_{2i-1}p: Mp]$) is the same for any $i$ with the property that $k\in K_i$ (resp $l\in L_i$) 
and any $p$ minimal projection in the $k$'th (resp $l$'th) summand of $M'\cap M_{2i}$ (resp. of $M'\cap M_{2i-1}$).  
One defines $v_k= [pM_{2i}p: Mp]^{1/2}$ and $u_l=[pM_{2i-1}p: Mp]^{1/2}$. Thus, $v_{*}=1$ and if $N'\cap M=\Bbb C$ 
then the single point set $L_1=\{l_1\}$ satisfies $u_{l_1}=[M:N]^{1/2}$. Also,  $\vec{u}=\Gamma^t (\vec{v})$, $\Gamma\Gamma^t(\vec{v})=\lambda^{-1}\vec{v}$, 
where $\Gamma=\Gamma_{N\subset M}=(a_{kl})_{k,l}$ is now viewed as a $K \times L$ matrix. 

The trace of a minimal projection $p$ in the $k$'th summand of $M'\cap M_{2i}$ (for $i$ large enough so that $k\in K_i$) is then given by $\tau(p)=\lambda^i v_k$. 
Similarly, if $q$ is a minimal projection in the $l$'th summand of $M'\cap M_{2i+1}$ then $\tau(q)=\lambda^{i+\frac{1}{2}} u_k$.  

The norm of the standard graph of the subfactor coincides with the growth rate of the higher relative commutants 
$\|\Gamma_{N\subset M}\|=\underset{n\rightarrow \infty}\to{\lim} (\text{\rm dim}(M'\cap M_n))^{1/n}$ and satsfies 
the estimate $\|\Gamma_{N\subset M}\|^2  \leq [M:N]$ (see 1.3.5 in [P92a]).

\vskip.05in
\noindent
{\bf 2.4. Alternative definitions of the index}. Let us recall two alternative ways of defining the index, from ([PP84]). 

Given a von Neumann subalgebra $B$ of a tracial von Neumann algebra $(M, \tau)$, there exists a family of 
elements $\{\xi_i\}_i\subset L^2M$ such that $E_B(\xi_j^*\xi_i)=\delta_{ij}p_i\in \Cal P(B)$ and $x=\sum_i \xi_i E_B(\xi_i^*x)$ for all $x\in M$, where the convergence 
is in $L^2M$. Equivalently, $\{\xi_i\}_i \subset L^2M$ are so that $\xi_i e_B$ are partial isometries with $\sum_i \xi_i e_B \xi_i^*=1$. 

Such $\{\xi_i\}_i \subset L^2M$ is called an {\it orthonormal basis} of $M$ over $B$, and it is unique in an appropriate sense (cf. Proposition 1.3 in [PP84]). 
The finite partial sums in $\sum_i \xi_i\xi_i^*$ lie in $L^1M_+$ and they form an increasing net with ``limit'' 
$Z=\sum_i \xi_i \xi_i^*=Z_0p_0+(\infty)(1-p_0)$, where $p_0\in \Cal P(\Cal Z(M))$ 
and $Z_0$ is a positive densely defined operator affiliated with $\Cal Z(M)$,  with $Z=Z(B\subset M)$ independent 
of the choice of the o.b. Thus, if $M$ is a II$_1$ factor, then $Z=\sum_i \xi_i \xi_i^* = \alpha 1$ with 
$\alpha\in [1, \infty]$ only depending on the isomorphism class of $B\subset M$. 

In case $N\subset M$ is an inclusion of II$_1$ factors, then by (1.3.3 in [PP84]) one has $[M:N]=Z(N\subset M)$, in other words 
$[M:N]=\sum_i m_im_i^*$ for any o.b. $\{m_i\}_i$ of $M$ over $N$.  Thus, if $[M:N]<\infty$, then any 
o.b. is made up of ``bounded elements'' $m_i \in M$, with $\|m_i\|\leq [M:N]^{1/2}$, 
and one can in fact choose them all but possibly one so that $E_N(m_i^*m_i)=1$. 

Another characterization of  the Jones index is given in (Theorem 2.2 in [PP84]), where it is shown that 
the quantity $\lambda(E_N):=\sup \{c\geq 0 \mid E_N(x)\geq c x, x\in M_+\}$, which measures the  ``flattening''  
of the expected $N$-value of positive elements in $M$, satisfies $\lambda(E_N)=[M:N]^{-1}$. This characterization 
is key to calculating the Connes-St\o rmer {\it relative entropy} $H(M|N)$ ([CS75]) of a subfactor $N\subset M$ (Theorem 4.6 in [PP84]), 
which in particular provides  the characterization ``$N\subset M$ extremal iff  $H(M|N)=\ln [M:N]$''.

\vskip.05in

\noindent
{\bf 2.5. $W^*$-inclusions with finite index}. The above two formulas for the Jones index of II$_1$ subfactors, which only depend  
on the expectation from the ambient algebra onto its subalgebra, give the possibility of defining the index for arbitrary inclusions of von Neumann algebras 
(or a $W^*$-{\it inclusion}) with conditional expectation, $\Cal N \subset^{\Cal E} \Cal M$, in two alternative ways. 

The easiest to define is as follows. Denote $\lambda(\Cal E)\overset{\text{\rm def}}\to{=}  
\sup \{c\geq 0 \mid \Cal E(x)\geq c x, \forall x\in \Cal M_+\}$ and define Ind$(\Cal E)= \lambda(\Cal E)^{-1}$. We call the latter  the 
{\it index} of $\Cal E$ (also referred to as the {\it probabilistic index}, or the [PP84]-{\it index}, of the expectation $\Cal E$). 

Note that if Ind$(\Cal E) < \infty$ then $\Cal E$ is automatically normal and faithful (see e.g. Lemma 1.1 in [P97b]). As we have seen in 2.4,   
if $N\subset M$ are II$_1$ factors with $\Cal E=E_N$ the trace preserving expectation, then Ind$(E_N)=[M:N]$. 

This definition has many advantages for various limiting arguments (see below). But it is useful to view it in combination with the 
definition of the index based on o.b. of $\Cal M$ over $\Cal N$, with respect to the expectation $\Cal E$. 
Since this is closely related to the notion of basic construction for arbitrary $W^*$inclusions $\Cal N\subset^{\Cal E}\Cal M$, 
we recall from [T72] that if $\Cal M\subset \Cal B(\Cal H)$ 
is the standard representation of  the von Neumann algebra $\Cal M$, then there exists a ``canonical'' projection $e=e_{\Cal N} \in \Cal B(\Cal H)$ 
such that 

$(i)$ $exe=\Cal E(x)e$, $\forall x\in \Cal M$; 

$(ii)$ $\Cal N= \{e\}'\cap \Cal M$; $x\in \Cal N$, $xe=0$ iff $x=0$; 

$(iii)$ $\vee \{ueu^* \mid u\in \Cal U(\Cal M)\}=1$;  

\noindent
with $e$ unique with these properties, up to spatial isomorphism. One denotes by $\langle \Cal M, e\rangle$ the von Neumann algebra generated by $\Cal M$ and $e$, 
which has sp$\{xey\mid x, y\in \Cal M\}$ as a wo-dense $^*$-subalgebra.   This is the {\it basic construction} for arbitrary $W^*$-inclusions with expectation.  
We still call $e$ the {\it Jones projection} implementing $\Cal E$. We see by this definition that $e\langle \Cal M, e\rangle e =\Cal Ne\simeq \Cal N$, 
so a ``corner'' of central support 1 of $ \langle \Cal M, e\rangle$ is equal to $\Cal N$. 

There does exist in this generality an analogue of o.b. of $\Cal M$ over $\Cal N$ with respect to $\Cal E$. For our purposes, it is sufficient to discuss 
this under two types of assumptions: when either Ind$(\Cal E)< \infty$; or when $\Cal M$ is finitely generated as a right $\Cal N$ module. 
In both cases it is immediate to show that there exists a family $\{m_j\}_j\subset \Cal M$ 
such that $\Cal E(m_i^*m_j)=\delta_{ij}p_j \in \Cal P(\Cal N)$ and $\sum_j m_jem_j^*=1$, that we  call an {\it orthonormal basis} (o.b.) of $\Cal N \subset^{\Cal E} \Cal M$. 

One has $\text{\rm Ind}(\Cal E) \leq \|\sum_j m_jm_j^*\| \leq (\text{\rm Ind}(\Cal E))^2$ (cf. 1.1.6 in [P93b]). 
Thus, Ind$(\Cal E)<\infty$ iff $Z(\Cal E) := \sum_j m_jm_j^*$ is bounded, in which case 
this element, which lies in $\Cal Z(\Cal M)$ and is $\geq 1$, doesn't depend on the o.b.. 

The  o.b. properties imply that each $X\in 
\langle \Cal M, e\rangle$ can be written as $X=\sum_{i,j}m_i y_{ij}em_j^*$, for some unique $y_{ij}\in p_i\Cal Np_j$, i.e., 
$\langle \Cal M, e\rangle$ is an ``amplification'' of $\Cal N$, with $\Cal Z(\Cal N)$ naturally identifying with $\Cal Z(\langle \Cal M, e\rangle)$, 
via the map $\Cal Z(\Cal N)\ni z \mapsto z'$, where $z'$ is the unique element in $\Cal Z(\langle \Cal M, e\rangle)$ such that $z'e=ze$. 

If either $\Cal N, \Cal M$ are properly infinite, or if they are type II$_1$, or if they are both type I$_{fin}$ but with each 
type I$_n$ summand of $\Cal N$ having multiplicity $\leq n$ in each homogeneous type I$_m$ summand of $\Cal M$, then one actually has 
Ind$(\Cal E)=\|\sum_j m_j m_j^*\|$, for any o.b. $\{m_j\}_j$ of $\Cal M$ relative to $\Cal N$ (see Theorem 1.1.6 in [P93b]). 

Thus, under this dimension condition,  
one has Ind$(\Cal E)=\text{\rm Ind}(\Cal E \otimes id)$  for the inclusion $(\Cal N\subset \Cal M)\overline{\otimes} \Cal B(\Cal H)$. 
Note that this stability of the index  implies the ``stability'' of the inequality $\Cal E(x)\geq \lambda x$, $\forall x\in \Cal M_+$, where $\lambda=\lambda(\Cal E)$, 
in the sense that more than being positive, the map $(\Cal E - \lambda id_{\Cal M}): \Cal M\rightarrow \Cal M$ is completely positive.  

If $Z=\sum_i m_i m_i^*$ is bounded, then $\Cal E_1:\langle \Cal M, e\rangle \rightarrow \Cal M$ defined by $\Cal E_1(xey)=xyZ^{-1}$, 
defines a normal conditional expectation with Ind$(\Cal E_1)<\infty$. More precisely, by 
(Lemma 1.2.1 in [P93b]) one has: $\text{\rm Ind}(\Cal E_1)=\|\Cal E(Z)\|$; $\{m_ieZ^{1/2}\}_i$ is an o.b. for $\Cal E_1$;  $\|\sum_j m_j e Z e m_j^*\|= 
\|\Cal E(Z)\|$.  

An important feature of the probabilistic index Ind$(\Cal E)$ of an  expectation $\Cal E$ is that it behaves well to ``limit operations''.  
For instance, if $(\Cal N_n \subset \Cal M_n)$ is a sequence of $W^*$-inclusions that are embedded into $\Cal N \subset^{\Cal E} \Cal M$  such that 
$\Cal E(\Cal M_n)=\Cal N_n$, $\forall n$, and $\Cal M_n \nearrow \Cal M$, $\Cal N_n \nearrow \Cal N$, 
then $\lim_n \text{\rm Ind}(\Cal E_n)=\text{\rm Ind}(\Cal E)$ and if $f_n \in \Cal P(\Cal N)$, $f_n \nearrow 1$, 
then $\lim_n \text{\rm Ind}(\Cal E(f_n \cdot f_n))=\text{\rm Ind}(\Cal E)$ (see [PP84], [PP86]). 
Also, if $\Cal E_i:\Cal M\rightarrow \Cal N$ is a net of expectations with 
decreasing (finite) index, equivalently $\lambda(\Cal E_i)\nearrow \lambda_0>0$, 
then any Banach-limit $\Cal E$ of the $\Cal E_i$ will be an expectation of $\Cal M$ onto $\Cal N$ that satisfies $\Cal E(x)\geq \lambda_0 x$, $\forall x\in \Cal M_+$ 
(so it is automatically normal and faithful). 
Thus,  if one defines the {\it minimal index} Ind$_{min}(\Cal N\subset \Cal M)\in (0, \infty]$ 
of an arbitrary inclusion of von Neumann algebras as the infimum of Ind$(\Cal E)$, over all normal expectations $\Cal E$ of $\Cal M$ onto $\Cal N$ 
(with the convention that it is equal to $\infty$ if there exists none), then Ind$_{min}(\Cal N\subset \Cal M) < \infty$ implies there exists an expectation $\Cal E_0$ with Ind$(\Cal E_0)= 
\text{\rm Ind}_{min}(\Cal N \subset \Cal M)$. In other words, the minimal index ``is attained'', if finite. Note that if $\Cal E_i$ are $\lambda_i$-Markov, 
then $\Cal E_0$ follows $\lambda_0=\lim_i \lambda_i$ Markov. It is easy to see that if $Q\subset P$ are finite dimensional as above, with connected inclusion 
bipartite graph $\Lambda=\Lambda_{Q\subset P}$, then there exists a unique expectation $E:P \rightarrow Q$ with Ind$(E)=\text{\rm Ind}_{min}(Q\subset P)$, 
and it is exactly the expectation preserving the trace on $P$ implemented by the Perron-Frobenius eigenvector $\vec{t}$ of $\Lambda^t\Lambda$, 
thus being $\lambda=\|\Lambda\|^{-2}$ Markov (see e.g. Section 1.1.7  in [P92a]). In particular, $\|\Lambda\|^2 \leq \text{\rm Ind}(E)$ for any expectation $E:P \rightarrow Q.$ 

Another feature of the probabilistic index when considered for an expectation of a C$^*$-inclusion $C\subset^E B$, 
is that it allows taking the bidual $W^*$-inclusion $C^{**}\subset^{E^{**}} B^{**}$, which satisfies $\text{\rm Ind}(E^{**})=\text{\rm Ind}(E)<\infty$ 
(see Section 2.4 in [P92a] and Section 3.3 in this paper).  The [PP84]-index is also key in relating the relative Dixmier property  for W$^*$ and C$^*$-inclusions with the finiteness of the index ([P97], [P98]), as well as in the calculation of the relative entropy for finite dimensional $W^*$-inclusions $Q\subset P$ (see Sec. 6 in [PP84]) 
and for inductive limits of such inclusions (see [PP88]).  

Finally, let us mention that if the $W^*$-inclusion $\Cal N\subset^\Cal E \Cal M$ is {\it irreducible}, i.e., if one has trivial relative commutant $\Cal N'\cap \Cal M$, then $\Cal E$ 
is the unique expectation of $\Cal M$ onto $\Cal N$, so one can as well use the Jones notation for the index $[\Cal M: \Cal N]$, even when the factors $\Cal N, \Cal M$ are 
infinite (non-tracial).  For arbitrary $W^*$-inclusions of infinite factors with finite index expectation $\Cal N\subset^\Cal E\Cal M$ ([K86]), 
one usually considers $\Cal E$ to be the expectation of minimal index ([Hi88]),  determined uniquely by the fact that the (scalar) 
values $\Cal E(q)$ of $\Cal E$ on minimal projections $q\in \Cal N'\cap \Cal M$ are proportional to $[q\Cal Mq: \Cal Nq]^{1/2}$ ([Hi88]). 
As pointed out in ([K86]), in this case one can take an o.b. of just one element $\{m\} \subset \Cal M$, which thus satisfies $X=m\Cal E(m^*X)$, $\forall X\in \Cal M$ 
and $\Cal E(m^*m)=1$, $mm^*=\lambda^{-1}1$.

\vskip.05in

\noindent
{\bf 2.6. $\lambda$-Markov inclusion, the associated tower and enveloping algebra}. Following (3.3.1 in [J82]; 1.1.4 in [P92a]; 1.2 in [P93b]), 
a $W^*$-inclusion with expectation $\Cal N \subset^{\Cal E} \Cal M$ that has finite index and satisfies 
$\sum_j m_jm_j^*=\lambda^{-1}1\in \Bbb C1$ for some (thus any) o.b. is called a $\lambda$-{\it Markov inclusion}. 

If $\Cal N \subset^{\Cal E} \Cal M$ is $\lambda$-Markov and $\Cal M\subset^{e} \langle \Cal M, e\rangle$ is its basic construction, 
then $\Cal E_1: \Cal M_1 \rightarrow \Cal M$ defined by $\Cal E_1(xey)=\lambda xy$, $x, y \in \Cal M$, 
is a conditional expectation of $\Cal M_1=\langle \Cal M, e \rangle$ onto $\Cal M_0=\Cal M$ 
and $\Cal M_0\subset^{\Cal E_1} \Cal M_1$ is again $\lambda$-Markov, with o.b. $\{\lambda^{-1/2}m_ie\}_i$. 

Thus, like in the II$_1$ factor case, one can iterate this construction and obtain a whole {\it tower of $\lambda$-Markov inclusions} $\Cal M_{n-1}\subset^{\Cal E_n}_{e_{n-1}} \Cal M_{n}$, 
$n\geq 1$, where we have put $e_0=e$. Moreover,   
the composition of expectations $\Cal E_1 \circ ....\circ \Cal E_n \circ ...  $ implements a 
trace $\tau$ on $Alg(\{e_n\}_{n\geq 0}$  with $\tau(e_n)=\lambda$, $\forall n$,  and the sequence of projections $e_0, e_1, ....$ satisfies properties 2.2 $(a')-(c')$ with respect to this trace and $\lambda$. Thus, $\lambda^{-1}$ lies in the Jones spectrum  and $(\{e_n\}_{n\geq 0}, \tau)$ gives a $\lambda$-sequence of projections. 

Moreover, the inductive limit of the 
Jones tower $\Cal M_n$ associated with a $\lambda$-Markov $W^*$-inclusion $\Cal N \subset^{\Cal E} \Cal M$ gives rise to a canonical 
{\it enveloping} von Neumann algebra, denoted $\Cal M_\infty$. 

\proclaim{2.6.1. Lemma} Let $\Cal N\subset^{\Cal E} \Cal M\subset_{e_0}^{\Cal E_1} \Cal M_1 \subset ...$ be a $\lambda$-Markov tower 
of $W^*$-inclusions. Let $\phi$ be a normal faithful state on $\Cal N$ and still denote by $\phi$ the state  on $\cup_n \Cal M_n$ which on $\Cal M_n$ is 
defined by $\phi(X)=\phi \circ \Cal E\circ \Cal E_1 \circ ... \circ \Cal E_n (X)$. Denote $(\Cal M_\infty, \Cal H_\phi)$ the GNS completion of 
$(\cup_n \Cal M_n, \phi)$. Then we have: 
\vskip.05in 

$(i)$ The spatial isomorphism class of $(\Cal M_\infty, \Cal H_\phi)$ doesn't depend on $\phi$. 

$(ii)$ The tower $\Cal M_n$ is naturally embedded into $\Cal M_\infty$ and there exist unique 
$\phi$-preserving conditional expectations $\tilde{\Cal E}_n: \Cal M_\infty \rightarrow \Cal M_{n-1}$ satisfying  
$\tilde{\Cal E}_{n|\Cal M_m}=\Cal E_n \circ ... \circ \Cal E_m$ for all $m \geq n\geq 0$, where $\Cal E_0=\Cal E$. 

\vskip.05in

Moreover, if $\text{\rm Tr}$ is an n.s.f. $($normal semifinite faithful$)$ trace on $\Cal M$ 
that's $\Cal E$-invariant, then $\text{\rm Tr}:=\text{\rm Tr}\circ \tilde{\Cal E}_0$ defines an n.s.f.  
trace on $\Cal M_\infty$ and  all $\tilde{\Cal E}_n$ leave invariant this trace. Also, if $\lambda\neq 1$ then we have: 
\vskip.05in 

$(a)$ If $\Cal M$ has no finite direct summand, then $\Cal M_\infty$ is of type $\text{\rm II}_\infty$. 

$(b)$ If $\text{\rm Tr}(1)< \infty$, then  $\Cal M_\infty$ is $\text{\rm II}_1$. If in addition $\Cal N \subset \Cal M$ are finite dimensional with irreducible 
inclusion graph, then $\Cal M_\infty$ is a type $\text{\rm II}_1$ factor.  
\endproclaim
\noindent
{\it Proof}. Parts $(i), (ii)$ are standard applications of Takesaki's classic 
results in [T72]. The fact that $\text{\rm Tr}$ n.s.f. implies $\text{\rm Tr} \circ \Cal E\circ \Cal E_1 \circ ... \circ \Cal N_n$ is 
n.s.f. on $\Cal M_n$ follows easily from the fact that Ind$(\Cal N \subset \Cal M_n) = \lambda^{-n-1}<\infty$.  
The last part is trivial and we leave it as an exercise. 
\hfill $\square$

\vskip.05in
Note that if $\Cal M$ is a factor, then any $\Cal N \subset^{\Cal E} \Cal M$ with finite index is Markov. In particular, if $N\subset M$ are II$_1$ factors 
with $[M:N]<\infty$ and $\lambda=[M:N]^{-1}$, then $N\subset M$ is $\lambda$-Markov. 

If $\Cal N\subset^{\Cal E}\Cal M$ is a $W^*$-inclusion with Ind$(\Cal E)=\lambda^{-1}<\infty$, then any other normal expectation of $\Cal M$ onto $\Cal N$ 
is of the form $\Cal E_A=\Cal E(A^{1/2} \cdot A^{1/2})$, for some $A\in (\Cal N'\cap \Cal M)_+$ with $\Cal E(A)=1$. Thus, such $A$ satisfies $1=\Cal E(A)\geq \lambda A$, 
so $A\leq \lambda^{-1} 1$. This same argument applied to $\Cal E_A$ shows that if $\lambda_A=\lambda(\Cal E_A)>0$ then $A\geq \lambda_A 1$. Note that $\{m_j\}_j$ is o.b. for 
$\Cal E$ iff $\{m_jA^{-1/2}\}_j$ is o.b. for $\Cal E_A$. 

A $W^*$-inclusion $\Cal N\subset^{\Cal E}\Cal M$  may have finite index without being Markov, 
yet  for some $A$ the expectation $\Cal E_A$ become Markov. For this, one needs $A\in (\Cal N'\cap \Cal M)_+$ to satisfy $\sum_j m_j A^{-1}m_j^* \in \Bbb C1$. 
For instance, by (Theorem 3.2 in [J82]), if $\oplus_{i\in I}\Bbb M_{n_i}(\Bbb C) = Q\subset^E P=\oplus_{j\in J} \Bbb M_{m_j}(\Bbb C)$ 
is a finite dimensional $W^*$-inclusion with $E$ preserving a normal faithful trace state $\tau$ on $P$, given by the weight vector 
$\vec{t}=(t_j)_j$, then $E$ is $\lambda$-Markov iff $\vec{t}$ is a Perron-Frobenius eigenvector for 
$\Lambda^t \Lambda$, corresponding to $\lambda^{-1}=\|\Lambda^t\Lambda\|=\|\Lambda\|^2$, 
where $\Lambda=\Lambda_{Q\subset P}$ is the inclusion 
bipartite graph of $Q\subset P$, viewed as $I \times J$ matrix, i.e., $\Lambda^t\Lambda \vec{t}=\lambda^1\vec{t}$.

\vskip.05in

\noindent
{\bf 2.7. Atomic  $W^*$-inclusions}. Let $\Cal N \subset^{\Cal E} \Cal M$ be an inclusion of atomic von Neumann algebras with Ind$(\Cal E)=\lambda^{-1}<\infty$.  
Thus, both $\Cal N, \Cal M$ are  direct sums of type I factors, $\Cal N=\oplus_{i\in I} \Cal B(\Cal K_i)$, $\Cal M=\oplus_{j\in J} \Cal B(\Cal H_j)$. 
It is trivial to see that if $b_{ij}$ denotes the multiplicity of $\Cal B(\Cal K_i)$ in $\Cal B(\Cal H_j)$ (i.e., $b_{ij}^2 = \text{\rm dim}(\Cal B(\Cal K_i)'\cap \Cal B(\Cal H_j)$), 
then $b_{ij} \leq \lambda^{-1}$ (if $b_{ij}\leq \text{\rm dim}(\Cal K_i)$ then one actually has $b_{ij}\leq \lambda^{-1/2}$, see Section 6 in [PP84] or page 200 in [P92a]). 
The finiteness of the index also implies that for each $j\in J$ 
(resp. $i\in I$) the number of $i\in I$ (resp. $j\in J$) with $b_{ij}\neq 0$ is finite. 

Thus, such an inclusion $\Cal N \subset \Cal M$ is described by a bipartite graph (diagram) 
$\Lambda_{\Cal N\subset \Cal M}=(b_{ij})_{i\in I, j\in J}$, where the number $b_{ij}$ of edges between the vertices $i$ and $j$ is equal to the multiplicity of 
$\Cal B(\Cal K_i)$ in $\Cal B(\Cal H_j)$.  We alternatively view $\Lambda$ as an $I \times J$ matrix with integer non-negative entries $b_{ij}$.  One trivially has 
that $\Cal Z(\Cal N) \cap \Cal Z(\Cal M)=\Bbb C1$ iff $\Lambda$ is connected (equivalently, $\Lambda$ is irreducible as a matrix). 

We are particularly interested  in the case $\Cal N, \Cal M$ are of type I$_\infty$ 
and $\Cal N \subset^{\Cal E} \Cal M$ is $\lambda$-Markov, where $\lambda^{-1}=\text{\rm Ind}(\Cal E)$. 
Assume this is the case and let $\Cal M \subset_e^{\Cal E_1} \Cal M_1=\langle \Cal M, e\rangle$ 
be the associated basic construction, with $e$ denoting the Jones projection. Then $\Cal M_1$ is an amplification of $\Cal N$, so it is again an atomic 
von Neumann algebra and by (3.3.2 in [J82]) there is a natural isomorphism $\Cal Z(\Cal N) \ni z \mapsto z' \in \Cal Z(\Cal M_1)$, where $z'$ is the 
unique element in $\Cal Z(\Cal M_1)$ with $z'e=ze$. Moreover, if one identifies the set labeling the direct summands of $\Cal M_1$ with $I$, via this identification, then 
the $J \times I$ bipartite graph $\Lambda_{\Cal M \subset \Cal M_1}$ is given by $(\Lambda_{\Cal N \subset \Cal M})^t=(b'_{ji})_{j\in J, i\in I}$, where $b'_{ji}=b_{ij}$. 

An important case is when the atomic $W^*$-inclusion $\Cal N\subset^{\Cal E}\Cal M$ is both Markov and 
$\Cal E$ leaves invariant some normal semifinite faithful (abbreviated hereafter n.s.f.) trace Tr on $\Cal M$. A $W^*$-inclusion $\Cal N\subset^{\Cal E}\Cal M$ with the property that 
$\Cal M$ admits a $\Cal E$-invariant n.s.f. trace, is called {\it Tracial}. If $\Cal N\subset^{\Cal E} \Cal M$ is Tracial and $\lambda$-Markov, 
then an $\Cal E$-invariant n.s.f. trace Tr on $\Cal M$ is called a $\lambda$-{\it Markov Trace}.

The proof of (Theorem 3.3.2 in [J82]) can be easily adapted to general atomic inclusions to  completely characterize when a Tracial $\Cal E$ is Markov:

\proclaim{2.7.1. Lemma}  Let $\oplus_{i\in I} \Cal B(\Cal K_i)=\Cal N \subset^{\Cal E} \Cal M=\oplus_{j\in J} \Cal B(\Cal H_j)$ be an atomic $W^*$-inclusion 
of finite index with inclusion graph $\Lambda=\Lambda_{\Cal N\subset \Cal M}=(b_{ij})_{i\in I, j\in J}$. Then we have: 

$1^\circ$ $\|\Lambda\|^2 \leq \text{\rm Ind}(\Cal E)$.  

$2^\circ$ Assume $\Cal E$ preserves a n.s.f. trace $\text{\rm Tr}$ on $\Cal M$, with 
$\vec{t}=(t_j)_j$ its weight vector,  i.e., $t_j$ is the trace of a minimal projection in $\Cal B(\Cal H_j)$. Denote $\lambda=\text{\rm Ind}(\Cal E)^{-1}$.   
Then $\tilde{\text{\rm Tr}}(xey)=\lambda \text{\rm Tr}(xy)$, $x, y\in \Cal M$, defines a n.s.f. trace on $\Cal M_1$ and 
the following conditions are equivalent: 
\vskip.05in

$(i)$ $\Cal N\subset^{\Cal E}\Cal M$ is $\lambda$-Markov; 

$(ii)$  $\tilde{\text{\rm Tr}}_{|\Cal M}=\text{\rm Tr}$; 

$(iii)$  $\Lambda^t\Lambda(\vec{t})=\lambda^{-1}\vec{t}$. 
\vskip.05in
 
$3^\circ$ If the equivalent conditions $(i)-(iii)$ in $2^\circ$ above are satisfied and $\Lambda$ is finite, then $\lambda^{-1}=\|\Lambda\|^2$. 
\endproclaim
\noindent
{\it Proof}. $1^\circ$ Let $f_i \in \Cal P(\Cal N)$ be an increasing net of finite rank projections such that $f_i \nearrow 1$ and denote $\Lambda_i$ 
the inclusion graph of the finite dimensional inclusions $f_i\Cal Nf_i \subset f_i\Cal M f_i$. 
We clearly have 
$\lim_{i}\|\Lambda_i\|=\|\Lambda\|$ and $\text{Ind}\Cal E(f_i \cdot f_i) \leq \text{Ind}\Cal E$, $\forall i$. We already pointed out in Section 2.5 that 
in the finite dimensional case ones has Ind$\Cal E(f_i \cdot f_i)\geq \|\Lambda_i\|^2$. Thus, 
$\|\Lambda\|^2=\lim_i\|\Lambda_i\|^2 \leq \limsup_i \text{Ind}\Cal E(f_i \cdot f_i) \leq \text{Ind}\Cal E$. 

$2^\circ$ The fact that $\tilde{\text{\rm Tr}}$ is a trace follows from the fact 
that $(x_1ey_1)(x_2ey_2)=x_1\Cal E(y_1x_2)ey_2$, $(x_2ey_2)(x_1ey_1)=x_2\Cal E(y_2x_1)ey_1$ so applying $\tilde{\text{\rm Tr}}$ gives 
$$
\tilde{\text{\rm Tr}}((x_1ey_1)(x_2ey_2))=\lambda\text{\rm Tr}(x_1\Cal E(y_1x_2)y_2)=\lambda\text{\rm Tr}(\Cal E(y_1x_2)y_2x_1)
$$
$$
= \lambda\text{\rm Tr}(\Cal E(y_1x_2)\Cal E(y_2x_1))=\lambda\text{\rm Tr}(x_2\Cal E(y_2x_1)y_1)= \tilde{\text{\rm Tr}}((x_2ey_2)(x_1ey_1)).
$$

$(i) \Rightarrow (ii)$ Let $\{m_k\}_k\subset \Cal M$ be an o.b. with respect to $\Cal E$. 
If $x\in \Cal M$ is finite rank then by applying  $\tilde{\text{\rm Tr}}$ to $x=x\sum_k m_k e m_k^*$ we get 
$$
\tilde{\text{\rm Tr}}(x)=\sum_k \tilde{\text{\rm Tr}}(xm_k e m_k^*)=\lambda \text{\rm Tr}(x\sum_k m_k m_k^*)=\lambda^{-1}\lambda \text{\rm Tr}(x) =\text{\rm Tr}(x). 
$$

$(ii) \Rightarrow (iii)$ The equality $\tilde{\text{\rm Tr}}(x)= \text{\rm Tr}(x)$, $\forall x\in \Cal M$, implies $\tilde{\text{\rm Tr}}(ex)=\lambda\text{\rm Tr}(x)=
\lambda\tilde{\text{\rm Tr}}(x)$ and thus $\Cal E_1$ is the $\tilde{\text{\rm Tr}}$-preserving expectation of $\Cal M_1$ onto $\Cal M$.   

Denote $\vec{s}=\Lambda(\vec{t})$ and note that $s_i$ gives the trace Tr (thus also trace $\tilde{\text{\rm Tr}}$) of any minimal projection $f_i$ in $\Cal B(\Cal K_i)$. 
Then  $f_ie$ is a minimal projection in the $i$'th summand of $\Cal M_1$ and from the above it has trace $\tilde{\text{\rm Tr}}(f_ie)=\lambda s_i$  
One thus gets $\Lambda\Lambda^t (\lambda \vec{s})=\vec{s}$ and $\Lambda^t(\lambda \vec{s})=\vec{t}$. 
Thus $\Lambda^t\Lambda(\vec{t})=\Lambda^t(\vec{s})=\lambda^{-1}\vec{t}$. 

$(iii) \Rightarrow (i)$ The fact that $\vec{s}=\Lambda(\vec{t})$ together with $(iii)$ shows that the trace $\text{\rm Tr}_1$ on $\Cal M_1$ 
given by the vector $\lambda\vec{s}$ has the property that restricted to $\Cal M, \Cal N$ coincides with $\text{\rm Tr}$ and that $\text{\rm Tr}_1(ex)=\lambda \text{\rm Tr}(x)$, 
for all finite rank $x\in \Cal M$. This in turn implies that the unique $\text{\rm Tr}_1$-preserving expectation $\Cal E_1'$ of $\Cal M_1$ onto $\Cal M$ 
satisfies $\Cal E_1'(e)=\lambda 1$ and thus coincides with $\Cal E_1$. 

$3^\circ$ Since $\Lambda^t\Lambda(\vec{t})=\lambda^{-1}\vec{t}$ and $\vec{t}$ has positive entries, if $\Lambda$ is finite then one necessarily has 
$\lambda^{-1}=\|\Lambda^t\Lambda\|$, by the Perron-Frobenius theorem (see e.g., [G60]). 
\hfill $\square$

\vskip.05in

\noindent
{\bf 2.7.2. Definition.}  An  atomic properly infinite $W^*$-inclusion $(\Cal N\subset^{\Cal E} \Cal M, \text{\rm Tr})$ with an $\Cal E$-invariant semifinite trace Tr 
that satisfies  the equivalent conditions $(i)-(iii)$ in 2.7.1.2$^\circ$ above is called a {\it Tracial $\lambda$-Markov} atomic inclusion. 

\vskip.05in 

The above lemma shows that such an object is completely determined by its inclusion bipartite graph $\Lambda=\Lambda_{\Cal N \subset \Cal M}$ 
and the trace vector $\vec{t}$. Indeed, by $2.7.1.2^\circ$, given any pair $(\Lambda, \vec{t})$ 
where $\Lambda=(b_{ij})_{i\in I, J\in J}$ is a bipartite graph and $\vec{t}=(t_j)_{j\in J}$ has positive entries 
and satisfies $\Lambda^t\Lambda (\vec{t})=\lambda^{-1}\vec{t}$ there exists a unique inclusion of type I$_\infty$ atomic von Neumann 
algebras with a n.s.f. trace Tr and Tr-preserving expectation $\Cal E$ 
such that $(\Cal N\subset^\Cal E \Cal M, \text{\rm Tr})$ is a Tracial $\lambda$-Markov with $\Lambda_{\Cal N\subset \Cal M}=\Lambda$ and 
Tr given by $\vec{t}$. 

We call such a pair $(\Lambda, \vec{t})$ a {\it Markov weighted bipartite graph}. If in addition 
we fix an even vertex $j_0\in J$ and renormalize $\vec{t}$ so that $t_{j_0}=1$, then $(\Lambda, j_0, \vec{t})$ is called a {\it pointed} Markov weighted graph. 
When specifying the scalar $\lambda$ with $\Lambda^t\Lambda(\vec{t})=\lambda^{-1} \vec{t}$, we call it a 
(pointed) $\lambda$-Markov weighted graph.  

We'll show in 2.7.5 below that  such an object $(\Lambda, j_0, \vec{t})$ is also equivalent to a sequence of inclusions of finite dimensional 
C$^*$-algebras $\{A_n\}_{n\geq 0}$  with a trace state $\tau$ and a representation of the Jones $\lambda$-projections $\{e_n\}_{n\geq 1}$,  
satisfying certain properties. It turns out that this set of characterizing properties (axioms)   
is surprisingly minimal (cf. Remark 1.4.3  and Theorem 1.5 in [P00]).

\proclaim{2.7.3. Proposition} Let $\Bbb C=A_0 \subset A_1 \subset A_2 \subset ...$ be a sequence of inclusions of finite dimensional $C^*$-algebras,  
with a faithful trace $\tau$ on $\cup_n A_n$ and a representation of the Jones $\lambda$-sequence of projections $\{e_n\}_{n\geq 1}\subset \cup_n A_n$, such that: 

\vskip.05in

$(i)$ $e_n \in A_{n-1}'\cap A_{n+1}$, $\forall n\geq 1$; 

$(ii)$ $e_nxe_n=E_{A_{n-1}}(x)e_n$, $\forall x\in A_n$, $n\geq 1$. 

\vskip.05in

Then we have: 

\vskip.05in 

$(a)$ For each $n\geq 1$ and $z\in \Cal Z(A_{n-1})$ there exists a unique $z'\in \Cal Z(A_{n+1})$ such 
that $ze_n=z'e_n$, and the resulting map $z\mapsto z'$  implements an embedding $\Cal Z(A_{n-1}) \hookrightarrow \Cal Z(A_{n+1})$. 

$(b)$ If $J_n$ $($respectively $I_n)$ labels the set of simple summands in $A_{2n}$ $($resp. $A_{2n+1})$ and we identify 
$J_n$ $($resp. $I_n)$ with a subset of $J_{n+1}$ $($resp. $I_{n+1})$ and let $J=\cup_n J_n$, $I=\cup_n I_n$, 
then there exists a unique pointed $J\times I$ bipartite graph 
$(\Gamma, \{j_0\})$, where $\{j_0\}=J_0\subset J$, such that $\Lambda_{A_{2n}\subset A_{2n+1}}=_{K_n}\Gamma$ and 
$\Lambda_{A_{2n+1}\subset A_{2n+2}}=_{L_n}\Gamma^t$, $n \geq 0$. 

$(c)$ There exist unique vectors $\vec{t}=(t_j)_{j\in J}$, $\vec{s}=(s_i)_{i\in I}$ such that $t_{j_0}=1$, $\Gamma^t(\vec{t})=\vec{s}$, $\Gamma\Gamma^t(\vec{t})=\lambda^{-1}\vec{t}$, 
and $\lambda^n t_j$ gives the trace of a minimal projection in the $j$'th summand of $A_{2n}$, while $\lambda^ns_i$ gives the trace  of a 
minimal projection in the $i$'th summand of $A_{2n+1}$. 
\endproclaim
\noindent
{\it Proof}.  Let $B_{0n}=A_n, n\geq 0$ and define $B_{11}=\Bbb C=B_{12}$, $B_{1n}=Alg\{1, e_k\mid 2\leq k \leq n-1\} $ for $n\geq 3$. 
Then $(B_{ij})_{j\geq i; i=0,1}$ is clearly a generalized $\lambda$ sequence of commuting squares, in the sense 
of (Definition 1.3 in [P00]). Thus, by (Theorem 1.5 in [P00]),  it follows that Ind$(E^{A_n}_{A_{n-1}})\leq \lambda^{-1}$ for all $n\geq 1$, 
and that there exists a Jones $\lambda$-tower of factors $M_{-1} \subset M_0\subset M_1 \subset_{e_1} M_2  \subset  ...$, with $A_n\subset M_n$ 
satisfying $E_{A_n}E_{M_{n-1}}=E_{A_{n-1}}$, $n \geq 1$ (the commuting square relation, see Section 2.8 below).  

But then the proof of the existence of a unique graph and weight vector in (Proposition 2.1, Corollary 2.2 in [P89]) 
for a sequence of inclusions of finite dimensional C$^*$-algebras $\Bbb C=A_0 \subset A_1 \subset A_2 \subset ...$ with Jones $\lambda$-projections 
works exactly the same in this more general case.  
\hfill $\square$

\vskip.05in 
\noindent
{\bf 2.7.4. Definition}. Following (1.4.3 in [P00]), a sequence of  inclusions of finite dimensional C$^*$-algebras with a faithful trace and a representation of the 
Jones $\lambda$-sequence of projections $(\{A_n\}_{n\geq 0}, \{e_n\}_{n\geq 1}, \tau)$ satisfying conditions $(i), (ii)$ in 2.7.3 above  
is called a $\lambda$-{\it sequence of inclusions}, with $(\Gamma, j_0, \vec{t})$ in $(c)$ above being its associated pointed weighted graph. 

\vskip.05in 

\proclaim{2.7.5. Proposition} Let $(\Lambda, j_0, \vec{t})$ be a pointed bipartite graph with a $\lambda$-Markov weight vector. Let 
$(\Cal N\subset^\Cal E \Cal M, \text{\rm Tr})$ be the associated Tracial $\lambda$-Markov atomic $W^*$-inclusion and 
$\Cal N \subset \Cal M \subset^{\Cal E_1}_{e_0} \Cal M_1 \subset^{\Cal E_2}_{e_1} \Cal M_2 ...$ 
its Jones tower, with the semifinite trace $\text{\rm Tr}$ on $\cup_n \Cal M_n$ that's invariant to all $\Cal E_n$.  
Let $p=p_{j_0}$ be a minimal projection in $\Cal B(\Cal H_{j_0})$. Then the sequence of inclusions 
$$
(A_0 \subset A_1 \subset_{e'_1} A_2 \subset ...)\simeq (p\Cal M p \subset p\Cal M_1p\subset_{e_1p} p\Cal M_2p \subset ...)
$$ 
with the trace state $\tau=\text{\rm Tr}(p \cdot p)$ and Jones projections $e_i'=e_ip$, is a $\lambda$-sequence of inclusions, with its graph $\Gamma$ given by $\Lambda^t$. 
\endproclaim 
\noindent
{\it Proof}. This is trivial by the properties of Tracial $\lambda$-Markov inclusion and its Jones tower. 
\hfill $\square$

\vskip.05in 

\noindent
{\bf 2.8. Commuting square embeddings}. If $\Cal Q \subset^{\Cal F} \Cal P$, $\Cal N \subset^{\Cal E} \Cal M$ are $W^*$-inclusions with expectation then 
a {\it commuting square} ({\it c.sq.}) {\it embedding} 
(or simply an {\it embedding}) of $\Cal Q \subset^{\Cal F} \Cal P$ into $\Cal N \subset^{\Cal E} \Cal M$ is a von Neumann algebra inclusion  $\Cal P \subset \Cal M$ 
so that $\Cal Q \subset \Cal N$, $\Cal E(\Cal P)=\Cal Q$, $\Cal E_{|\Cal P}=\Cal F$. Note that this trivially implies $\text{\rm Ind}(\Cal F)\leq \text{\rm Ind}(\Cal E)$ 
and that if these indices are finite, then any o.b. $\{\xi_i\}_i \subset \Cal P$ for $\Cal F$ is 
an orthonormal system for $\Cal E$, which can be completed to an o.b. $\{\eta_j\}_j\subset \Cal M$ for $\Cal E$, thus $\sum_i \xi_i\xi_i^*  \leq \sum_j \eta_j\eta_j^*$

This property was first considered in (page 29 and 1.2.2 in [P82]; cf. also [P81]) to study the ``interaction'' between subalgebras of a tracial von Neumann algebra $M$, 
calculate relative commutants and normalizers of subalgebras. If $P, Q\subset M$  are von Neumann subalgebras, 
and $E_P, E_Q$ denote as usual the trace preserving expectations onto them, then the commuting square relation 
amounts to $E_P E_Q= E_Q E_P=E_{P\cap Q}$. This is equivalent to $E_P(Q)\subset P$ and also 
to $e_P e_Q$ being a  projection in $\Cal B(L^2M)$. In case $M$ is the von Neumann algebra $L\Gamma$ of a discrete group $\Gamma$, then any two subgroups 
$G, H\subset \Gamma$ give rise to subalgebras $P=LG, Q=LH$ satisfying this relation. So one can view the c.sq.  relation for subalgebras of a II$_1$ factor as 
a natural ``lattice-like'' condition. 

Commuting square embeddings are the  natural ``morphisms'' between $W^*$-inclu- sions with expectation. 
All such morphisms considered here will be taken between inclusions of ``same index'', a property 
that we call  ``non-degeneracy'' (following 1.1.5 in [P92a]).  More precisely, the commuting square embedding of 
$\Cal Q \subset^{\Cal F} \Cal P$ into $\Cal N \subset^{\Cal E} \Cal M$
is {\it non-degenerate} if  one has that $\text{\rm sp} \Cal P\Cal N$ is weakly dense in $\Cal M$ (N.B.: in all cases of interest for us $\Cal  P$ is finitely generated right $\Cal Q$-module, 
where the condition becomes $\text{\rm sp}\Cal P\Cal N=\Cal M$). 

Note that if Ind$(\Cal E)<\infty$, or if all algebras involved are tracial and the expectations involved 
are trace preserving (with respect to the trace on the largest algebra $\Cal M$), 
then this is equivalent to saying that any o.b. for $\Cal F$ is an o.b. for $\Cal E$. So one necessarily have Ind$(\Cal E)=\text{\rm Ind}(\Cal F)$.  
Moreover, if $\Cal M \subset \Cal M_1=\langle \Cal M, e_{\Cal N}\rangle$ is the basic construction for $\Cal E$ then sp$\{xe_{\Cal N} y \mid x, y \in \Cal P\}$ 
is a $^*$-subalgebra of $\Cal M_1$ and its weak closure contains $\Cal P$ and identifies naturally with the basic construction 
$\Cal Q \subset^{\Cal F}\Cal P \subset \Cal P_1:=\langle \Cal P, e_{\Cal Q}\rangle$. 

If the above c.sq. is non-degenerate, it follows trivially that $\Cal Q \subset^{\Cal F} \Cal P$ is $\lambda$-Markov iff $\Cal N \subset^{\Cal E} \Cal M$ is 
$\lambda$-Markov, with $\lambda=\lambda(\Cal E)=\lambda(\Cal F)$. 
Also, if one has a c.sq. embedding $(\Cal Q\subset^{\Cal F} \Cal P) 
\subset(\Cal N\subset^{\Cal E}\Cal M)$ and both $\Cal Q \subset \Cal P$ and $\Cal N \subset \Cal M$ 
are $\lambda$-Markov (same $\lambda$) then the commuting square follows nondegenerate. 

If these conditions are satisfied, then we say that 
$(\Cal Q\subset^{\Cal F} \Cal P) 
\subset(\Cal N\subset^{\Cal E}\Cal M)$ is a $\lambda$-{\it Markov commuting square}. 
By the above observations  about basic construction for non-degenerate commuting squares, such c. sq. gives rise 
to a whole tower of $\lambda$-Markov c.sq. embeddings
$$
\CD
\Cal N\ @.\subset^{\Cal E}\ @.\Cal M\ @.\subset^{\Cal E_1}_{e_0}\ @.\Cal M_1 \subset^{\Cal E_2}_{e_1} ...   \\
\noalign{\vskip-6pt}
\cup\ @.\ @.\cup\ @.\ @.\cup  \\
\noalign{\vskip-6pt}
\Cal Q\ @.\subset^{\Cal F}\ @.\Cal P\ @.\subset^{\Cal F_1}_{e_0}\  @.\Cal P_1  \subset^{\Cal F_2}_{e_1} ....  
\endCD
$$
where $e_i$ implements both $\Cal E_i$ and $\Cal F_i$, $\Cal E_{i+1}(e_i)=\lambda$, $\Cal M_{i+1}=\text{\rm sp}\Cal M_ie_i\Cal M_i$ 
and $\Cal P_{i+1}=\text{\rm sp}\Cal P_i e_i\Cal P_i$.

Finally, let us recall from (Proposition 1.1.6 of [P92a]) that if $M$ is tracial and $N, P\subset M$ are von Neumann subalgebras such that 
$E_PE_N=E_NE_P=E_Q$, where $Q=N \cap P$  (but no finite index assumption on any of the expectations) 
then the resulting commuting square embedding $(Q\subset P) \subset (N\subset M)$ is nondegenerate iff the embedding $(Q\subset N)\subset (P\subset M)$ is nondegenerate. 

From a remark above, if such a commuting square is non-degenerate then one can take its basic construction both ``horizontally'' 
$(Q\subset P \subset \langle P, Q\rangle) \subset (N\subset M \subset \langle M, N \rangle)$ and ``vertically'' 
$(Q\subset N\subset \langle N, Q\rangle) \subset (P \subset M \subset \langle M, P\rangle)$, with the canonical trace Tr$_{|\langle M, N\rangle}$ restricted to $\langle P, Q\rangle$ 
equal to Tr$_{|\langle P, Q\rangle}$, resp.  Tr$_{|\langle M, P\rangle}$ restricted to $\langle N, Q\rangle$ 
equal to Tr$_{|\langle N, Q\rangle}$. Moreover,  
one has a canonical Tr-preserving expectation $E^{\langle M, N\rangle}_{\langle P, Q\rangle}$ (resp. $E^{\langle M, P\rangle}_{\langle N, Q\rangle}$) between the corresponding extension 
algebras and its restriction  to $M$ is equal to $E^M_P$ (resp. $E^M_N$). In addition, 
the resulting embeddings $(N\subset M) \subset (\langle N, Q\rangle \subset \langle M, P \rangle)$ 
and $(P\subset M) \subset (\langle P, Q\rangle \subset \langle M, N \rangle)$ are non-degenerate. Also, if the two rows (resp. columns) of the initial 
commuting square are Markov, then so is the row (resp. column) basic construction inclusion. 

\vskip.05in 

\noindent
{\bf 2.9. Subfactors of $R$ from Markov commuting squares}. An important class of Markov commuting squares are the ones with all the algebras involved 
finite dimensional  and all expectations trace preserving. We explain here how such an object produces an extremal hyperfinite subfactor 
as an inductive limit of the associated Jones tower (cf. [W88], [GHJ88]; this construction was found independently in early 1984 by Pimsner-Popa). 
 
Let us consider a finite dimensional c.sq. embedding $(P_{00}\subset P_{01}) 
\subset (P_{10}\subset P_{11})$, with a (faithful) trace state $\tau$ on the largest algebra $P_{11}$ and all expectations 
involved being $\tau$-preserving. The $\lambda$-Markov condition means 
that the c.sq. is non-degenerate and the row inclusions are $\lambda$-Markov. In addition we assume this c.sq. is so that both row 
inclusion graphs $\Lambda_{P_{00}\subset P_{01}}$ and $\Lambda_{P_{10}\subset P_{11}}$ are irreducible. We call such an object 
a  {\it Markov cell}. 

As shown in 2.8 and 2.6.1 $(b)$, such an object $(P_{00}\subset P_{01})\subset (P_{10}\subset P_{11})$ gives rise to a 
tower of $\lambda$-Markov c.sq., with row enveloping II$_1$ factors $P_{1\infty}$, $P_{0\infty}$ 
$$
\CD
P_{10}\ @.\subset\ @.P_{11}\ @.\subset\  @.P_{12}\ ...@.\  \nearrow @.P_{1,\infty}\\
\noalign{\vskip-6pt}
\cup\ @.\ @.\cup\ @.\  @.\cup\  @.\  @.\cup\ \\
\noalign{\vskip-6pt}
P_{00}\ @.\subset\ @.P_{01}\ @.\subset\  @.P_{02}\ ....@.\  \nearrow @.P_{0,\infty}
\endCD
$$
Moreover, due to the non-degeneracy at each step, the commuting square 
$$
\CD
P_{10}\ @.\subset\ @.P_{1\infty}\\
\noalign{\vskip-6pt}
\cup\ @.\ @.\cup\ \\
\noalign{\vskip-6pt}
P_{00}\ @.\subset\ @.P_{0\infty}
\endCD
$$
follows non-degenerate. This implies 
$\text{\rm Ind}(E^{P_{10}}_{P_{00}})=\text{\rm Ind}(E^{P_{1\infty}}_{P_{0\infty}})=[P_{1\infty}: P_{0\infty}]$. 
In addition, since the algebras involved are II$_1$ factors, the inclusion $P_{0\infty}\subset P_{1\infty}$ is $\lambda_{01}$-Markov, 
where $\lambda_{01}=[P_{1\infty}:P_{0\infty}]^{-1}$. So $P_{00}\subset P_{10}$ follows $\lambda_{01}$-Markov as well, and so do all 
the vertical inclusions $P_{0n}\subset P_{1n}$. 

Thus, $[P_{1\infty}:P_{0\infty}]=\lambda_{01}^{-1}=\|\Lambda_{P_{0n}\subset P_{1n}}\|^2$, $\forall n\geq 0$.  
Moreover, by the formulas for relative entropy in [PP88], it follows that $H(P_{1\infty}| P_{0\infty})=\lim_n H(P_{1n}| P_{0n}) = -\ln \lambda_{01}$, 
and hence $H(P_{1\infty}| P_{0\infty})=\ln [P_{1\infty}:P_{0\infty}]$. By (Corollary 4.5 in [PP84]),  this implies $P_{0\infty} \subset P_{1\infty}$ is extremal.

In conclusion, a Markov cell 
$$
\CD
P_{10}\ @.\subset\ @.P_{11}\\
\noalign{\vskip-6pt}
\cup\ @.\ @.\cup\ \\
\noalign{\vskip-6pt}
P_{00}\ @.\subset\ @.P_{01}
\endCD
$$
which by the definition is only assumed Markov ``horizontally'', follows Markov ``vertically'' as well, 
with the inclusion graphs satisfying $\Lambda_{P_{10}\subset P_{11}}\circ \Lambda_{P_{00}\subset P_{10}}
=\Lambda_{P_{01}\subset P_{11}} \circ \Lambda_{P_{00}\subset P_{01}}$ and $\|\Lambda_{P_{10}\subset P_{11}}\|=\|\Lambda_{P_{00}\subset P_{01}}\|$,  
$\|\Lambda_{P_{00}\subset P_{10}}\|=\|\Lambda_{P_{01}\subset P_{11}}\|$. Moreover, the  resulting inclusion of enveloping II$_1$ factors $P_{0\infty}\subset P_{1\infty}$ 
gives an extremal hyperfinite subfactor of index $\alpha=\|\Lambda_{P_{00}\subset P_{10}}\|^2$.  As we've seen in 2.2, if $4< \alpha < 3+2\sqrt{2}$, 
the subfactor follows irreducible. 

Given a finite bipartite graph $\Lambda$, {\it the commuting square problem for $\Lambda$} consists in constructing a Markov cell   
with the vertical inclusion graph $\Lambda_{P_{00}\subset P_{10}}$ equal to $\Lambda$.

\vskip.05in 

\noindent
{\bf 2.10.  The sets $\mycal C(R)$ and $\Bbb E^2$}. Following ([GHJ88]), 
we denote by $\Bbb E$ the set of norms of bipartite graphs (finite or infinite) and let $\Bbb E^2:=\{\alpha^2 \mid \alpha \in \Bbb E\}$. 
Equivalently, $\Bbb E^2$ is the set of square norms of matrices with non-negative integer entries.  
The set $\Bbb E$ was first described in ([H072], [CDG82]). A detailed account can be found in  (Appendix A in [GHJ88]). We 
recall some properties, using the notations therein.  

First of all, notice that $\Bbb E^2$ is a closed set consisting 
of an increasing sequence of  accumulation points, followed by the half line $[2+\sqrt{5}, \infty]$.  

One has $\Bbb E^2 \cap (0, 4] = \{4\cos^2(\pi/n)\mid n\geq 3\} \cup \{4\}$, the only bipartite graphs 
of square norm  $<4$ being the Coxeter graphs $A_n, D_n, E_6, E_7, E_8$. The bipartite graphs of square norm $4$ are $A_\infty, D_\infty, A_{-\infty, \infty}$, 
$D_n^{(1)}, E_6^{(1)}, E_7^{(1)}, E_8^{(1)}=E_9$.  

There is a gap right after the first accumulation point $4$, with $\|E_{10}\|^2 \approx 4.0265...$ being the first (smallest) element in $\Bbb E^2 \cap (4, \infty)$. 
Then one has an increasing sequence of accumulation points $c_n^2$, $n\geq 3$,  converging to $c^2_\infty = 2+\sqrt{5}$. 
Each $c_n^2$ for $n\geq 3$ is an accumulation point both from below and from above. The set of $\Bbb E_0^2$ of square norms of finite bipartite graphs 
is a dense subset of $\Bbb E^2$. 

Following ([J82], [J86]), given a II$_1$ factor $M$ one denotes $\mycal C(M)$  the set of indices of irreducible 
subfactors of $M$. We also let $\mycal E(M)$ be the set of indices of extremal subfactors of $M$. One obviously has $\mycal C(M^t)=\mycal C(M)$, $\mycal E(M^t)=\mycal E(M)$, $\forall t>0$.

Also, one has $\mycal C(M_1)\cdot \mycal C(M_2)\subset \mycal C(M_1\overline{\otimes}M_2)$ and similarly for $\mycal E$. 
Since $R\overline{\otimes} R\simeq R$ this implies $\mycal C(R)$, $\mycal E(R)$ are multiplicative semigroups inside the Jones spectrum $\{4\cos^2(\pi/n) \mid n\geq 3\}\cup [4, \infty)$, that contains the integers.

By Jones Theorem, $\{4\cos^2(\pi/n) \mid n\geq 3\}\subset \mycal C(R)$ and, as we have seen in 2.9 above, several other 
$\alpha\in \Bbb E_0^2$ have been shown to be contained in $\mycal C(R)$ as well. 
Solving the commuting square problem for all finite bipartite graphs would show that in fact $\Bbb E_0^2\subset \mycal E(R)$.

\heading 3. $W^*$-representations of subfactors \endheading

In this section we'll define the analogue for  a finite index subfactor $N\subset M$ of the notion of Hilbert-module 
(or $W^*$-representation) of a single II$_1$ factor $M$. Roughly speaking, this will be an 
inclusion of ``multi Hilbert modules'' ${\oplus_i} (_N\Cal K_i)
\hookrightarrow {\oplus_j}(_M\Cal H_j)$, a structure 
that's rigorously described as a nondegenerate embedding of $N\subset M$ into the atomic $W^*$-inclusion $\oplus_i \Cal B(\Cal K_i)\subset^{\Cal E} \oplus_j \Cal B(\Cal H_j)$.  

While in the case of a single factor $M$ a (left) Hilbert $M$-module 
$_M\Cal H$ (or $W^*$-representation $M \hookrightarrow \Cal B(\Cal H)$) comes 
with the Murray-von Neumann dimension $\text{\rm dim}(_M\Cal H)$ (or coupling constant $c_{M, M'}$ of $M \hookrightarrow\Cal B(\Cal H)$), 
the role of this for an inclusion of II$_1$ factors will be played by the dimension/coupling vector $(\text{\rm dim}(_M\Cal H_j))_j$. 

\vskip.05in

\noindent
{\bf 3.1. Some basic definitions}. A non-degenerate (normal) embedding  of a finite index extremal subfactor $N\subset M$  
 into an atomic $W^*$-inclusion $\Cal N \subset^{\Cal E} \Cal M$ is called 
a {\it $W^*$-representation} of  $N \subset M$. Thus, $\Cal N, \Cal M$ are direct sums of type I$_\infty$ factors, i.e., 
$\oplus_i \Cal B(\Cal K_i) = \Cal N\subset^{\Cal E} \Cal M=\oplus_j \Cal B(\Cal H_j)$. 
Note that by 2.8,  $\Cal N\subset^{\Cal E} \Cal M$ follows automatically $\lambda$-Markov for $\lambda=[M:N]^{-1}$. Thus, by Proposition 2.7.1, 
its {\it inclusion} (bipartite) {\it graph} $\Lambda=\Lambda_{\Cal N\subset \Cal M}$ satisfies $\|\Lambda\| \leq [M:N]^{1/2}.$  

The representation $\Cal N \subset^{\Cal E} \Cal M$ is {\it irreducible} if $\Cal M'\cap \Cal N=\Cal Z(\Cal M)\cap \Cal Z(\Cal N)=\Bbb C$, or equivalently if $\Lambda_{\Cal N \subset \Cal M}$ 
is connected as a graph (irreducible as a matrix). 

Two representations $\Cal N_l\subset^{\Cal E_l} \Cal M_l$, $l=0, 1$, are {\it equivalent} (or {\it isomorphic}) if there 
exists an isomorphism $\theta:\Cal M_0 \simeq \Cal M_1$, with $\theta(\Cal N_0)=\Cal N_1$, $\theta\circ \Cal E_0 = \Cal E_1\circ \theta$, that 
intertwines the corresponding embeddings of $N\subset M$. The two representations are {\it stably equivalent} (or {\it stably isomorphic}) if there exist 
projections $p_l\in M'\cap \Cal N_l$ such that $p_l\Cal N_lp_l \subset^{\Cal E_l(p_l\cdot p_l)} p_l\Cal M_l p_l$, $l=0, 1$, are equivalent. 

Stable isomorphism of $W^*$-representations involves ``reducing'' the commuting square embedding $(N\subset M)\subset (\Cal N\subset^{\Cal E}\Cal M)$ 
by  a projection in $M'\cap \Cal N$. This latter algebra is of course an isomorphism invariant for the representation. 

We call   $M'\cap \Cal N$ the {\it RC-algebra} (abbreviated for ``relative commutant $W$*-algebra'') of the $W^*$-representation. The 
interesting case is when $M'\cap \Cal N$ is a factor. If this is the case, then we say that $\Cal N\subset^\Cal E\Cal M$ is {\it RC-factorial}. 
Note that the  Murray-von Neumann {\it type}  of the RC-algebra/factor is an isomorphism invariant of the representation. More on this in Section 3.6.

Since $N\subset M$ is $\lambda$-Markov for $\lambda=[M:N]^{-1}$ and $W^*$-representations are non-degenerate embeddings, 
a representation $(N\subset M)\subset (\Cal N \subset^{\Cal E}\Cal M)$ 
is a $\lambda$-Markov commuting square embedding, so by Section 2.8 it gives rise 
to a tower of representations $(N\subset M \subset_{e_0}M_1 \subset_{e_1}  ...) \subset (\Cal N\subset^{\Cal E} \Cal M \subset_{e_0}^{\Cal E_1}\Cal M_1 \subset_{e_1}  ...)$. 

Since $M_{-1}=N\subset M=M_0$ are II$_1$ factors, by using the {\it downward basic construction} in (1.2.3 of [P92a]), one can 
choose (up to conjugacy by a unitary in $N=M_{-1}$) a projection $e_{-1}\in M_0$ with $E_N(e_{-1})=\lambda1$ and define $M_{-2}:=\{e_{-1}\}'\cap M_{-1}$.   
Then $M_{-2}$ follows a II$_1$ subfactor of index $\lambda^{-1}$ with the property that $M_{-2}\subset M_{-1}\subset_{e_{-1}} M_0$ is a basic construction. 
One can make such choices of Jones projections recusively, thus obtaining a {\it tunnel} of factors 
$ ... \subset_{e_{-2}} M_{-1} \subset_{e_{-1}} M_0$. 

Moreover, by (1.2.6-1.2.9 in [P93b]), if for each $i\leq -1$ we define $\Cal M_{i-1} = \{e_{i}\}'\cap \Cal M_i$ and $\Cal E_i:\Cal M_i \rightarrow \Cal M_{i-1}$ 
by $\Cal E_i(X)=\lambda \sum_j m^{i}_j X {m^i_j}^*$, where $\{m^i_j\}_j$ is an o.b. of $\{e_n\}_{n\geq i}''$ over $\{e_n\}_{n\geq i+1}''$ (in $\Cal M_\infty$), 
then  $(M_{i-1}\subset M_i \subset_{e_{i-1}} M_{i+1}) \subset (M_{-n-1}\subset M_{-n} \subset_{e_{-n}} M_{-n+1})$, $n\in \Bbb Z$, are all 
representations. We call such a double sequence of representations a {\it tower-tunnel} of representations. 

\vskip.05in
\noindent
{\bf 3.2. Two classes of examples}. This concept of $W^*$-representation of a subfactor was introduced in  (Section 2 in [P92a]), 
where one also notices the following class of examples (see Proposition 2.1 in [P92a]): 

\vskip.05in
\noindent 
{\bf 3.2.2. Example}. Let $N\subset M$ be an extremal inclusion of II$_1$ factors  with finite Jones index. 
Let $(Q\subset P)\subset (N\subset M)$ be a nondegenerate 
commuting square with  $Q \subset P$ finite dimensional. 
We will call such $Q\subset P$ a {\it graphage} of $N\subset M$. 
Since $N\subset M$ is $\lambda=[M:N]^{-1}$ Markov and the commuting square is non-degenerate, $Q\subset P$ follows $\lambda$-Markov as well. 
By the remark at the end of Section 2.8, if one takes the basic construction of this commuting square vertically and one denotes $\Cal N = \langle N, Q\rangle 
\subset^\Cal E \langle M, P\rangle=\Cal M$, where $\Cal E= E^{\langle M, P\rangle}_{\langle N, Q\rangle}$, then $\Cal E_{|M}=E^M_N$ 
and $(N\subset M)\subset (\Cal N\subset^\Cal E \Cal M)$ is a non-degenerate $W^*$-embedding. 
Since $\Cal N \subset \Cal M$ is an amplification of 
$Q\subset P$, it follows that $\Cal N, \Cal M$ are  atomic, so $\Cal N\subset^\Cal E \Cal M$ is a representation  of $N\subset M$. Moreover, 
since $\Cal E=E^{\langle M, P\rangle}_{\langle N, Q\rangle}$ preserves the canonical trace 
$\text{\rm Tr}_{\langle M, P\rangle}$, we have that $\Cal N\subset^\Cal E \Cal M$ is both $\lambda$-Markov and Tracial. Also, since $\Cal N\subset \Cal M$ is an amplification 
of $Q\subset P$,  the  
bipartite graph $\Lambda_{\Cal N\subset \Cal M}$ identifies naturally 
with $\Lambda_{Q\subset P}$, so the representation is irreducible iff $\Lambda_{Q\subset P}$ is irreducible (equivalently, $\Cal Z(P)\cap \Cal Z(Q)=\Bbb C$). 
Also,  we have $[M:N]=\|\Lambda_{\Cal N\subset \Cal M}\|^2=\|\Lambda_{Q\subset P}\|^2$.  

\vskip.05in 

Another class of examples of $W^*$-representations of a given subfactor 
$N\subset M$ comes from the following trivial observation: 

\proclaim{3.2.2. Lemma}  If $(N\subset M) \subset (\tilde{N}\subset^{\tilde{E}} \tilde{M})$ 
is a non-degenerate commuting square embedding of the extremal inclusion of $\text{\rm II}_1$ factors 
$N\subset M$ into another inclusion of factors with expectation, 
then any representation $(\tilde{N}\subset^{\tilde{E}} \tilde{M})\subset (\Cal N\subset^\Cal E \Cal M)$ $($that is, a non-degenerate 
embedding of $\tilde{N}\subset^{\tilde{E}}\tilde{M}$ into an atomic $W^*$-inclusion with expectation $\Cal N\subset^{\Cal E}\Cal M)$ gives a $W^*$-representation 
$(N\subset M)\subset (\Cal N\subset^{\Cal E}\Cal M)$, by composing the embeddings. 
\endproclaim

\vskip.05in

\noindent
{\bf 3.3. Tracial representations}. A representation $\Cal N \subset^{\Cal E} \Cal M$ of a subfactor $N\subset M$ is {\it Tracial} 
if the atomic $W^*$-inclusion $\Cal N\subset^{\Cal E} \Cal M$ is Tracial, i.e., there exists an n.s.f. trace Tr on $\Cal M$ such that 
Tr$\circ \Cal E=$Tr. As seen in 3.1.2 above, a representation arising from a graphage $(Q\subset P)\subset (N\subset M)$ does have this property. 
The existence of a graphage  is a rather strong structural property of $N\subset M$, which in particular implies $[M:N]=\|\Lambda_{Q\subset P}\|^2 \in \Bbb E^2$. 
So a II$_1$ subfactor $N\subset M$ with $A_\infty$-graph and index in the set $(4, 2+\sqrt{5})\setminus \Bbb E^2$ (which by [P90] exists  
for any $\alpha$ lying in this set) does not have any graphage. A Tracial representation can be viewed as a 
``dim graphage'' of $N\subset M$. 

Note that the tower/tunnel of reps associated with a Tracial $W^*$-representation $(N\subset M)\subset (\Cal N\subset^{\Cal E}\Cal M)$ are all Tracial.

\proclaim{3.3.1. Proposition} Let $(N\subset M) \subset (\Cal N\subset^{\Cal E}\Cal M)$ 
be an irreducible representation with finite inclusion graph $\Lambda=\Lambda_{\Cal N\subset \Cal M}$. 

\vskip.05in

$1^\circ$ If $\text{\rm Tr}$ is a n.s.f.  $\lambda=[M:N]^{-1}$ Markov trace on $\Cal N \subset \Cal M$, then it is necessarily $\Cal E$-invariant. 

$2^\circ$  $\Cal N \subset^{\Cal E} \Cal M$ is Tracial if and only if $\|\Lambda\|^2=[M:N]$.
\endproclaim
\noindent
{\it Proof}. By Proposition 2.7.1, any n.s.f. trace Tr on $\Cal M$ for which the Tr-preserving expectation $\Cal E'$ is $\lambda'$-Markov, for some $\lambda'>0$, 
forces $\lambda'=\|\Lambda\|^{-2}$ and Tr be given by a weight vector proportional to the (unique) Perron-Frobenius eigenvector of $\Lambda^t\Lambda$ 
corresponding to eigenvalue $\|\Lambda\|^2$. Thus, condition $1^\circ$ implies $\|\Lambda\|^2=[M:N]$. In particular,  Ind$\Cal E = \text{\rm Ind}(\Cal E')$. Since $\Cal E'$ is the unique 
expectation with index $\|\Lambda\|^2$, this implies $\Cal E=\Cal E'$. 

This also proves $\Leftarrow$ in $2^\circ$, while the opposite implication follows from $2.6.1$. 
\hfill $\square$

\vskip.05in

\noindent
{\bf 3.4. The coupling vector of a representation}. Given an $(N\subset M)$-representa- tion 
$ \oplus_{i\in I} \Cal B(\Cal K_i) = \Cal N\subset^{\Cal E} \Cal M=\oplus_{j\in J} \Cal B(\Cal H_j)$, 
we denote $\vec{d}_{M}(\Cal N \subset \Cal M)=(d_M(j))_{j\in J}$, 
resp. $\vec{d}_N(\Cal N \subset \Cal M)=(d_N(i))_{i\in I}$, the vectors  with entries $d_M(j)=\text{\rm dim}(_M\Cal H_j)\in (0, \infty]$, 
$d_N(i)=\text{\rm dim}(_N\Cal K_i)\in (0, \infty]$, and call them    
the {\it coupling vectors} (or {\it dimension vectors}) of the representation. If 
$d_M(j)<\infty$, $d_N(i)<\infty$, $\forall i, j$, then we say that the representation has 
{\it finite couplings} (or {\it finite dimension vectors}). 

We'll next prove that if this is the case, then $\Cal N\subset^{\Cal E}\Cal M$ is automatically Tracial whenever $\Lambda_{\Cal N\subset \Cal M}$ is finite, 
with the dimension vector $\vec{d}_M$ giving the weights of the $\Cal E$-invariant n.s.f. trace Tr. Another important class of Tracial representations 
with finite coupling vector giving the weights of Tr will be discussed in Sec 3.8. 

Recall first some well known facts about the dimension of Hilbert modules over a II$_1$ factor and the way 
it relates to Jones index (see e.g. [J82]). 

\proclaim{3.4.1. Lemma} Let $N\subset M$ be an inclusion of $\text{\rm II}_1$ factors and $\Cal H, \Cal H'$ some $($left$)$ 
Hilbert $M$-modules. Then we have: 
\vskip.05in 
$(a)$ $\text{\rm dim}(_M(\Cal H\oplus \Cal H'))=\text{\rm dim}(_M\Cal H)+\text{\rm dim}(_M\Cal H')$.

$(b)$ When viewing $\Cal H$ as an $N$-module, one has $\text{\rm dim}(_N\Cal H)=[M:N]\text{\rm dim}(_M\Cal H)$. 

$(c)$ If $p$ is a projection in $M$ then $\text{\rm dim}(_{pMp}p(\Cal H))=\tau(p)^{-1}\text{\rm dim}(_M\Cal H)$. 
 \endproclaim

\proclaim{3.4.2. Proposition} Let $\oplus_l \Cal B(\Cal K_i) = \Cal N\subset^{\Cal E} \Cal M=\oplus_j \Cal B(\Cal H_j)$ be an irreducible representation  of 
$N\subset M$ with inclusion graph  $\Lambda=\Lambda_{\Cal N\subset \Cal M}$ and dimension vectors $\vec{d}_{M}=(d_M(j))_{j\in J}$, 
$\vec{d}_N=(d_N(i))_{i\in I}$. If $d_M(j_0)<\infty$,  or $d_N(i_0)<\infty$, for some $j_0\in J$, or $i_0\in I$, then  the representation has  finite couplings  
and we have:

\vskip.05in
$1^\circ$ $\vec{d}_N=\Lambda(\vec{d}_M)$, $\Lambda^t(\vec{d}_N)=[M:N]\vec{d}_M$, $\Lambda^t\Lambda(\vec{d}_M)=[M:N]\vec{d}_M$. 

$2^\circ$ The n.s.f. trace $\text{\rm Tr}$ on $\Cal M$ given by the weight vector $\vec{d}_M$ is a $\lambda=[M:N]^{-1}$ Markov n.s.f. trace for $\Cal N \subset \Cal M$. 

\vskip.05in
$3^\circ$ If $\Lambda$ is finite, then the n.s.f. trace $\text{\rm Tr}$ on $\Cal M$ given by the weights $\vec{d}_M$ is $\Cal E$-invariant. 

\endproclaim 
\noindent
{\it Proof}. Note that proving 
the ``entry by entry'' equalities $\Lambda(\vec{d_M})=\vec{d}_N$ 
and $\Lambda^t(\vec{d}_N)=[M:N]\vec{d}_M$ with the entries being in $(0, \infty]$ (so apriori not all finite) implies both the fact that 
``$d_M(j_0)<\infty$, for some $j_0\in J$, or $d_N(i_0)<\infty$, for some $i_0\in I$, implies all entries of both $\vec{d}_M, \vec{d}_N$ are finite'' and the fact 
that $\Lambda^t\Lambda(\vec{d}_M)=[M:N]\vec{d}_M$. 

Let $\Lambda=(b_{ij})_{i\in I, j\in J}$ as usual. Since $\Cal H_j = \oplus_i \Cal K_i^{\oplus b_{ij}}$, by parts $(a), (b)$ in Lemma 3.4.1 we have 
$[M:N]\text{\rm dim}(_M\Cal H_j)=\text{\rm dim}(_N\Cal H_j)=\sum_i b_{ij} \text{\rm dim}(_N\Cal K_i)$, showing that  
$\Lambda^t(\vec{d}_N)=[M:N]\vec{d}_M$. 

Let $\Cal N\subset\Cal M \subset^{\Cal E_1}_e \Cal M_1=\oplus_{i\in I}\Cal B(\Cal K_i')$ be the basic construction for 
$\Cal N \subset^{\Cal E}\Cal M$, with $N\subset M \subset_eM_1$ represented in it. Note that $(N\subset \Cal N)$ 
is isomorphic to $e(M_1\subset \Cal M_1)e=(eM_1e\subset e\Cal M_1e)$ and that $e\Cal B(\Cal K'_i)e=\Cal B(e(\Cal K'_i))$, 
so that $N\subset \Cal B(\Cal K_i)$ is ``same as'' $eM_1e\subset \Cal B(e(\Cal K_i'))$. By part $(c)$ in Lemma 3.4.1, since 
$\tau_{M_1}(e)=[M:N]^{-1}$, it follows that 
$$
\text{\rm dim}(_{M_1}\Cal K'_i)=\tau_{M_1}(e)\text{\rm dim}(_{eM_1e}e(\Cal K'_i))
=[M:N]^{-1}\text{\rm dim}(_N\Cal K_i) \tag 3.4.2.1
$$ 
implying that for the dimension vectors we have 
$$
\vec{d}_{M_1}=[M:N]^{-1} \vec{d}_N \tag 3.4.2.2
$$
By applying the first part of the proof to the representation of  $M\subset M_1$ into $\oplus_j \Cal B(\Cal H_j)=\Cal M \subset \Cal M_1=\oplus_i \Cal B(\Cal K'_i)$, 
with its inclusion bipartite graph/matrix $\Lambda_{\Cal M \subset \Cal M_1}$ identified with $\Lambda^t=((b_{ij})_{i,j})^t$, 
it follows that $\Lambda(\vec{d}_M)=(\Lambda_{\Cal M\subset \Cal M_1})^t \vec{d}_M=[M_1:M]\vec{d}_{M_1}$, so  
by $(3.4.2.2)$ and by using the fact that $[M_1:M]=[M:N]$, we get $\vec{d}_N=\Lambda(\vec{d}_M)$, end the rest of the equalities in $1^\circ$.  

Part 2$^\circ$ is now an immediate consequence of part 1$^\circ$ and of Lemma 2.7.1.2$^\circ$, while Part $3^\circ$ follows from Part $2^\circ$ 
and Proposition 3.3.1. 

\hfill $\square$

\noindent
{\bf 3.4.3. Definition}. If a $W^*$-representation $(N\subset M)\subset (\Cal N\subset^{\Cal E} \Cal M)$ has finite couplings and 
the n.s.f. trace Tr implemented by its coupling vector is $\Cal E$-invariant, then we say that it is {\it canonically Tracial}. 

\vskip.05in

\noindent
{\bf 3.5. Smooth representations}. A representation $\Cal N\subset^{\Cal E} \Cal M$ of a subfactor $N\subset M$ is called {\it smooth} 
if in the associated tower of representations $(N\subset M \subset_{e_0} M_1 \subset ...)\subset (\Cal N\subset^{\Cal E}\Cal M \subset_{e_0}^{\Cal E_1} \Cal M_1 \subset ...)$ 
one has $N'\cap M_j \subset \Cal N'\cap \Cal M_j$ (see Section 2.3 in [P92a]). It is trivial to see that this condition 
implies $M_i'\cap M_j \subset \Cal M_i'\cap \Cal M_j$, $\forall j\geq i \geq -1$ 
(where $\Cal M_{-1}=\Cal N$, $\Cal M_0=\Cal M$, $M_{-1}=N, M_0=M$) and that these relations are equivalent to $\Cal M_i'\cap M_j =M_i'\cap M_j$, $\forall i,j$. 

This compatibility relation between the higher relative commutants of $N\subset M$ and $\Cal N \subset \Cal M$ is natural to impose,  
a fact that's amply emphasized by results in ([P92a], [P97a]). More generally, 
a non-degenerate commuting square embedding of $N\subset M$ into an arbitrary $W^*$-inclusion with expectation $\Cal N\subset^{\Cal E}\Cal M$ 
is smooth if $M_i'\cap M_j \subset \Cal M_i'\cap \Cal M_j$, $\forall j\geq i \geq -1$, with the same convention of notations as above.

\vskip.05in

\noindent
{\bf 3.6. Exact representations}. Following (Section 2.4 in [P92a]),  
a representation $\Cal N\subset^{\Cal E}\Cal M$ for a subfactor $N\subset M$ is  called {\it exact},   
if $\Cal M=M\vee {M'\cap \Cal N}$. We then also say that it has the {\it relative bi-commutant} property. 

We briefly recall below some facts about exact representations, referring to (Section 2.4 in [P92a]) for more details.

Note that an exact representation is irreducible iff $\Cal P=M'\cap \Cal N$ is a factor 
and that if this is the case then each factorial type I$_\infty$ direct summand $\Cal B(\Cal H)$ of $\Cal M$ (resp. $\Cal B(\Cal K)$ of $\Cal N$) 
is an irreducible binormal representation of $M \otimes \Cal P$ (resp. $N\otimes \Cal P$), equivalently an irreducible Hilbert bimodule 
$_M\Cal H_{\Cal P^{op}}$ (resp. $_N\Cal K_{\Cal P^{op}}$). We'll call $\Cal P=M'\cap \Cal N$ the {\it RC-factor} (or {\it exacting factor}) 
of the exact representation. 

Let $\Cal N \subset \Cal M \subset_{e_0}^{\Cal E_1} \Cal M_1 \subset ...$ be the 
Jones tower for $\Cal N \subset \Cal M$ and denote by $\Cal M_\infty$ its enveloping von Neumann algebra. We then  have 
$M_n'\cap  \Cal N=M'\cap \Cal N = \Cal P$ for all $0\leq n \leq \infty$ and $\Cal M_n = M_n \vee \Cal P$, for all $-1\leq n \leq \infty$. 
In particular, this shows that an exact representation is smooth.

Exact representations arise  concretely as follows. 
For simplicity, let us assume $N, M$ are separable II$_1$ factors and we only look for separable exact representations of $N\subset M$. 
Take $P$ to be a subfactor of $M^{op}\overline{\otimes} \Cal B(\ell^2\Bbb N)$, 
like for instance $P=M^{op}$. Denote by $M \otimes_{bin} P$ the completion of $M \otimes P$ in the maximal binormal C$^*$-norm, obtained from all binormal representations of $M\otimes P$. 
It is easy to see that the norm it induces on the subalgebra $N\otimes P$ is equal to its own maximal binormal C$^*$-norm, 
thus identifying $N\otimes_{bin} P$ as a C$^*$-subalgebra of $M\otimes_{bin}P$. 
Moreover, the map $E_N\otimes id_P: M\otimes P \rightarrow N\otimes P$ 
extends to a conditional expectation, still denoted $E_N \otimes id$, of $M \otimes_{bin} P$ onto $N\otimes_{bin} P$. Since 
the inequality $E_N(x)\geq \lambda x$, $\forall x\in M_+$, is stable (i.e., difference is completely positive), 
it follows that one still has $E_N\otimes id_P(X)\geq \lambda X$, $\forall X\in (M\otimes_{bin}P)_+$. 
This expectation extends to an expectation from the von Neumann algebras they entail. To see this rigorously, 
note first that by taking biduals of $N\otimes_{bin} P \subset^{E_N\otimes id}M \otimes_{bin} P$ one gets an inclusion of von Neumann algebras with expectation 
$(N\otimes_{bin} P)^{**} \subset^{\Cal E}(M \otimes_{bin} P)^{**}$, where $\Cal E=(E_N\otimes id)^{**}$ still satisfies $\Cal E(X)\geq \lambda X$  for all $X\geq 0$. 
Due to weak density of $N\otimes_{bin} P \subset M \otimes_{bin} P$ inside it, for which one has the formula $X=\sum_j m_j \Cal E(m_j^*X)$, 
$\forall X \in M \otimes_{bin} P$, one has this formula for all $X \in (M \otimes_{bin} P)^{**}$. 

By (Section 2.4 in [P92a], or 1.1.2$(iii)$ in [P93b]), 
the inequality $\Cal E(X)\geq \lambda X$, $\forall X\geq 0$, insures that the central support of the atomic parts of $(N\otimes_{bin} P)^{**}, (M \otimes_{bin} P)^{**}$ coincide, 
and so do the central supports of the parts where $M, P$ (resp. $N, P$) are represented normally, and for which each factorial direct summand is separable. Thus, 
if one denotes $\Cal N^u_P \subset^{\Cal E} \Cal M^u_P$ this atomic inclusion, then each factorial direct summand of $\Cal M^u_P$ (resp. $\Cal N^u_P$) 
is an irreducible binormal representation  of $M \otimes_{bin}P$ (resp. $N\otimes_{bin} P$).  

Obviously, any irreducible separable exact 
representation of $N\subset M$ arises this way.  Doing this construction for all subfactors $P\subset M^{op}\overline{\otimes}\Cal B(\ell^2\Bbb N)$, 
then choosing one irreducible representation for each isomorphism class of such a rep., then taking direct sum, 
gives a representation of $N\subset M$ that we denote $\Cal N^{u} \subset^{\Cal E}\Cal M^{u}$ and that we call the {\it universal exact} (or {\it binormal})  
$W^*$-representation of $N\subset M$.   

One should note that if $\Cal N \subset \Cal M$ arises as an irreducible $\Cal N_P \subset \Cal M_P$, for some $P\subset M^{op}\overline{\otimes} \Cal B(\ell^2\Bbb N)$, 
then one has $\Cal P=M'\cap \Cal N\supset P$ but the inclusion can be strict (see Remark 3.6.1 for concrete examples). 
We call $\Cal P$ the {\it RC-envelope} of $P$.

We denote $\Lambda^{u}=\Lambda^u_{N\subset M}=(b_{lk})_{l\in L^u, k\in K^u}$ the inclusion graph of $\Cal N^u \subset \Cal M^u$,  where  
$K^u$ labels the set of atoms in $\Cal Z(\Cal M^{u})$, $L^u$ labels the set of atoms of $\Cal Z(\Cal N^{u})$. 
One also denotes by $\Cal H_k$ (resp $\Cal K_l$) the irreducible binormal $M-P$ (resp. $N-P$) Hilbert bimodule corresponding to the atom $k\in K^u$ 
(resp. $l\in L^u$).  

Note that $\Cal Z(\Cal M^{u})\cap \Cal Z(\Cal N^{u})$ is atomic, in fact any atom of $\Cal Z(\Cal M^{u}), \Cal Z(\Cal N^{u})$ 
is majorized by a (unique) atom  of this intersection, corresponding to a connected component of $\Lambda^{u}$. Each such 
direct summand $(\Cal N\subset \Cal M)=(\Cal N^{u}\subset \Cal M^{u})s(q_0)$ 
gives an irreducible representation of $N\subset M$, i.e. with $\Cal Z(\Cal N)\cap \Cal Z(\Cal M)=\Bbb C$. 

The construction of an irreducible sub-representation $\Cal N \subset \Cal M$ of  $(\Cal N^{u} \subset^{\Cal E}\Cal M^{u})$ can 
be obtained more directly as follows. Start with  a (separable) irreducible Hilbert $M-P$ bimodule $\Cal H_0$, corresponding to some central atom 
$q_0 \in \Cal Z(\Cal M^{u})$, labelled by $0\in K^u$ (N.B. by Connes theory of correspondences, this is equivalent to  
an embedding with trivial relative commutant $P\hookrightarrow (M^{op})^\alpha$, for some $0< \alpha \leq \infty=\aleph_0$, see [P86]). Set $J_0=\{0\}$. 
Then take all irreducible $N-P$ Hilbert subbimodules $\Cal K_i$ appearing in $_N(\Cal H_0)_P$, indexed by  $i\in   I_1$, denoting $b_{i0}$ its multiplicity. 
Then take all irreducible $M-P$ Hilbert bimodules $\Cal H_j$ with $_N(\Cal K_i)_P\leq$ $_N(\Cal H_j)_P$, for some $i\in I_1$,  
indexed by $j\in J_1$, or equivalently the irreducible direct summands of 
$_M(L^2M\otimes_N \Cal H_0)_P$. Also, denote $b_{ij}$ the corresponding multiplicity. One continues recursively, with $_N(\Cal K_i)_P, i\in I_n$, being  
the irreducible  $N-P$ submodules of $_N(L^2M^{\otimes_N (n-1)}\otimes_N \Cal H_{0})_P$ and $_M(\Cal H_j)_P, j \in J_n,$ 
the irreducible $M-P$ submodules of $_M(L^2M^{\otimes_N n}\otimes_N \Cal H_0)_P$, while denoting $b_{ij}$ the multiplicity of this sub bimodule. 
One has natural  identifications $I_n\hookrightarrow I_{n+1}$,  
$J_n\hookrightarrow J_{n+1}$, due to 2-periodicity in the Jones tower and the fact that $_{Q}{L^2(M_n)}_N$ is isomorphic to $_{Q}(L^2M^{{\otimes_N} (n+1)})_N$, 
for any $n\geq 0$, where $Q \in \{N, M\}$. If one denotes $J=\cup_n J_n, I=\cup_n I_n$, then it is immediate to see 
that $(\Cal N \subset \Cal M)=(\oplus_{i\in I}\Cal B(\Cal K_i) \subset \oplus_{j\in J} \Cal B(\Cal H_j))$, with the inclusion bipartite graph given by $\Lambda=(b_{ij})_{i \in I, j\in J}$  
(see Section 2.4 in [P92a] for more on this).

\vskip.05in

\noindent
{\bf 3.7. The standard $W^*$-representations.} (Section 2.4 in [P92a] and Section 5 in [P97a]). The ``concrete construction'' 
of exact representations of $N\subset M$ from irreducible subfactors $P$ of amplifications of $M$, shows 
that one has ``plenty'' of examples of such representations, 
as many as irreducible subfactors  of $M^\alpha, 0< \alpha \leq \infty$, one can have (see however 
Problems 6.1.6, 6.1.7 and accompanying remarks). 

The simplest case of such a construction is when 
$P=M$ and $\Cal H_0$ is the Hilbert bimodule $_ML^2M_M$, i.e., the standard representation of $M$. 
We call the corresponding irreducible direct summand of $\Cal N^u \subset \Cal M^u$ the {\it standard $W^*$-representation}  
of $N\subset M$ and denote it $\Cal N^{st}\subset^{\Cal E^{st}}\Cal M^{st}$. So the RC-factor $\Cal P$ for the standard representation 
of $N\subset M$ is $M^{op}$ itself. 

One quick way to describe this representation is to consider in $\Cal B(L^2(M_\infty))$ 
the inclusion $\Cal N^{st}=N\vee M^{op} \subset M\vee M^{op}=\Cal M^{st}$, with the expectation $\Cal E^{st}$ given by 
$X \mapsto \sum_j m_j e_0Xe_0m_j^*$, for $X\in M\vee M^{op}$, where $\{m_j\}_j$ is here an o.b. of 
$\{e_n\}_{n\geq 1}''$ over $\{e_n\}_{n\geq 2}''$ and as usual $N\subset M\subset_{e_0}M_1 \subset_{e_1}M_2 \subset ... \nearrow M_\infty$ 
is the Jones tower for $N\subset M$. This coincides withe the unique expectation extending $E^M_N \otimes id_{M^{op}}$. 

One can show that, up to isomorphism of representations, $\Cal N^{st}\subset^{\Cal E^{st}}\Cal M^{st}$ is 
the unique irreducible exact representation $\Cal N \subset^{\Cal E} \Cal M$  of $N\subset M$ 
for which there exists a direct summand $\Cal B(\Cal H_0)$ of $\Cal M$  such that 
dim$(_M\Cal H_0)=1$ (equivalently, $_M\Cal H_0=_ML^2M$)  and such that  after 
cutting by the support projection of $\Cal B(\Cal H_0)$ in $\Cal Z(\Cal M)$, one has   $\Cal P=M'\cap \Cal B(\Cal H_0) 
\simeq M'\cap \Cal B(L^2M)=M^{op}$. 

Recall from (Section 2.4 in [P92a]) that $\Cal N^{st}\subset^{\Cal E^{st}}\Cal M^{st}$ is Tracial, with inclusion graph $\Lambda=\Lambda_{\Cal N^{st}\subset \Cal M^{st}}$ 
naturally identifying with the transpose of the standard graph of $N\subset M$, $\Gamma_{N\subset M} = (a_{kl})_{k\in K, l\in L}$. Namely, $\Lambda=(b_{lk})_{l\in L, k\in K}$, 
where $b_{lk}=a_{kl}$, the ``pointed'' set $*\in K$ labeling the irreducible subfactors in the ``even levels'' of the Jones tower, $M \subset M_{2n}$, 
and $L$ labeling the  irreducible subfactors in the ``odd levels'' of the Jones tower, $N \subset M_{2n}$. More precisely, one has 
$\oplus_{l\in \Cal K_l} \Cal B(\Cal K_l) = \Cal N^{st} \subset \Cal M^{st}=\oplus_{k\in K}\Cal B(\Cal H_k)$, where 
$\{\Cal H_k\}_{k\in K}$  (resp. $\{\Cal K_l\}_{l\in L}$) is the list of all irreducible Hilbert $M-M$ (resp. $N-M$) bimodules in $_ML^2(M_n)_M$ 
(resp. $_NL^2(M_n)_M$), $n\geq 0$. The weight vectors giving the canonical $\Cal E^{st}$-preserving n.s.f. trace Tr 
were denoted in (2.4 of [P92a]) as $\vec{v}=(v_k)_{k\in K}$ and $\vec{u}=(u_l)_{l\in L}$, with $v_k$ (resp. $u_l$) 
the square root of the index of the irreducible subfactor $Mp\subset pM_{2n}p$ (resp. $Nq \subset qM_{2n}q$), where $p$  (resp. $q$) is a minimal 
projection in $M'\cap M_{2n}$ (resp $N'\cap M_{2n})$ labelled by $k\in K_n$ (resp. $l\in L_n)$,  
the significance of $K_n \subset K, L_n \subset L$ being as in Section 2.3. 

Note that the considerations in (Section 2.4 of [P92a]) show that $\vec{v}, \vec{u}$ coincide with  the dimension vectors $\vec{d}_M, \vec{d}_N$ 
of the representation $(N\subset M) \subset (\Cal N^{st} \subset^{\Cal E^{st}}\Cal M^{st})$, i.e., 
$v_k=\text{\rm dim}(_M\Cal H_k)$, $u_l=\text{\rm dim}(_N\Cal K_l)$, $\forall k \in K, l\in L$. Thus, with the terminology we introduced in Sec. 3.3, 
the standard representation  has finite couplings and it is  canonically Tracial, i.e., the n.s.f. trace Tr given by 
the coupling vector is $\Cal E^{st}$-invariant. 

We also consider the {\it dual standard $W^*$-representation} of $N\subset M$, defined  again on $L^2(M_\infty)$ by $(N\vee N^{op}=\Cal N^{st'}\subset^{\Cal E^{st'}} \Cal M^{st'} = M\vee N^{op})$,   
with the expectation $\Cal E^{st'}$ being the unique expectation extending $E^M_N \otimes id_{N^{op}}$. With the same reasoning as above, 
its inclusion graph $\Lambda'=\Lambda_{\Cal N^{st'} \subset \Cal M^{st'}}$ identifies naturally with the dual standard graph $\Gamma'_{N\subset M}$, which in turn coincides with the standard graph $\Gamma_{N^{op}\subset M^{op}}$, 
of the opposite subfactor $(N\subset M)^{op}=(N^{op}\subset M^{op})$  (see Section 2.4 in [P92a]).  Also, this representation has finite couplings and it is canonically Tracial, 
with the coupling vectors given by the canonical weights $\vec{u'}, \vec{v'}$ of $\Gamma'_{N\subset M}$.

\vskip.05in
\noindent
{\bf 3.7.1. Definition}. The atomic $W^*$-inclusions involved in the standard representation and its dual form a natural commuting square embedding $(\Cal N^{st'}\subset^{\Cal E^{st'}}\Cal M^{st'})
\subset (\Cal N^{st}\subset^{\Cal E^{st}}\Cal M^{st})$ given by 
$$
\CD
N\vee M^{op}\ @.\subset\ @.M\vee M^{op}\\
\noalign{\vskip-6pt}
\cup\ @.\ @.\cup\ \\
\noalign{\vskip-6pt}
N\vee N^{op}\ @.\subset\ @.M\vee N^{op}
\endCD
$$
with the vertical expectations being extensions of $id_{N}\otimes E^{M^{op}}_{N^{op}}$ resp. $id_M\otimes E^{M^{op}}_{N^{op}}$, with all expectations involved preserving 
the canonical n.s.f. trace Tr on $M\vee M^{op}$,  and with vertical 
inclusion graphs given by $\Lambda_{N\vee N^{op}\subset N\vee M^{op}}=\Gamma_{N\subset M}$, respectively $\Lambda_{M\vee N^{op}\subset M\vee M^{op}}=(\Gamma'_{N\subset M})^t$. 
This object is obviously an isomorphism  invariant of the subfactor $N\subset M$. We call it the {\it standard $\lambda$-commuting square} 
(or {\it standard $\lambda$-cell}) of $N\subset M$, and denote it $\Cal C_{N\subset M}$. 

\vskip.05in

Let us note that, as an invariant of $N\subset M$, the standard commuting square $\Cal C_{N\subset M}$ 
contains the same amount of information as (i.e., it is equivalent to) the standard invariant $\Cal G_{N\subset M}$.  
Recall in this respect that  an isomorphism 
of standard invariants (viewed as abstract objects like in [P94]) 
means a trace preserving isomorphism between the union of  the finite dimensional algebras involved  which takes the $_{ij}$  
algebras of the lattices one onto the other and the Jones-sequences of 
projections one onto the other. For standard $\lambda$-cells, an isomorphism is the usual notion of Tr-preserving isomorphism of commuting squares.

\proclaim{3.7.2. Theorem} $1^\circ$ Let $N\subset M$ be an extremal inclusion of $\text{\rm II}_1$ factors. Then  $\Cal C_{N\subset M}$ 
identifies naturally with the Tracial Markov commuting square $(\Cal A_{-1}^{-1}\subset \Cal A_0^{-1})\subset (\Cal A_{-1}^0\subset \Cal A_0^0)$ 
in Section $2$ of $\text{\rm [PS01]}$. 

$2^\circ$ If  $(Q\subset P)$ is another extremal inclusions of $\text{\rm II}_1$ factors, then 
$\Cal C_{N\subset M}\simeq \Cal C_{Q\subset P}$ if and only if $\Cal G_{N\subset M}\simeq \Cal G_{Q\subset P}$. 
\endproclaim
\noindent
{\it Proof}. Part $1^\circ$ is trivial by the definition of $\Cal C_{N\subset M}$ and 
by the way the commuting square $(\Cal A_{-1}^{-1}\subset \Cal A_0^{-1})\subset (\Cal A_{-1}^0\subset \Cal A_0^0)$ is constructed in ([PS01]).

$2^\circ$ Letting $M_{-1}=N \subset M=M_0\subset_{e_0}M_1 \subset_{e_1} M_2 \subset ...$ be the Jones tower for $N\subset M$, 
one has a Jones tower for the (Markov) standard $\lambda$-cell in 3.8.1 
$$
\CD
M_{-1}\vee M^{op}\ @.\subset\ @.M_0\vee M^{op}\ @.\subset_{e_0}\  @.M_1 \vee M^{op}\ ...@.\  \\
\noalign{\vskip-6pt}
\cup\ @.\ @.\cup\ @.\  @.\cup\  @.\  \ \\
\noalign{\vskip-6pt}
M_{-1} \vee N^{op}\ @.\subset\ @.M_0\vee N^{op}\ @.\subset_{e_0}\  @.M_1\vee N^{op}\ ....@.\  
\endCD
$$
with the Jones tower $M_{-1}\subset M_0 \subset_{e_0} M_1 \subset_{e_1} ...$  represented in the bottom row (and in the top row 
by composition of embeddings). It is immediate to see that $(M_i\vee M^{op})'\cap (M_j \vee N^{op})=M_i'\cap M_j$, with the Jones 
projections appearing in the algebras $M_i'\cap M_j$ same way as they do in the higher relative commutants. Moreover, 
it is easy to see that the compositions of expectations $\Cal E^{st}_i, i\geq 0$ for the top row implement on these 
higher relative commutants  the trace state $\tau$ they inherit from $M_\infty$. Since all these data comes from 
$\Cal C_{N\subset M}$ viewed as an abstract object, it follows that $\Cal C_{N\subset M}$ uniquely determines $\Cal G_{N\subset M}$. 

Conversely, if $\Cal G_{N\subset M}$ is given, then Lemma 2.1 in [PS01] associates to it in a canonical way a Tracial Markov commuting 
$(\Cal A_{-1}^{-1}\subset \Cal A_0^{-1})\subset (\Cal A_{-1}^0\subset \Cal A_0^0)$, which by part $1^\circ$ is the same as $\Cal C_{N\subset M}$. 
\hfill $\square$

\vskip.05in 

\noindent
{\bf 3.7.3. Remark}. We mentioned that one can have exact representations $\Cal N_P \subset \Cal M_P$ where $P$ is strictly smaller than 
the RC-factor/envelope. Indeed, if $N\subset M$ is any extremal inclusion of separable II$_1$ factors with finite index and 
Jones tower $N\subset M\subset M_1 \subset M_2 ...\nearrow M_\infty$, then by ([P18b]) there exists a decreasing sequence of hyperfinite subfactors 
$R_n\subset M$ such that $\cap_n R_n = \Bbb C1$ and $R_n'\cap M_j=M'\cap M_j$, for all $j\leq \infty$. Thus, given any $P=R_n$ 
we have $(N\subset M) \subset (\Cal N_P\subset \Cal M_P)=(N\vee P^{op}\subset M\vee P^{op})$ is isomorphic to the standard representation 
$(N\subset M) \subset (\Cal N^{st}\subset \Cal M^{st})=(N\vee M^{op}\subset M\vee M^{op})$.

\vskip.05in

\vskip.05in

\noindent
{\bf 3.8. Exact representations with finite couplings.} We discuss in this section the class of exact representations with 
finite couplings, showing  they are canonically Tracial. We also prove   
that, after an appropriate amplification, the RC-factor $\Cal P$ can be identified 
with the mirror image $P^{op}$ of an irreducible subfactor $P\subset M$ 
satisfying the bicommutant condition $(P'\cap M_\infty)'\cap M=P$.  
The standard representation is a particular case of this class of exact representations, corresponding to $P=M$, while the 
dual standard representation corresponds to the case $P=N$. 

Thus, let $\oplus_{i\in I} \Cal B(\Cal K_i)=\Cal N\subset^{\Cal E} \Cal M=\oplus_{j\in J} \Cal B(\Cal H_j)$ be an irreducible exact representation of 
the subfactor $N\subset M$, with inclusion bipartite graph/matrix 
$\Lambda=\Lambda_{\Cal N\subset \Cal M}=(b_{ij})_{i\in I, j\in J}$. As usual,  
denote $\Cal P=M'\cap \Cal N$. Recall from Proposition 3.4.2 that $\Cal N\subset \Cal M$ has finite dimension vectors $\vec{d}_M, \vec{d}_N$ 
iff there exists $j_0\in J$ such that $\alpha=d_M(j_0)=\text{\rm dim}(_M\Cal H_{j_0})< \infty$. Note that if this is the case then $\Cal P$  is a II$_1$ factor
which identifies with an irreducible II$_1$ subfactor of the $\alpha$-amplification of $M^{op}$. 

Let $t> 0$. Replacing $\Cal N\subset^{\Cal E} \Cal M$ 
by $\Cal N \otimes \Bbb M_n(\Bbb C) \subset^{\Cal E \otimes id} \Cal M \otimes \Bbb M_n(\Bbb C)$ for $n\geq t$ and then taking  
$p$ to be a projection in the II$_1$ factor $\Cal P\otimes \Bbb M_n(\Bbb C)$ of normalized trace $\tau(p) = t/n$, 
one obtains a representation of $N\subset M$ into $p(\Cal N \otimes \Bbb M_n(\Bbb C) \subset \Cal M \otimes \Bbb M_n(\Bbb C))p$ 
with exactness factor $p(\Cal P\otimes \Bbb M_n(\Bbb C))p=\Cal P^t$ (the $t$-amplification of $\Cal P$) and expectation
$\Cal E\otimes id(p \cdot p)$. It is easy to see that the isomorphism class of this representation does not depend on the choice of $n$ and $p$. 

We call this representation of $N\subset M$ the $t$-{\it amplification} of $\Cal N\subset^{\Cal E}\Cal M$ and denote it 
$\Cal N^t \subset^{\Cal E^t}\Cal M^t$. Its inclusion graph $\Lambda_{\Cal N^t \subset \Cal M^t}$ identifies naturally with $\Lambda_{\Cal N\subset \Cal M}$ 
and if one denotes $\oplus_{i\in I} \Cal B(\Cal H^t_i) = \Cal N^t \subset \Cal M^t =\oplus_{j\in J} \Cal B(\Cal H^t_j)$, then 
the entries of its dimension vectors are given by $\text{\rm dim}(_M\Cal H_j^t)=t \text{\rm dim}(_M\Cal H_j)$,  $\text{\rm dim}(_N\Cal K_i^t)=t \text{\rm dim}(_N\Cal K_i)$. 
The amplified representations $\Cal N^t \subset \Cal M^t$, $t >0$, are obviously stably isomorphic. 

Thus, if $\oplus_{i\in I} \Cal B(\Cal K_i)=\Cal N \subset^{\Cal E} \Cal M=\oplus_{j\in J} \Cal B(\Cal H_j)=\Cal M$ is an exact representation 
with finite couplings and $j_0\in J$ then, modulo stable equivalence, we may assume $\text{\rm dim}(_M\Cal H_{j_0})=1$. 
Thus, there exists a unique irreducible subfactor $P\subset M$ such that after reducing with the central support $p_{j_0}$ of $\Cal B(\Cal H_{j_0})\simeq \Cal B(L^2M)$ 
one has $P^{op}p_{j_0}=\Cal Pp_{j_0}$ and $\Cal N\subset \Cal M$ can be reconstructed from the Hilbert $M-P$ bimodule $_M(\Cal H_{j_0})_P=_ML^2M_P$, 
as explained in the last paragraph of 3.5, in which case we denote it $\Cal N_P \subset \Cal M_P$.

\vskip.05in
\noindent
{\bf 3.8.1. Notation}. We denote by $\Cal N^{u,fc} \subset \Cal M^{u, fc}$ the subrepresentation 
of $\Cal N^u \subset \Cal M^u$ obtained as the direct sum of irreducible representations that  
contain a type I$_\infty$ factor summand $\Cal B(\Cal H_j)$  of $\Cal M^u$ with dim$(_M\Cal H_j)=1$. Thus, 
 $\Cal N^{u,fc} \subset^{\Cal E} \Cal M^{u, fc}$ is of the form $\oplus_P (\Cal N_P \subset \Cal M_P)$, where the sum is 
over some irreducible subfactors $P\subset M$. We denote the inclusion graph of $\Cal N^{u,fc} \subset \Cal M^{u, fc}$  
by $\Lambda^{u,fc}_{_{N\subset M}}$. 

\vskip.05in

Taking into account that $_M(\Cal H_j)_P$, $_N(\Cal K_i)_P$ are all sub bimodules of $L^2(M_n), n\geq 0$, whose union is dense in $L^2(M_\infty)$, it follows 
that  $_ML^2(M_\infty)_P$ (resp. $_N(L^2(M_\infty)_P$) is the direct sum of irreducible $M-P$ (respectively $N-P$) Hilbert bimodules 
$\{\Cal H_j\}_j$ (respectively $\{\Cal K_i\}_i$), which appear with infinite multiplicity. Thus, when viewed as a direct summand of 
the universal exact representation $\Cal N^u \subset \Cal M^u$, $\Cal N_P \subset \Cal M_P$  is isomorphic to the 
$W^*$-inclusion $N \vee P^{op} \subset M\vee P^{op}$, acting on $L^2(M_\infty)$. 

In fact, the  entire Jones tower 
for $\Cal N_P=\Cal M_{-1} \subset \Cal M_0=\Cal M_P$ is represented on $L^2(M_\infty)$, with the consecutive inclusions 
$\Cal M_n\subset_{e_{n}} \Cal M_{n+1} =\langle \Cal M_n, e_n\rangle $ given by $M_n \vee P^{op}\subset_{e_n} M_{n+1} \vee P^{op}$,  
and expectation $\Cal E_{n+1}: \Cal M_{n+1}\rightarrow \Cal M_n$ given by $\Cal E_{n+1}(X)= \sum_km_k^n e_{n+1}Xe_{n+1} {m_k^n}^*$, 
$X\in \Cal M_{n+1}$, where $\{m_k^n\}_k$ is o.b. for $(\{e_m\}_{m\geq n+2})''$ over $(\{e_m\}_{m\geq n+3})''$, $N\subset M \subset_{e_0}M_1 \subset_{e_1} ...$ being the 
Jones tower and sequence of projections for $N\subset M$, acting (from the left) on $L^2(M_\infty)$.

Moreover, if one denotes $B=P'\cap M_\infty$, then the $\Cal M_n$'s are all contained in the II$_\infty$ von Neumann algebra $\langle M_\infty, e_B\rangle = (J_{M_\infty}B J_{M_\infty})'\cap \Cal B(L^2M_\infty)$,  with the canonical n.s.f. trace Tr=Tr$_{\langle M_\infty, B\rangle}$ (defined as usual by Tr$(xe_B y)=\tau(xy), x, y\in M_\infty$) being semifinite on each $\Cal M_n$ and preserving $\Cal E_n$, $n \geq 0$. 
In addition, in $\Cal B(L^2M_\infty)$ one has $\overline{\cup_n \Cal M_n}=\langle M_\infty, e_B \rangle $  and $\langle M_\infty, e_B \rangle$ naturally identifies this way with the enveloping von Neumann algebra $\Cal M_\infty$ 
of the Jones tower $\{\Cal M_n\}_n$, i.e., with the GNS completion of $(\cup_n \Cal M_n, \Phi)$, where $\Phi = \phi \circ \Cal E \circ \Cal E_1 \circ ... $, 
with $\phi$ any n.s.f. weight on $\Cal N$ (e.g. tracial; or a normal faithful state). 

To see this, it is sufficient to prove that $e_B$ is contained 
in $M\vee P^{op}$. Indeed, because then $ue_Bu^*$, $u\in \Cal U(M)$ are in $M\vee P^{op}$ and $\vee_u ue_Bu^*=1$, 
showing that Tr is semifinite on $M\vee P^{op}$. Then recursively, using that $\Cal E_n$ have finite index, Tr follows semifinite on $\Cal M_n$, $\forall n$. 

To see that $e_B\in M\vee P^{op}$, we show that the cyclic projection $[(M\vee P^{op})'(\hat{1})]$, which belongs to $\Cal M_P=M\vee P^{op}$, 
is equal to $e_B$. Indeed, by cutting with 
the orthogonal projection $e_{M_n}$, of $L^2(M_\infty)$ onto $L^2(M_n)$, which is invariant to $\Cal M_P=M\vee P^{op}$ and commutes with $e_B$ (because one has  
the commuting square relation $E_{M_n}E_{P'\cap M_\infty}=E_{P'\cap M_n}$), one gets $e_{M_n}(M\vee P^{op})'e_{M_n}(\hat{1})=(P'\cap M_n)(\hat{1}) \subset B(\hat{1})$ 
and the equality follows by letting $n\rightarrow \infty$.  

Note that by the commuting square relation $e_{M_n}e_B=e_Be_{M_n}=e_{P'\cap M_n}$,  by the formula for the canonical trace Tr on $\langle M_\infty, e_B\rangle$, 
and by the fact that $e_{M_n}\in (M \vee P^{op})'$, it follows that Tr restricted to $(M \vee P^{op})z_n$ (where $z_n\in \Cal Z(\Cal M_P)$ is the support of 
$\Cal M_P\ni x \mapsto xe_{M_n}$) coincides with Tr$_{\langle M, e_{P'\cap M_n}\rangle}$ and thus the minimal projections of its direct summands 
are proportional with the trace $\tau$ of the minimal projections in $P'\cap M_n \subset M_n$, which by  Section 3.3 are proportional to the entries 
of $\vec{d}_M$ supported by $z_n$. This shows that for any finite subset $J_0\subset J$  the restriction of 
Tr to $\{p_j\}_{j\in J_0}$, with $p_j$ minimal projection in $\Cal B(\Cal H_j)$,  is proportional to $\{\text{\rm dim}(_M\Cal H_j)\}_{j\in J_0}$. This 
in turn clearly implies that the weight vector for Tr on $M \vee P^{op}$ is proportional with the dimension vector $\vec{d}_M$. 

We summarize all these fact in the next statement. 

\proclaim{3.8.2. Proposition} Let $N\subset M$ be an extremal subfactor of finite index. Let $P\subset M$ be an irreducible subfactor and denote 
as above by $\Cal N_P\subset^{\Cal E} \Cal M_P$ the associated exact $(N\subset M)$-representation with finite couplings 
$\oplus_{i\in I} \Cal B(\Cal H_i)=N\vee P^{op} \subset^{\Cal E} M\vee P^{op}= \oplus_{j\in J} \Cal B(\Cal H_j)$, viewed as acting on $L^2(M_\infty)$, 
with $P^{op}=J_{M_\infty}PJ_{M_\infty}$.  Then we have: 
\vskip.05in 
$1^\circ$ $\Cal N_P \subset^{\Cal E} \Cal M_P$ is canonically Tracial. Moreover, the sequence 
of inclusions with Jones projections $\Bbb C=P'\cap M \subset P'\cap M_1 \subset_{e_1} M_2 \subset ...$ and the trace state $\tau$ inherited 
from $M_\infty$ forms a $\lambda$-sequence of inclusions, in the sense of Definition $2.7.4$ and Proposition $2.7.5$, with its pointed bipartite graph $\Gamma$ given by 
$\Lambda^t$ and weight vector $\vec{t}=(t_j)_j$ given by $t_j = \text{\rm dim}(_M\Cal H_j)$.

$2^\circ$ The Jones tower $\Cal N_P \subset^{\Cal E} \Cal M_P \subset^{\Cal E_1}_{e_0} \Cal M_1 \subset ... $ identifies with 
the tower $M_n \vee P^{op} \subset \Cal B(L^2M_\infty)$ with the enveloping von Neumann algebra $(\Cal M_\infty, \text{\rm Tr})$ 
identifying with $(\langle M_{\infty}, B \rangle, \text{\rm Tr}_{\langle M_{\infty}, B \rangle})$, where $B=P'\cap M_{\infty}$. 

$3^\circ$ $P^{op}$ is the RC-factor of $\Cal N_P \subset \Cal M_P$ $($i.e., $M'\cap \Cal N_P=P^{op})$ iff 
$P$ satisfies the bicommutant condition $(P'\cap M_\infty)'\cap M=P$. 
\endproclaim
\noindent
{\it Proof}. We already proved the first part of $1^\circ$ and part $2^\circ$. The second part of $1^\circ$ is an immediate consequence 
of Propositions 2.7.3 and 3.4.2. 

To prove $3^\circ$, recall from Sec. 3.6 
that $M'\cap (N \vee M^{op})=M^{op}$ (because $M^{op}$ is the RC-factor of the standard rep.). Thus 
$$
M'\cap (N\vee P^{op})=(M'\cap (N\vee M^{op}))\cap (N \vee P^{op})
$$
$$
=M^{op}\cap (N \vee P^{op})=M^{op}\cap (N\vee P^{op}) \cap ((P'\cap M_\infty)'\cap M_\infty)^{op}
$$
$$
= ((P'\cap M_\infty)'\cap M)^{op} \cap (N\vee P^{op})
$$
Thus, if $(P'\cap M_\infty)'\cap M=P$ then the last term is equal to $P^{op}$, showing that $M'\cap (N\vee P^{op})=P^{op}$. 

Conversely,  if $M'\cap (N \vee P^{op})=\Cal P$, with $P^{op}\subset \Cal P \subset M^{op}$ but $P^{op}\neq \Cal P$, then 
$\Cal P$ commutes with $(P'\cap M_\infty)^{op}$, so $((P'\cap M_\infty)'\cap M)^{op} \supset \Cal P$, showing that $(P'\cap M_\infty)'\cap M$ 
is strictly larger than $P$.  
\hfill $\square$

\vskip.05in

\noindent
{\bf 3.9. $A_\infty$-subfactors.} We say that $N\subset M$  
is an $A_\infty$-{\it subfactor} if its standard graph $\Gamma_{N\subset M}$ is equal to the bipartite graph 
$A_\infty$. This is equivalent to the ``dual'' standard graph $\Gamma'_{N\subset M}$ being equal to $A_\infty$, and also to the fact  
 that $[M:N]\geq 4$ and the higher relative commutants in the tower 
$N \subset M \subset_{e_0} M_1 \subset_{e_1} ...$ are generated by the Jones projections, i.e., 
$M_i'\cap M_j=\{1, e_{i+1}, ..., e_{j-1}\}''$, for all $j\geq i+2$. 

The existence of  subfactors with $A_\infty$ graph and index equal to $4$ of the hyperfinite II$_1$ factor $R$ was shown by Vaughan Jones in his original paper [J82], 
who gave two alternate constructions: $(a)$ as the inclusion $\{e_n\}_{n\geq 1}'' \subset \{e_n\}_{n\geq 0}''$, 
where $e_n$ is the sequence of Jones projections for the tower of $\Bbb C \subset \Bbb M_2(\Bbb C)$; $(b)$ as the inclusion of fixed point 
algebras of the product action of $SU(2)$ on $\overline{\otimes}_{n\geq 1}(\Bbb M_2(\Bbb C), tr)_n \subset  \overline{\otimes}_{n\geq 0}(\Bbb M_2(\Bbb C), tr)_n$. 
By (Corollary 5.2.1 in [P92a]) there exists in fact a unique index 4 subfactor of $R$ with $A_\infty$ graph,  up to conjugacy by an automorphism.  

The problem of the existence of $A_\infty$-subfactors with arbitrary index $\lambda^{-1}> 4$, which is closely related to the restriction problem 
on the Jones index for irreducible subfactors  beyond $4$, remained open for some time.  
This was solved in [P90], 
where it was shown that any $\lambda^{-1}>4$ can occur as the index of an $A_\infty$-subfactor. 
This was surprising, because initially it was believed that the set of all possible indices of  irreducible subfactors has gaps beyond $4$ (cf. page 940 of [J86]).

The $A_\infty$-subfactors in ([P90]) are canonical, universal objects, obtained through a tracial amalgamated 
free product construction, which we briefly recall.  Let $(\{e_n\}_{n\geq 0}, \tau)$ be a sequence of $\lambda$-Jones projections. Such a sequence  
exists by [J82],  because of the existence of subfactors of arbitrary index $\lambda^{-1}>4$ constructed as locally trivial, non-extremal subfactors 
of $R$, defined by $R_\lambda=\{x + \theta(x) \mid x\in pRp\}$, where $\theta:pRp\simeq (1-p)R(1-p)$, with $\tau(p)\tau(1-p)=\lambda$  
(due to the fact that by [MvN43]  the fundamental group of $R$ is the entire multiplicative group $(0, \infty)$). 

Let $Q$ be a diffuse tracial von Neumann 
algebra (taken as ``initial data''). Denote $A_{i, \infty}=\{e_j \mid j\geq i\}''$ and let $M^\lambda_\infty(Q)= Q \overline{\otimes} A_{1, \infty} *_{A_{1,\infty}} A_{0, \infty}$. 
Then define $M^\lambda(Q)$ as the smallest von Neumann algebra in $M^\lambda_\infty(Q)$ that contains $Q$ and is stable to 
the u.c.p. map $x \mapsto \sum_j m_j e_0xe_0m_j^*$, where $\{m_j\}_j$ is o.b. of $A_{1, \infty}$ over $A_{2, \infty}$, and let $N^\lambda(Q)=\{e_0\}'\cap M^\lambda(Q)$. 
Theorem 5.2 in [P90] then shows that $N^\lambda(Q)\subset M^\lambda(Q)$ is an inclusion of factors of index $\lambda^{-1}$ 
and standard graphs equal to $A_\infty$. 

By [PS01],  if $Q=L(\Bbb F_\infty)$ then $N^\lambda(Q),  M^\lambda(Q) \simeq L(\Bbb F_\infty)$, 
so $L(\Bbb F_\infty)$ contains $A_\infty$-subfactors of any index $\geq 4$. It is pointed out in (8.1 in [P90]; cf. Theorem in [P92b]; see also Theorem 4.5 in  [P00]) that 
given any extremal subfactor $N\subset M$ of index $\lambda^{-1}$, any free ultrafilter $\omega$ and any $Q\subset M^{\omega}$, 
the subfactor  $N^\lambda(Q)\subset M^\lambda(Q)$ can be realized as commuting square embedding into $N^\omega \subset M^\omega$. 

The construction in ([P90]) was further refined in [P94] to obtain a characterization of all lattices of tracial finite dimensional algebras with $\lambda$-Jones 
projections $\Cal G=(\{A_{ij}\}_{j\geq i\geq -1}, \{e_j\}_{j\geq 0}, \tau)$ that can occur as the standard invariant of a subfactor of index $\lambda^{-1}$, 
i.e., for which there exists an extremal subfactor $N\subset M$ with $\Cal G_{N\subset M}=\Cal G$. The 
abstract objects $\Cal G$ are called {\it standard $\lambda$-lattices}. Thus, the result in [P90] states that the ``minimal'' lattice,  
consisting of the algebras generated by the Jones projections, $A_{ij}=Alg(\{1, e_k \mid i+1 \leq k \leq j-1\})$, $\forall j\geq i \geq -1$, 
is a standard $\lambda$-lattice. This is called the {\it Temperley-Lieb-Jones} (abbreviated {\it TLJ}) 
$\lambda$-{\it lattice} and denoted $\Cal G^\lambda$. Thus, $A_\infty$-subfactors are the subfactors that have TLJ standard invariant. 

But the main interest for us in this paper is  the study  of $A_\infty$-subfactors of the hyperfinite factor $R$. 
By results in [H93], any irreducible subfactor $N\subset M$ of index $[M:N]\in (4, \frac{5+\sqrt{13}}{2})$ has $A_\infty$ graph. 
In particular, if $4< [M:N]<2+\sqrt{5}$ then $\Gamma_{N\subset M}=A_\infty$. So if $\alpha \in (4, 2+\sqrt{5})$ and $N\subset R$ is an 
irreducible subfactor of index $\alpha$ of the hyperfinite factor, then $\Gamma_{N\subset R}=A_\infty$. 
For $\alpha=\|E_{10}\|^2 \approx 4.062...$, which we saw in Section 2.10 is the first number in $\Bbb E^2$ that's larger than $4$, Ocneanu was able to solve the 
commuting square problem for $E_{10}$ (see 2.9), thus obtaining an example of irreducible hyperfinite subfactor of index $\alpha$ through the method described 
in Section 2.9 (see [S90] for details). 
 
The next two results detail some classes of representations that any $A_\infty$-subfactor has, starting with the standard representation:  

\proclaim{3.9.1. Proposition} Let $N\subset M$ be an $A_\infty$-subfactor and denote 
$\oplus_{l\in L}\Cal B(\Cal K_l)=\Cal N^{st}\subset \Cal M^{st} =\oplus_{k\in K} \Cal B(\Cal H_k)$ its standard representation, 
with its inclusion graph $\Lambda=\Lambda_{\Cal N^{st}\subset \Cal M^{st}}$.  
Then $\Lambda^t=A_\infty=\Gamma_{N\subset M}= (a_{kl})_{k\in K, l\in L}$. 
If one denotes $K=\{*=0, 2, 4, ...\}$ the ``even'' vertices of $\Gamma_{N\subset M}$ and $L=\{1, 3, ...\}$ its ``odd'' vertices, 
then the entries are given by $a_{2n, 2n+1}=1=a_{2n+2,2n+1}$, $n\geq 0$, with all other $a_{kl}=0$. 
The corresponding couplings are given by $d_n=\sqrt{P_n(\lambda)/\lambda P_{n-1}(\lambda)}$, for all $n\geq 0$, 
where $P_{-1}(\lambda)=1$, $P_0(\lambda)=1$, $P_{n+1}(\lambda)=P_n(\lambda)-\lambda P_{n-1}(\lambda)$, 
$\forall n \geq 1$. 
\endproclaim
\noindent
{\it Proof}. The formulas for the entries of the weights $d_n$, $n \geq 0$, are well known (see e.g., [GHJ88]). 
\hfill $\square$ 

\vskip.05in

Recall from ([J82]) that if $(\{e_n\}_{n\in \Bbb Z}, \tau)$ is 
the two-sided $\lambda$-sequence of Jones projections with $\lambda^{-1}>4$ and we denote 
$R_n=\{e_i \mid i\leq n-1\}''$ (resp. $P_n=\{e_i\mid i\geq n+1\}''$), $n \in \Bbb Z$, then  $(R_n\subset_{e_n} R_{n+1})_{n\in \Bbb Z}$
(resp. $(P_{n+1}\subset_{e_n} P_{n})_{n\in \Bbb Z}$) is  a Jones  tower/tunnel of II$_1$ factors of index $\lambda^{-1}$. Moreover, by 
(Corollary 5.4 in [PP84]; see also Corollary 3.3 in [P88]), $R_n\subset R_{n+1}$ (resp. $P_{n+1}\subset P_n$) is a 
locally trivial subfactor with $R_{n-1}'\cap R_n=\Bbb Cf_n +\Bbb C(1-f_n)$ (resp. $P_{n+1}'\cap P_{n}=\Bbb Cf'_n +\Bbb C(1-f'_n)$),   
where $f_n, f'_n$ are projections of trace $t<1/2$ and $t(1-t)=\lambda$. 

\proclaim{3.9.2. Proposition} Let $N\subset M$ be an $A_\infty$ subfactor and $N\subset M\subset_{e_0} M_1 \subset_{e_1} M_2 \subset_{e_2} ... \nearrow M_\infty$ 
its Jones tower and enveloping $\text{\rm II}_1$ factor. 
For each $i\geq 0$ denote $\tilde{M_i}=\{e_n\mid n\geq i+1\}'\cap M_\infty=P_i'\cap M_\infty=(M_i'\cap M_\infty)'\cap M_\infty$, with 
$\tilde{M}_0=\tilde{M}$, $\tilde{M}_{-1}=\tilde{N}$. 
Denote by $\tilde{E}_n: \tilde{M}_n \rightarrow \tilde{M}_{n-1}$ the map given by $x\mapsto \sum_jm^n_je_nxe_n{m_j^n}^*$, $x\in \tilde{M}_n$, where $\{m^n_j\}_j$ 
is an o.b. of $P_n$ over $P_{n+1}$. 
\vskip.05in
$1^\circ$ $\tilde{M}_n$ are $\text{\rm II}_1$ factors with $\tilde{M}_n=(M_n'\cap M_\infty)'\cap M_\infty$ and one has: 

\noindent
$(i)$ $\tilde{M}_{n-1}\subset  \tilde{M}_n$  is locally trivial with $\tilde{M}_{n-1}'\cap \tilde{M}_n=\Bbb Cf'_n + \Bbb C(1-f'_n)$; 

\noindent
$(ii)$ $M_i'\cap   \tilde{M}_j =  \tilde{M}_i'\cap   \tilde{M}_j =Alg(\{1, e_i, ..., e_j; f'_i, ..., f'_j\})$; 

\noindent
$(iii)$ $\tilde{E}_n$ is  a conditional expectation satisfying $\tilde{E}_n(f'_n)=1-t$; 

\noindent
$(iv)$ $\{\tilde{M}_n \subset_{e_n}^{\tilde{E}_{n+1}}\tilde{M}_{n+1}\}_{n\in \Bbb Z}$ is a Jones tower-tunnel; 

\noindent 
$(v)$ $\tilde{M}_n=vN(M_n, f'_n)$; 

\noindent
$(vi)$ $(N\subset^{E_N} M)\subset (\tilde{N} \subset^{\tilde{E}}\tilde{M})$ is a commuting square embedding, i.e., $\tilde{E}_{|M}=E_N$.

$2^\circ$ Let $P=\tilde{M}$ and consider the Hilbert $M-P$ bimodule $_M \Cal H_P=_ML^2(\tilde{M})_{\tilde{M}}$. The 
exact irreducible representation $\Cal N_P \subset^{\Cal E}\Cal M_P$ has graph $A_{-\infty, \infty}$ and is not Tracial. 

$3^\circ$ There exists a choice of a tunnel $....  \subset_{e_{-2}} M_{-1} \subset_{e_{-1}} M_0$ for $N=M_{-1}\subset M_0=M$ 
such that $R_0=\{e_j\mid j\leq -1\}''$ satisfies $R_0'\cap M=\Bbb C$ and more generally $R_0'\cap M_\infty=R_0'\cap R_\infty$. 
If one denotes by $\Cal N_{R_0} \subset \Cal M_{R_0}$ the exact  irreducible representation associated with the irreducible 
Hilbert $M-R_0$ bimodule $_ML^2M_{R_0}$, then $\Cal N_{R_0} \subset \Cal M_{R_0}$ has inclusion graph $\Lambda=A_{-\infty, \infty}$ and  
has finite couplings. Moreover, if one labels by $\Bbb Z$ the consecutive vertices of $\Lambda$, with even vertices $J=\{2n\mid n\in \Bbb Z\}$, 
odd vertices $I=\{2n+1 \mid n\in \Bbb Z\}$, then the vector $\vec{d}_M$ of $M$-couplings at even levels is given by $d_{2n}=(\frac{1-t}{t})^n$, $n\in \Bbb Z$. 
\endproclaim 
\noindent
{\it Proof}. Part $1^\circ$ is essentially (Theorem 7.8 in [P90]). 

$2^\circ$ The fact that if $P=\tilde{M}$ the inclusion graph of $\Cal N_P\subset^\Cal E \Cal M_P$,  
is equal to $A_{-\infty, \infty}$ follows immediately from $1^\circ$. The corresponding $(N\subset M)$-representation is non-Tracial 
because one has $\tilde{N}=\tilde{M}_{-1} \subset \tilde{M}_0=\tilde{M}$ as 
an intermediate inclusion, i.e., one has $(N\subset M) \subset (\tilde{N}\subset \tilde{M}) \subset (\Cal N_P \subset^{\Cal E} \Cal M_P$,  
with $\Cal E_{|\tilde{M}}=\tilde{E}_0$ and $f_0'\in \tilde{N}'\cap \tilde{M}$ satisfying $\tilde{E}_0(f_0')=1-t$, while $\tau(f_0')=t$. 

$3^\circ$ The existence of the choice of the tunnel $M=M_0\supset M_{-1} \supset ...$ so that $R_0=\overline{\cup_n M_{-n}'\cap M}=\{e_j \mid j\leq -1\}''$  
has trivial relative commutant in $M$, and more generally $R_0'\cap M_\infty=R_0'\cap R_\infty$, follows from (Theorem 4.10.$(b)$ in [P97a]; 
see also Lemma 2.7 in [P18a]). 

By part $1^\circ$ of Proposition 3.8.2, the calculation of the $M$-couplings amounts to the calculation of the trace vector for the 
$\lambda$-sequence of inclusions $\Bbb C = R_0'\cap M \subset R_0'\cap M_1 \subset_{e_1} R_0'\cap M_2 \subset ...$, which 
coincides with $R_0'\cap R_0 \subset R_0'\cap R_1 \subset R_0'\cap R_2 \subset ...$, for which this calculation was done in (Section 5 in [PP84]). 
\hfill $\square$ 

\vskip.05in
\noindent
{\bf 3.9.3. Remark.} Given an $A_\infty$-subfactor $N\subset M$, the irreducible inclusion of II$_1$ factors  $M\subset \tilde{M}=P_0'\cap M_\infty=(M'\cap M_\infty)'\cap M_\infty$ 
in Proposition 3.9.2 is quasi-regular, in the sense of (Definition 4.9 in [P97a], or Section 1.4.2 in [P01]). This follows trivialy from the fact that $_ML^2(\tilde{M})_M$ 
is a submodule of $_ML^2(M_\infty)_M$, which is a direct sum of finite index bimodules.  Using an argument similar to (proof of Theorem 4.5 in [P97a]), 
one can in fact show that if 
$M\subset \tilde{M}\subset \langle \tilde{M}, e\rangle \simeq M^\infty$ denotes the  basic construction, with its canonical n.s.f. trace Tr, 
then  $M'\cap \langle \tilde{M}, e\rangle$ is discrete abelian, generated by 
minimal projections $\{f_k\}_{k\in K}$, with $\text{\rm Tr}(f_k)=v_k$, where $(v_k)_k$ is the standard vector at even levels of  
the standard graph $\Gamma_{N\subset M}$, and that $_ML^2(\tilde{M})_M=\oplus_{k\in K} (_M{\Cal H_k}_M)$. 

\vskip.05in 

\noindent
{\bf 3.10. Untamed representations.} By Lemma 3.2.2, any non-degenerate embedding of a given subfactor $N\subset M$ into another 
subfactor $\tilde{N}\subset^{\tilde{E}} \tilde{M}$ composed with the standard representation of the latter, gives a 
representation of $N\subset M$. When combined with (Theorem 4.10 in [P00]), this gives the following class of representations, 
which we loosely call {\it untamed}:

\proclaim{3.10.1. Theorem} Let $N\subset M$ be an extremal inclusion of separable $\text{\rm II}_1$ factors of index $4<\lambda^{-1} = [M:N]<\infty$. 
Given any standard graph $(\Gamma, \vec{v})$ of index $\lambda^{-1}$, 
there exists a separable Tracial representation $(N\subset M)\subset (\Cal N\subset^{\Cal E}\Cal M)$ having 
$\lambda$-Markov weighted graph $(\Lambda_{\Cal N\subset \Cal M}, \vec{t})$ 
given by $({\Gamma}^t, \vec{v})$. 
\endproclaim
\noindent
{\it Proof}. Let $Q\subset P$ be an extremal inclusion of II$_1$ factors of index $\lambda^{-1}$ and 
weighted standard graph ($\Gamma, \vec{v})$. Both $\Cal G_{N\subset M}, \Cal G_{Q\subset P}$ contain the TLJ 
$\lambda$-lattice $\Cal G^\lambda$. So by (Therorem 4.10 in [P00]) there exists a separable extremal subfactor $\tilde{N}\subset \tilde{M}$ 
with standard invariant $\Cal G^\lambda$ in which both $N\subset M$ and $Q\subset P$ embed as non-degenerate 
commuting squares. Let $e=e^{\tilde{M}}_P$ denote the Jones projection for the basic construction 
$P\subset \tilde{M} \subset_e \langle \tilde{M}, P \rangle \simeq P^\infty$. By applying (Proposition 2.1 in [P92a]), 
one thus gets a non-degenerate embedding 
$$
(\tilde{N}\subset \tilde{M})\subset (\langle \tilde{N}, e \rangle \subset \langle \tilde{M}, e\rangle)\simeq (Q^\infty \subset^{E^\infty} P^\infty)
$$
where $(Q^\infty \subset^{E^\infty} P^\infty)=(Q\subset^{E^P_Q} P)\overline{\otimes} \Cal B(\ell^2\Bbb N)$. Thus, this inclusion 
further embeds with non-degenerate commuting squares into 
$(\Cal N\subset^{\Cal E}\Cal M)\overset\text{def}\to{=}(\Cal Q^{st}\subset \Cal P^{st})\overline{\otimes} \Cal B(\ell^2\Bbb N)$, 
which is a Tracial $\lambda$-Markov atomic $W^*$-inclusion with weighted inclusion graph  given by 
the standard weigthted graph of $Q\subset P$, i.e., $(\Gamma, \vec{v})$. Composing the embeddings, we get 
the representation $(N\subset M)\subset (\Cal N\subset^{\Cal E}\Cal M)$ which has all required properties. 
\hfill $\square$

\proclaim{3.10.2. Corollary} Any extremal inclusion of separable $\text{\rm II}_1$ factors $N\subset M$ 
of index $4<\lambda^{-1} = [M:N]<\infty$ has a Tracial representation with inclusion graph equal to $A_\infty$ 
and $\lambda$-weights as given in $3.9.1$. 
\endproclaim

\heading 4.  Amenability for graphs, subfactors, and $\lambda$-lattices \endheading 

In this section we recall the definition of amenability  for weighted bipartite graphs, standard $\lambda$-lattices and subfactors from ([P92a], [P93a], [P97a]), 
and revisit some results in ([P97a]) about these notions.

\vskip.05in

\noindent 
{\bf 4.1. Definition}. A Markov weighted graph $(\Lambda, \vec{t}, \lambda)$  is {\it amenable} if $\|\Lambda\|^2=\lambda^{-1}$. 

\vskip.05in

The next result establishes a F\o lner-type characterization of amenability for Markov weighted graphs. Like in (Definition 3.1 in [93a]),  
 if $\Lambda=(b_{ij})_{i\in I, j\in J}$  is a  bipartite graph and $F\subset J$ is a non-empty set, then we let $\Lambda^t\Lambda(F)=
 \{j \in J \mid \exists j'\in F$ such that $\sum_i b_{ij}b_{ij'}\neq 0\}$ and denote $\partial F \overset{\text{\rm def}}\to{=}  \Lambda^t\Lambda(F) \setminus F$.

\proclaim{4.2. Theorem} Let $(\Lambda, \vec{t})$ be a Markov weighted graph, with $\Lambda=(b_{ij})_{i\in I, j\in J}$ $\vec{t}=(t_j)_{j\in J}$ 
and $\Lambda\Lambda^t(\vec{t})=\lambda^{-1}\vec{t}$. Then $(\Lambda, \vec{t})$ is amenable if and only if it satisfies the following 
F\o lner-type condition: 

\vskip.05in
\noindent
$(4.2.1)$ For any  $ \varepsilon>0$, there exists $F=F(\varepsilon)\subset J$ finite such that $\|\vec{t}_{|\partial F}\|_2 
< \varepsilon\|\vec{t}_{|F}\|_2$, i.e., $(\sum_{j\in \partial F} t_j^2)^{1/2}$ $ < \varepsilon (\sum_{j\in F} t_j^2)^{1/2}$. 
\endproclaim
\noindent
{\it Proof}. In the case $\Lambda$ is the standard graph of a subfactor, this result amounts to ($(1') \Leftrightarrow (2)$ 
of Theorem 5.3 in [P97a]) and ($(i) \Leftrightarrow (ii)$ of Theorem 3.5 in [P93a]),  
whose proof only uses the fact that $\Lambda$ has non-negative integers as entries, the weight-vector $\vec{t}$ has strictly positive entries, 
and the fact that $\Lambda^t\Lambda(\vec{t})=\lambda^{-1}\vec{t}$. We revisit that argument, for convenience. 

For $a=(a_j)_j, b=(b_j)_j \in \Bbb C^J$, with $b$ finitely supported, we write $\langle a, b \rangle=\sum_j a_j \overline{b_j}$. 
 If condition $(4.2.1)$ is satisfied, then for any $\varepsilon >0$, the 
finite set $F=F(\varepsilon)\subset J$ gives rise to the finitely supported element $t_F=(t_j)_{j\in F} \in \ell^2J$ and  
if we denote $F'=\Lambda^t\Lambda(F)$, then we have 
$$
\lambda^{-1}\sum_{j\in F} t_j^2 = \langle \Lambda^t\Lambda(t_{F'}), t_F \rangle 
$$
$$
\leq \|\Lambda\|^2 \|t_{F'}\|_2 \|t_F\|_2\leq \|\Lambda\|^2 (1+\varepsilon)^{1/2} \|t_F\|_2^2, 
$$
where $\Lambda^t \Lambda$ is viewed here as an operator on $\ell^2(J)$ and $\| \ \|_2$ denotes the norm on this Hilbert space. 
Letting $\varepsilon \rightarrow 0$, this shows that $\|\Lambda\|^2\geq \lambda^{-1}$. Since we also have $\|\Lambda\|^2\leq \lambda^{-1}$   
(see Section 2.7), we get $\|\Lambda\|^2=\lambda^{-1}$. 

Conversely, assume $\|\Lambda\|^2=\lambda^{-1}$. 
Let $\Phi=\lambda T^{-1}\Lambda^t\Lambda T$, viewed as a $J \times J$ matrix with non-negative entries, 
where $T$ is the diagonal matrix with entries $\vec{t}=(t_j)_j$. Note that $\Phi$ defines a unital positive linear map from 
the semifinite von Neumann algebra $P=\ell^\infty J$ into itself. We endow $P$ with the n.s.f. trace 
Tr which on a finitely supported $a=(a_j)_j\in \ell^\infty J=P$ 
is given by Tr$(a)=\sum_j a_jt_j^2$. The $\lambda$-Markovianity condition for $(\Lambda, \vec{t})$ trivially implies $\text{\rm Tr}(\Phi(b)) 
=\text{\rm Tr}(b)$, $\forall b\in P$. By Kadison's inequality, this also implies $\|\Phi(b)\|_{2,\text{\rm Tr}} \leq \|b\|_{2,\text{\rm Tr}}$, 
$\forall b \in L^2(P, \text{\rm Tr})$. 

Since $\|\lambda\Lambda^t\Lambda\|=1$, it follows that for any $\delta>0$, there exists $F_0\subset J$ finite such that 
$T_0=_{F_0}(\lambda\Lambda^t\Lambda)_{F_0}$ satisfies $1 \geq \|T_0\|\geq 1-\delta^2/2$. Since $T_0$ is a  
symmetric $F_0\times F_0$ matrix with non-negative entries, by the classic Perron-Frobenius Theorem there exists 
$b_0 \in P_+$ supported by $F_0$ such that $\langle b_0, b_0 \rangle =1$ and $T_0b_0=\|T_0\|b_0\geq (1-\delta^2/2)b_0$. 
Thus, $\lambda\Lambda^t\Lambda (b_0) \geq (1-\delta^2/2)b_0$. 

Denote $b=T^{-1}(b_0)\in P_+$ and notice that $\|b\|^2_{2,\text{\rm Tr}}=1$. Also, we have:
$$
\|\Phi(b)-b\|_{2,\text{\rm Tr}}^2 \leq 2-2\text{\rm Tr}(\Phi(b)b) 
$$
$$
=2 -2 \langle \lambda\Lambda^t\Lambda(b_0), b_0 \rangle \leq 2-2(1-\delta^2/2)=\delta^2.
$$

Thus, $\Phi$ is a Tr-preserving unital c.p. map on $P$ and $b\in L^2(P, \text{\rm Tr})_+$ is a unit vector satisfying $\|\Phi(b)-b\|_{2,\text{\rm Tr}}
\leq \delta$, $\|\Phi(b)\|_{2, \text{\rm Tr}}\geq 1-\delta^2/2$. 

By (Appendix A.2 in [P97a]), when $\delta< 10^{-4}$ this implies there exists a spectral projection 
$e$ of $b$, corresponding to an interval $(c, \infty)$ for some $c>0$, such that $\|\Phi(e)-e\|_{2,\text{\rm Tr}}< \delta^{1/4} \|e\|_{2,\text{\rm Tr}}$. 
The latter inequality implies $\|(1-e)\Phi(e)\|_{2,\text{\rm Tr}}< \delta^{1/4}\|e\|_{2,\text{\rm Tr}}$. If one denotes by $F\subset J$ the support 
of $e\in P=\ell^\infty J$, this is easily seen to imply $\sum_{j\in \partial F} t_j^2 < \lambda^{-4} \delta^{1/4} \sum_{j\in F} t_j^2$. Thus, if one choses 
$\delta\leq (\lambda^4 \varepsilon^2)^4$ at the beginning, then $F$ satisfies $(4.2.1)$ for the given $\varepsilon >0$.   
\hfill $\square$ 

\vskip.05in

\noindent 
{\bf 4.3. Definitions}.  $1^\circ$ Let $N\subset M$ be an extremal inclusion of II$_1$ factors with finite index and $\Cal N^{st}\subset^{\Cal E^{st}} \Cal M^{st}$ 
its standard representation.  $N\subset M$ is {\it injective} if 
there exists a norm-one projection $\Phi:\Cal M^{st}\rightarrow M$ such that $\Phi(\Cal N^{st})=N$. It is easy to see that this is equivalent to the condition $\Cal E^{st} \circ \Phi 
=\Phi \circ \Cal E^{st}$.  We then also say that $\Phi$ is a norm-one projection of $\Cal N^{st}\subset^{\Cal E^{st}}\Cal M^{st}$ onto $N\subset M$. 

$N\subset M$ is {\it amenable} if the standard representation has an $(N\subset M)$-{\it hypertrace} i.e, a state $\varphi$ on $\Cal M^{st}$ 
that has $M$ in its centralizer and is $\Cal E^{st}$-invariant. 

$2^\circ$ A standard $\lambda$-lattice $\Cal G$ (resp. standard invariant $\Cal G_{N\subset M}$ of an extremal subfactor $N\subset M$) 
is {\it amenable} if its standard graph $\Gamma_{\Cal G}$ 
satisfies the Kesten-type condition $\|\Gamma_{\Cal G}\|^2=\lambda^{-1}$ (resp. $\|\Gamma_{N\subset M}\|^2=[M:N]$).

\proclaim{4.4. Proposition}  Let  $(N\subset M)\subset (\Cal N \subset^{\Cal E} \Cal M)$ be a 
non-degenerate  commuting square embedding of $N\subset M$ into a $W^*$-inclusion with expectation. There exists a norm one projection $\Phi:\Cal M\rightarrow M$ 
such that $\Phi \circ \Cal E=E_N \circ \Phi$ iff there exists a state $\varphi$ on $\Cal M$ that's $\Cal E$-invariant and has $M$ in its centralizer.  
\endproclaim 
\noindent
{\it Proof}. If $\Phi:\Cal M^{st}\rightarrow M$ is a norm one projection commuting with $\Cal E$ then  it is $M$-bimodular 
by Tomiyama's theorem and thus for any $x\in M, X\in \Cal M$,  
the state  $\varphi=\tau \circ \Phi$ 
satisfies $\varphi(xX)=\tau(\Phi(xX))=\tau(x\Phi(X))=\tau(\Phi(X)x)=\tau(\Phi(Xx)=\varphi(Xx)$. Also, $\varphi(X)= 
\tau(\Phi(X))=\tau(E_N(\Phi(X))=\tau(\Phi(\Cal E(X))=\varphi(\Cal E(X))$. 

Conversely,  if $\varphi$ is a state on $\Cal M$ that has $M$ in its centralizer and commutes with $\Cal E$ then 
one constructs a conditional expectation $\Phi$ from $\Cal M$ onto $M$ in the usual way: if 
$X\in \Cal M$, then $\Phi(X)$ is the unique element in $M$ with the property that $\tau(\Phi(X)x)=\varphi(Xx)$, $\forall x\in M$ (see e.g. Proposition 3.2.2 in [P92a]). 
Since $\varphi=\varphi\circ \Cal E$, it follows that 
$$
\Phi(\Cal E(X))=\varphi(\Cal E(X) \cdot )= \varphi(\Cal E(\Cal E(X) \cdot ))
$$
$$
=\varphi(\Cal E(X)\Cal E( \cdot ) =\varphi(X \Cal E( \cdot))=\varphi(XE_N(\cdot)) = E_N(\Phi(X)). 
$$
\hfill $\square$

\proclaim{4.5. Theorem}  Let $N\subset M$ be an extremal inclusion of separable $\text{\rm II}_1$ factors with finite index. The following 
conditions are equivalent: 

\vskip.05in 
$(1)$ $N\subset M$ is amenable. 

$(1')$ $N\subset M$ is injective. 

$(2)$ Any smooth representation $(N\subset M)\subset (\Cal N\subset^{\Cal E}\Cal M)$ has an $(N\subset M)$-hypertrace, i.e. an $\Cal E$-invariant state 
on $\Cal M$ that has $M$ in its centralizer. 

$(2')$ Any smooth representation  $(N\subset M)\subset (\Cal N\subset^{\Cal E}\Cal M)$ admits a norm one 
projection of $\Cal M$ onto $M$ that commutes with $\Cal E$.

$(3)$ $N, M$ are amenable $\text{\rm II}_1$ factors $($equivalently, $N \simeq R \simeq M)$ and 
its standard invariant $\Cal G_{N\subset M}$ is amenable $($i.e.,  $\|\Gamma_{N\subset M}\|^2=[M:N])$. 

$(4)$ $N\simeq R \simeq M$ and $\|\Lambda^{u,f}_{N\subset M}\|^2=[M:N]$, where $\Lambda^{u,f}_{N\subset M}$ denotes the inclusion graph 
of the universal exact finite representation of $N\subset M$ $($arising as direct sum of $\Cal N_P \subset \Cal M_P$, 
with $\Cal M_P=\oplus_j \Cal B(\Cal H_j)$, $\text{\rm dim}(_M{\Cal H_j}_P)<\infty, \forall j)$.  

$(5)$ $N\simeq R \simeq M$ and $\|\Lambda\|^2=[M:N]$ for any connected component $\Lambda$ of $\Lambda^{u,f}_{N\subset M}$. 

$(6)$ Given any finite set $F\subset M$ and any $\varepsilon > 0$, there exists a subfactor of finite index $P\subset N$ 
such that $F\subset_\varepsilon P'\cap M$. 

$(7)$ There exists a sequence of subfactors with finite index $M\supset N \supset P_1 \supset P_2 ...$ such that $P_n'\cap M \nearrow M$. 

$(8)$ There exists an $(N\subset M)$-compatible tunnel $M\supset N \supset P_1 \supset P_2 ...$ $($in the sense of Definition $2.1$ in $\text{\rm [P18a]})$ such that 
$P_n'\cap M \nearrow M$. 

$(9)$ $N\subset M$ is isomorphic to a model hyperfinite subfactor  
$N(\Cal G)\subset M(\Cal G)$, obtained as an inductive limit of 
higher relative commutants $($in $N$ and resp. $M)$ of an $(N\subset M)$-compatible tunnel $M\supset N \supset P_1 \supset P_2 ...$, 
whose choice is dictated by the standard invariant $\Cal G=\Cal G_{N\subset M}$. 
\endproclaim 
\noindent
{\it Proof}. By Proposition 4.2, we have $(1)\Leftrightarrow (1')$, $(2)\Leftrightarrow (2')$. One obviously has  
$(2)\Rightarrow (1)$ while $(1)\Rightarrow (2)$ follows from the implication $2)\Rightarrow 1)$ of (Theorem 7.1 in [P97a], see page 720 for the proof).

The implications $(8)\Rightarrow (7)\Rightarrow (6)$ are trivial 
and $(6)\Rightarrow (1)$ is a consequence of $4)\Rightarrow 5)\Rightarrow 6) \Rightarrow 7)$ in (Theorem 7.1 of [P97a], see page 720 for the proofs). 

Since $\|\Gamma_{N\subset M}\|=\|\Lambda_{\Cal N^{st}\subset \Cal M^{st}}\|$, the 
implication $(1) \Rightarrow (3)$ is a consequence of the general  Theorem 5.3 hereafter (note that it also  a consequence of 
Theorem 4.4.1 in [P92a] combined with Theorem 3.1 or 3.3 in [PP88]).

$(3)\Rightarrow (8)$ By the implication  $2)\Rightarrow 3)$ in (Theorem 7.1 of [P97a], see  top of page 719 for the proof), 
one first gets that for any $F\subset M$ finite and  any $\varepsilon >0$, there exists  an $(N\subset M)$-compatible subfactor $P\subset N$ and a  
finite dimensional subfactor $Q_0\subset P$ such that $F\subset_\varepsilon (P'\cap M)\vee Q_0$. 
Since $\Gamma_{N\subset M}$ amenable  implies $\Gamma_{P\subset M}$ amenable (cf. Corollary $6.6 (ii)$ in [P97a]; see also Proposition 2.6 in [P18a] for 
an alternative proof), we can apply this recursively to get a tunnel of $(N\subset M)$-compatible 
subfactors $M\supset N \supset P_1 \supset P_2 \supset ...$ and a sequence of commuting finite dimensional factors 
$Q_1, Q_2, ..., Q_n \subset \cap_{i=1}^n P_i$, such that $(P_n'\cap M)\vee Q_1 \vee ... \vee Q_n \nearrow M$. 
This implies $(P_n'\cap N)\vee Q_1 \vee ... \vee Q_n \nearrow N$ as well. It also implies  
that $(N\subset M)$ (or any other inclusion of hyperfinite factors with amenable graph) splits off $R$. This means that if we denote 
$M^0:=\vee_n(P_n'\cap M)$ and $N^0:=\vee_n(P_n'\cap N)$, then $(N\subset M)\simeq (N^0\subset M^0)\overline{\otimes}R$. 
But $(N^0\subset M^0)\simeq (N^0\subset M^0)\overline{\otimes}R$ (because $N^0\subset M^0$ are hyperfinite with same 
standard graph as $N\subset M$, thus amenable!). Hence, $(N\subset M)\simeq (N^0\subset M^0)$, which implies the approximation 
by higher relative commutants of $(N\subset M)=(N^0\subset M^0)$-compatible tunnels, required in $(8)$. 

This shows that $(1), (1'), (2), (2'), (3), (6), (7), (8)$ are equivalent. 

We further have $(2)\Rightarrow (5)$  by Theorem 5.3 below. The implication $(5)\Rightarrow (4)$ is trivial and $(4) \Rightarrow (3)$ 
by (Theorem 6.5 in [P97a]).  Thus, $(1)-(8)$ are equivalent. 

The equivalence of $(8)$ and $(9)$ is proved in (Remark 7.2.1 of [P97a]). But let us give here a more elegant argument, based on 
(Theorem 2.9 of [P18a]) and its proof, which allows deducing $(9)$ directly from $(3)$. 

Thus, we assume $\Cal G$ is an amenable standard $\lambda$-lattice with standard graph $\Gamma=\Gamma_\Cal G$ 
and canonical weights $\vec{v}$. Thus, $\|\Gamma\|^2=\lambda^{-1}$ and $\Gamma^t\Gamma \vec{v}=\lambda^{-1}\vec{v}$.  
By Theorem 4.2,  this condition is equivalent to the F\o lner property $(4.2.1)$ of the Markov weighted graph $(\Gamma, \vec{v})$. 
The proof of (Theorem 2.9 in [P18a]) shows that, given any (separable) subfactor $N\subset M$ with standard graph $\Cal G_{N\subset M}=\Cal G$, 
there exists a choice of an $(N\subset M)$-compatible tunnel $M\supset N \supset P_1 \supset P_2 ...$, such that 
if one denotes $(Q\subset R)=(\overline{\cup_nP_n'\cap N)}\subset \overline{\cup_n P_n'\cap M})$, then 
$[R:Q]=[M:N]$  and the higher relative commutants of $Q\subset R$ and resp. $N\subset M$ coincide, in fact 
$$
N'\cap M_n = Q'\cap M_n = Q'\cap R_n, \forall n. \tag 4.5.1
$$ 
In particular, $\Cal G_{Q\subset R}=\Cal G=\Cal G_{N\subset M}$. 

In order to satisfy condition  $(4.5.1)$, 
the $(N\subset M)$-compatible tunnel $M\supset N \supset P_1 \supset P_2 \supset ...$ 
that one takes depends on two types of choices, 
at each step $n$. 

Thus, if $M\supset N \supset ... \supset P_n$ have been already chosen, one next  takes $P_{n+1}$ to be 
a downward basic construction $P_{n+1} \subset P_n \simeq P_nq \subset qM_mq$, with $m\geq 1$ 
and $q\in P_n'\cap M_m$ appropriately chosen. The choice of $P_{n+1}$ is up to conjugacy by a unitary in $P_n$, but the way 
$m\geq 1$ and $q\in P_n'\cap M_m$ are chosen depends only on the properties of $\Cal G$, more precisely  
on the F\o lner-constants for $(\Gamma_{\Cal G}, \vec{v})$. 
It is important to note that this second type of choice, which depends only on $\Cal G$, can be taken the same for any $N\subset M$. 
 
 We call $M\supset N \supset P_1 \supset ...$  a {\it $\Cal G$-compatible tunnel}, and for each given $\Cal G$ we make once for all  a 
 choice for it, which we call the {\it model $\Cal G$-compatible tunnel}. 

In particular,  the isomorphism class of $Q\subset R$ is completely determined by $\Cal G$ and  $Q\subset R$  itself admits a  
model $\Cal G$-compatible tunnel $R\supset Q \supset P_1 \supset P_2...$ such that $P_n'\cap R\nearrow R$, $P_n'\cap Q\nearrow Q$. 
We denote this subfactor by $N(\Cal G)\subset M(\Cal G)$,  
calling it the {\it model subfactor} associated with the amenable standard $\lambda$-lattice $\Cal G$. 

Note that in case $\Cal G$ is both amenable and has ergodic core in the sense of ([P92a]), i.e., when $\Cal G =(A_{ij})_{j\geq i}$ is so 
that $Q=A_{1\infty}\subset A_{0\infty}=R$ are factors, then one can take the model $\Cal G$-compatible tunnel to be the Jones tunnel.

Now, given any hyperfinite subfactor $N\subset M$ with amenable standard graph $\Cal G_{N\subset M}=\Cal G$, then exactly the same proof 
as (Theorem 4.1 in [P93a]), shows that there exists a choice of the model $\Cal G$-compatible tunnel $M\supset N \supset P_1 \supset P_2 ...$ 
such that if one denotes $R_0=\cap_n P_n$ then $(P_n'\cap M) \vee R_0 \nearrow M$. Thus, $(N\subset M) \simeq (N(\Cal G) \subset M(\Cal G))\overline{\otimes} R_0$. 
Since $N(\Cal G) \subset M(\Cal G)$ splits off $R_0$, this shows that $(N\subset M)\simeq (N(\Cal G)\subset M(\Cal G))$.   
\hfill $\square$

\vskip.05in 
\noindent
{\bf 4.6. Remarks}. $1^\circ$ It is shown in (Theorem 7.5 in [P97a]) that if an extremal subfactor $N\subset M$ is amenable, in the sense of Definition 4.3.2$^\circ$, 
then any extremal subfactor $Q\subset P$ that can be embedded into it as a commuting square (not necessarily non-degenerate!) is amenable as well. 
In other words, if $N\simeq R \simeq M$ and $\|\Gamma_{N\subset M}\|^2=[M:N]$, then any 
commuting square sub-inclusion $Q\subset P$ of $N\subset M$ satisfies $\|\Gamma_{Q\subset P}\|^2=[P:Q]$ 
(in fact, the sub-inclusion $Q\subset P$ doesn't even need to be extremal, in general one has Ind$_{min}(Q\subset P)=\|\Gamma_{Q\subset P}\|^2$). 

This hereditary property is somewhat surprising. A key fact that allows proving this result is (Lemma 7.3 in [P97a]), 
which shows that if $Q\subset P$ is any finite index subfactor and $B$ is an arbitrary tracial von Neumann algebra that contains $P$, 
then the $^*$-subalgebra $B_0\subset B$ generated by $P$ and $\cup_k Q_k'\cap B$ is equal to sp$P(\cup_k Q_k'\cap B)R$, where 
$P\supset Q \supset Q_1 \supset ...$ is a tunnel for $Q\subset P$ and $R=\overline{\cup_k Q_k'\cap P}$. This easily implies that $_PL^2(B_0)_P$ is contained in 
$(\oplus_{k\in K^{st}} \Cal H_k)^{\oplus \infty}$, where $\Cal H_k$ is the list of irreducible Hilbert bimodules in the standard representation 
$\Cal Q^{st}\subset \Cal P^{st}$  of $Q\subset P$,  
allowing to show that the $Q-P$ bimodules and $P-P$ bimodules in $L^2(M_\infty)=\overline{\cup_n L^2M_n}$ give rise to  
a multiple of the standard representation of $Q\subset P$. The amenability of $N\subset M$ is then used to prove that  
$P'\cap M_\infty$ is big enough so that its commutant  in $M_\infty$ is ``locally'' approximately equal to $P$, a fact that allows constructing 
the $(Q\subset P)$-hypertrace on $\Cal Q^{st}\subset \Cal P^{st}$. 

$2^\circ$ Theorem 7.6 in [P97a] states that for an extremal subfactor $N\subset M$ 
the following three conditions are equivalent (formulated as such in that theorem): 

\vskip.05in
$1)$ $N\subset M$ amenable; 

$2)$ $\forall \varepsilon >0, \exists P\supset M$ hyperfinite  such that $\|\Lambda_{M'\cap P \subset N'\cap P}\|^2 \geq [M:N]-\varepsilon$;  

$3)$ $\|\Lambda^{u,rf}_{N\subset M}\|^2=[M:N]$;  
\vskip.05in
\noindent 
where $\Lambda^{u,rf}_{N\subset M}$ is the inclusion graph of the ``universal right-finite exact representation''
$\Cal N^{u,rf}\subset \Cal M^{u,rf}$ of $N\subset M$, obtained as the direct sum of $\Cal N_P \subset \Cal M_P$, with $P$  II$_1$ factors that contain $M$ as 
an irreducible subfactor, see 6.1.5 below for more about this sub-representation of the universal exact representation of a subfactor. 

However, while the proof of 
$1)\Leftrightarrow 2)$ and $1)\Rightarrow 3)$ are correct in [P97a], the proof of $3) \Rightarrow 2)$ on (bottom of page 724 of [P97a]) 
uses the fact that all $P$ in this construction are hyperfinite. So in order for that proof to work, one needs to modify the statement of (Theorem 7.6 in [P97a]) 
by replacing $\Cal N^{u,rf}\subset \Cal M^{u,rf}$ in condition $3)$ with $\Cal N^{u,hrf}\subset \Cal M^{u,hrf}$, defined similarly but with all 
$P$ taken $\simeq R$. 

On the other hand, one can prove the equivalence of $1)$ above with the following  
\vskip.05in
$2')$ $\forall \varepsilon >0, \exists P\supset M$ such that $\|\Lambda_{M'\cap P \subset N'\cap P}\|^2 \geq [M:N]-\varepsilon$ and $N, M\simeq R$;  

$3')$ $\|\Lambda^{u,rf}_{N\subset M}\|^2=[M:N]$ and $N, M\simeq R$;  
\vskip.05in
\noindent
We will detail the proof in a future paper.

$3^\circ$ In ([P97a]) there are two other interesting characterizations of amenability 
for a subfactor $N\subset M$,  that we did not include in Theorem 4.5 above. 
Thus, it is shown in (Theorem 7.1 in [P97a]) that $N\subset M$ is amenable iff its {\it symmetric enveloping} II$_1$ {\it factor} 
$M \underset{e_N}\to{\boxtimes} M^{op}$ is amenable (so $M \underset{e_N}\to{\boxtimes} M^{op}\simeq R$, by [C76]).  
And it is shown in (Theorem 8.1 in [P97a]) that $N\subset M$ is amenable iff the C$^*$-algebra generated in $\Cal B(L^2M)$ 
by $M, M^{op}$ and $e_N$ is simple (this is an ``Effros-Lance type'' characterization of amenability of $N\subset M$). 

$4^\circ$ Using Theorem 4.5 it is immediate to see that if  $N\subset M$ is amenable 
then any smooth commuting square embedding of $N\subset M$ into an arbitrary $W^*$-inclusion 
$\Cal N\subset^\Cal E \Cal M$ (with $\Cal N, \Cal M$ not necessarily atomic) has an $(N\subset M)$-hypertrace, i.e. an $\Cal E$-invariant state 
on $\Cal M$ that has $M$ in its centralizer.  Equivalently, any smooth commuting square embedding $(N\subset M)\subset (\Cal N\subset^\Cal E \Cal M)$,  admits a 
$\Cal E$-invariant norm-one projection of $\Cal M$ onto $M$. This condition actually appears as one of the 
equivalences in  (Theorem 7.1 of ([P97a]).

\heading 5.  Weak amenability for subfactors   \endheading

\noindent
{\bf 5.1. Definition}. An extremal subfactor of finite index $N\subset M$ is {\it weakly amenable} (respectively {\it weakly injective}) 
if it admits a Tracial representation $\Cal N\subset^{\Cal E}\Cal M$ that has an $(N\subset M)$-hypertrace (resp. a norm one projection). 

\proclaim{5.2. Proposition} $1^\circ$ If $N\subset M$ is amenable then it is weakly amenable. 
 
$2^\circ$ If $(Q\subset P)\subset (N\subset M)$ is a non-degenerate commuting square embedding of $\text{\rm II}_1$ factors 
and $N\subset M$ is weakly amenable, then $Q\subset P$ is weakly amenable. 
If in addition $[M:P]<\infty$ then $Q\subset P$ weakly amenable implies $N\subset M$  weakly amenable. 

$3^\circ$ Weak amenability is a stable isomorphism invariant: If $N\subset M$ is weakly amenable, then $(N\subset M)^t$, $\forall t>0$, and $(N\subset M)\overline{\otimes}R$ 
are weakly amenable. 
\endproclaim 
\noindent
{\it Proof}. Clear by the definitions. 
\hfill $\square$ 
\vskip.05in

We next show that weak amenability/injectivity, which follow equivalent by Proposition 4.4, are also equivalent to a 
Connes-F\o lner type condition, and they imply  the index of the subfactor must be equal to the square norm of the inclusion bipartite graph 
of the representation, thus belonging to the set $\Bbb E^2$.

\proclaim{5.3. Theorem} Let $N\subset M$ be an extremal inclusion of  {\rm
II}$_1$ factors and $\Cal N\overset{\Cal E}\to\subset\Cal M$ a Tracial representation. The following conditions are equivalent: 
\vskip.05in 
$(1)$ There exists a norm one projection of $\Cal N \subset^{\Cal E} \Cal M$ onto $N\subset M$ 

$(2)$ There exists a $(N\subset M)$-hypertrace on $\Cal N\subset^{\Cal E}\Cal M$. 

$(3)$ Given any finite set $F\subset \Cal U(M)$ and any $\varepsilon >0$, there exists a finite rank projection $p\in \Cal N$ 
such that $\sum_{u\in F}\|upu^*-p\|_{2, \text{\rm Tr}}< \varepsilon \|p\|_{2, \text{\rm Tr}}$. 

\vskip.05in 

Moreover, if the above conditions hold for $(N\subset M)\subset (\Cal N\subset^{\Cal E}\Cal M)$ 
and we denote by $(M_{i-1}\subset M_i \subset_{e_i} M_{i+1})\subset (\Cal M_{i-1}\subset^{\Cal E_i} \Cal M_i 
\subset_{e_i} \Cal M_{i+1})$, $i\in \Bbb Z$, the tower-tunnel of representations associated with it, then conditions 
$(1)-(3)$ hold for $(M_i \subset M_j)\subset (\Cal M_i \subset^{\Cal E_{ij}} \Cal M_j)$ for any $j>i$,  
where $M_0=M, M_{-1}=N$, $\Cal M_0=\Cal M$, $\Cal M_{-1}=\Cal N$ and $\Cal E_{ij}=\Cal E_{i+1} \circ ... \circ \Cal E_j$. 
\endproclaim
\noindent
{\it Proof}. We already proved the equivalence of conditions $(1)$ and $(2)$ in Proposition 4.4. 

Assume $(2)$ holds true and let $\varphi$ be an $\Cal E$-invariant state on $\Cal M$ that has $M$ in its centralizer.  

Let $F=\{u_1, ..., u_n\}$. Denote $\Cal
L=\{(\psi-\psi\circ\Cal E$, $\psi-\psi(u^*_1\cdot
u_1)$, $\psi-\psi(u^*_2\cdot u_2$, $\cdots$, $\psi-\psi(u^*_n\cdot
u_n))\in(\Cal M_*)^{n+1}\mid\psi$ a state in $\Cal M_*\}$. Note that $\Cal
L$ is a bounded, convex subset in $(\Cal M^*)^{n+1}=(\Cal M^{n+1})^*$ and
since the states $\varphi\in\Cal M_*$ are $\sigma(\Cal M^*,\Cal M)$ dense in
$S(\Cal M)$ it follows that the $\sigma((\Cal M^*)^{n+1},\Cal M^{n+1})$
closure $\bar{\Cal L}$ of $\Cal L$ contains all $(n+1)$-tuples $(\psi-
\psi \circ\Cal E, \psi-\psi(u^*_1\cdot
u_1), \cdots, \psi - \psi(u^*_n\cdot u_n))$ with $\psi\in S(\Cal M)$.

Taking $\psi=\varphi$, it follows that $\bar{\Cal L}$ contains $(\varphi-
\varphi \circ\Cal E, \varphi-\varphi (u^*_1\cdot
u_1), \cdots, \varphi-\varphi(u^*_n\cdot u_n))$.
But $\varphi(u\cdot u^*)=\varphi$, $\forall\ u\in\Cal U(M)$, so in
particular $\varphi-\varphi(u^*_i\cdot u_i)=0$, $i=1,2,\ldots,n$. Since we also have $\varphi \circ \Cal E =\varphi$, it follows that 
$(0,\ldots,0)=(\varphi-\varphi \circ\Cal E,\varphi-\varphi(u^*_1\cdot u_1),\ldots,\varphi-\varphi(u^*_n\cdot
u_n))\in\bar{\Cal L}$.

But since both $(0,\ldots,0)$  and $\Cal L$ are in $(\Cal M_*)^{n+1}$ and
since the dual of $(\Cal M_*)^{n+1}$ is $\Cal M^{n+1}$, it follows that the
$\sigma((\Cal M_*)^{n+1},\Cal M^{n+1})$ closure of $\Cal L$ in $(\Cal
M_*)^{n+1}$ is equal to the norm closure of $\Cal L$ and thus, $(0,\ldots,0)$
is norm adherent to $\Cal L$.

It follows that $\forall\ \delta>0$ there exists a state $\psi_0\in\Cal M_*$
such that
$$
\|\psi_0-\psi_0\circ\Cal E\|<\delta/3
$$
$$
\|\psi_0-\psi_0(u^*_i\cdot u_i)\| <\delta/3,  \quad 1\leq i\leq n.
$$

By replacing $\psi_0$ with $\psi=\psi_0\circ\Cal E$, it
follows tthat there exists a state $\psi\in\Cal M_*$ satisfying
$$
\psi =\psi\circ\Cal E 
$$
$$
\|\psi-\psi(u^*_i\cdot u_i)\|<\delta, \quad 1\leq i\leq n.
$$

Since $\Cal M_*=L^1(\Cal M,\text{Tr})$ and $L^1(\Cal M,\text{Tr})\cap\Cal M$
is dense in $L^1(\Cal M,\text{Tr})$, it follows that we may 
in addition assume there exists $b\in
L^1(\Cal M,\text{Tr})\cap\Cal M_+$ such that $\psi=\text{Tr}(\cdot b)$. 

Thus $\text{Tr}(b)=1$, $\text{Tr}(Xb)=\text{Tr}(\Cal E(X))b)$, $\forall\
X\in\Cal M$ and $\|\text{Tr}(\cdot b)-\text{Tr}((u^*_i\cdot u_i)b)\|<\delta$,
$1\leq i\leq n$. Since $\text{Tr}(\Cal E(X))b)=\text{Tr}(X(\Cal E(b)))$ and $\text{Tr}((u^*_i\cdot u_i)b)=\text{Tr}(\cdot
u_ibu^*_i)$, it follows that $b=\Cal E(b)$ and
$\|u_ibu^*_i-b\|_{1,\text{Tr}}<\delta$, $1\leq i\leq n$.

The first relation shows that $b\in\Cal N_{+}$ and the 
Powers-St\o rmer inequality applied to
the second shows that $a=b^{1/2}\in\Cal N_{+}$ satisfies
$\|a\|_{2,\text{Tr}}=1$ and $\|u_iau^*_i-a\|_{2,\text{Tr}}<\delta^{1/2}$, $\forall i$.

The Connes-Namioka trick (Theorem 1.2.1 in [C76]; or Section 2.5 in [C75]) then yields a spectral projection $p$ of $a$ such that
$\|  u_ipu^*_i-p\|_{2,\text{Tr}}<\delta^{1/2}\|p\|_{2,\text{Tr}}$, while
$p\in\Cal N$ (because $a\in\Cal N$). Taking $\delta=\varepsilon^2/|F|$ ends the proof or $(2) \Rightarrow (3)$. 

$(3)\Rightarrow (2)$. Assuming $(3)$, it follows that for each $F\subset \Cal U(M)$ finite, there exists $p_F \subset \Cal P(\Cal N)$ finite rank  
such that $\sum_{u\in F} \|up_Fu^* - p_F\|_{2, \text{\rm Tr}}< \frac{1}{|F|} \|p\|_{2, \text{\rm Tr}}$. 
Define $\varphi$ as $\varphi(X)=\text{\rm Lim}_F \text{\rm Tr}(Xp_F)/\text{\rm Tr}(p_F)$, $\forall X\in \Cal M$, where Lim$_F$ is a Banach limit 
over an ultrafilter majorizing the filter of finite subsets $F\subset \Cal U(M)$. It is then immediate to see that $\varphi$ this way defined has all $\Cal U(M)$ (thus all $M$) 
in its centralizer. Moreover, since the states Tr$(\ \cdot \ p_F)/\text{\rm Tr}(p_F)\in S(\Cal M)$ are $\Cal E$-invariant, $\varphi$ follows $\Cal E$-invariant. 

The last part is trivial and we leave it as an exercise. 
\hfill $\square$

\proclaim{5.4. Theorem} Let $N\subset M$ be an extremal inclusion of type {\rm
II}$_1$ factors. If $(N\subset M)\subset (\Cal N\subset^\Cal E \Cal M)$ is a Tracial representation for which there exists a norm-one projection onto $(N\subset M)$,  
then $\|\Lambda_{\Cal N \subset \Cal M}\|^2=[M:N]$.
\endproclaim
\noindent
{\it Proof}. By Proposition 2.7.1 we have $\|\Lambda_{\Cal N \subset \Cal M}\|^2\leq \text{\rm Ind}(\Cal E)=[M:N]$, so we only need to prove the opposite
inequality.

Let $e=e_{-1}\in M$ be a Jones projection, $\Cal
N_1\overset\text{def}\to=\{e\}'\cap\Cal N$ and $\Cal E_{-1}$ the conditional
expectation of $\Cal N$ onto $\Cal N_{-1}$ implemented by $e$ (see e.g.,
Proposition 2.2.4  in [P92a]). 

Let  $\varepsilon>0$. By the relative Dixmier property for $N\subset M$
(see Appendix A.1 in [P97a], or Theorem 1.1 in [P97b]), there exist unitary elements $u_1,\ldots,u_n\in N$ such
that
$$
\left\|\frac{1}{n}\sum_{i=1}^nu_ieu^*_i-\lambda1\right\|<\lambda\varepsilon/2.
$$

By Theorem 5.3, given any $\delta>0$ there exists a finite rank projection $p\in \Cal N_1$ such that 
$\|[u_i, p]\|_{2, \text{\rm Tr}}< \delta\|p\|_{2, \text{\rm Tr}}$
Let $w_i=pu_ip\in p\Cal Np$. Then
$$
\|w_iw^*_i-p\|_{2,\text{Tr}}=\|pu_ipu^*_ip-p\|_{2,\text{Tr}}=\|p(u_ipu^*_i-p)
p\|_{2,\text{Tr}} 
$$
$$
\leq \| u_ipu_i^*-p\|_{2, \text{Tr}} = \|[u_i, p]\|_{2, \text{\rm Tr}}< \delta \|p\|_{2,\text{Tr}}.
$$ 

A standard perturbation
argument (cf. e.g., Lemma A.2.1 in [P92a]) then shows that there exist unitary elements $v_i\in
p\Cal Np$ such that $\|v_i-w_i\|_{2,\text{Tr}}\leq
f(\delta)\|p\|_{2,\text{Tr}}$, with $f(\delta)$ a constant depending only on
$\delta$ and satisfying $f(\delta)\to0$ as $\delta\to0$.
This shows that for any 
$\delta'>0$ there exists a finite rank projection $p\in\Cal N_1$, $0\ne\text{Tr}p<\infty$, and
unitary elements $v_1,\ldots,v_n\in p\Cal Np$ such that $\|pu_i-v_i\|_{2, \text{\rm Tr}}<\delta'\|p\|_{2,\text{\rm Tr}}$, $\forall i$. 

So if we choose $\delta'>0$ such that $\delta'<\lambda\varepsilon/4$, then $p\in \Cal P(\Cal N_1)$,
$v_i\in\Cal U(p\Cal Np)$  satisfy
$$
\align
\left\|\frac{1}{n}\sum_{i=1}^nv_i(ep)v^*_i-\lambda
p\right\|_{2,\text{Tr}}&\leq\max_i\|pu_i-v_i\|_{2,\text{Tr}}+\left\|\frac{1}{n}
\sum_{i=1}^npu_iev^*_i-\lambda p\right\|_{2,\text{Tr}}\\
&\leq 2\max_i\|pu_i-v_i\|_{2,\text{Tr}}+\left\|p\!\left(\frac1{n}\sum_{i=1}^n
u_ieu^*_i-\lambda1\right)\!p\right\|_{2,\text{Tr}}\\
&\leq 2\max_i\|pu_i-v_i\|_{2,\text{Tr}}+\left\|p\!\left(\frac1{n}\sum_{i=1}^n
u_ieu^*_i-\lambda1\right)\right\|  \ \left\|p\right\|_{2,\text{Tr}}\\
&\leq (\lambda\varepsilon/2)\|p\|_{2,\text{Tr}}+(\lambda\varepsilon/2)\|p\|_{2,\text{Tr}}=
\lambda\varepsilon \|p\|_{2,\text{Tr}}. \endalign
$$

It follows that $p\Cal Np\subset^{\Cal E'} p\Cal Mp$ is a  finite dimensional 
$W^*$-inclusion with trace state 
$\tau=\text{Tr}(p)^{-1}\text{Tr}$ and $\tau$-preserving expectation $\Cal E'=\Cal E(p \cdot p)$, which has a projection
$e'=ep$ satisfying $\Cal E'(e')=\lambda p=\lambda1_{p\Cal Mp}$ and such that
$$
\left\|\frac1{n}\sum_{i=1}^nv_ie'v^*_i-\lambda1\right\|_2<\lambda\varepsilon
$$
for some unitary elements $v_i\in\Cal U(p\Cal Np)$. By  (Lemma 4.2 and
the estimate at the bottom of page 79 in [PP84]) it follows that if
$H(p\Cal Mp\mid p\Cal Np)$ denotes as usual the Connes-St\o rmer relative entropy
then
$$
H(p\Cal Mp\mid p\Cal
Np) \geq (1+\varepsilon^{1/2})^{-1}\ln\lambda^{-1}-(1+\varepsilon^{1/2})
\lambda^{-1}\eta(\lambda\varepsilon) \tag 5.4.1 
$$
$$=(1-\varepsilon)(1+\varepsilon^{1/2})^{-1}\ln[M:N]-(1+\varepsilon^{1/2})^{-1}
\eta(\varepsilon),
$$
where $\eta$ denotes here the function on the positive reals $\eta(t)=-t\ln t$, $t>0$. 

But by (Theorem 2.6 in [PP88]), if we denote by $\Lambda'$ the inclusion matrix for $p\Cal
Np\subset p\Cal Mp$ (which is thus a restriction of the inclusion matrix $\Lambda_{\Cal N\subset\Cal M}$) then
$$
\|\Lambda'\|^2\geq\exp(H(p\Cal Mp\mid p\Cal Np)) \tag 5.4.2 
$$

Finally, since $\|\Lambda\|\geq\|\Lambda'\|$, since $\varepsilon>0$ can be
taken arbitrarily small and since
$(1-\varepsilon)(1+\varepsilon^{1/2})^{-1}\to1$,
$(1+\varepsilon^{1/2})^{-1}\eta(\varepsilon)\to0$ as $\varepsilon\to0$, from $(5.4.1)$ and $(5.4.2)$  it
follows that
$$
\|\Lambda_{\Cal N \subset \Cal M}\|^2\geq [M:N].
$$
\hfill $\square$

\proclaim{5.5. Corollary} An extremal $\text{\rm II}_1$ subfactor  is weakly amenable if and only if it is weakly injective, 
and if these conditions are satisfied then the index of the subfactor lies in the set $\Bbb E^2$. 
\endproclaim

\noindent
{\bf 5.6. Remark}. A $W^*$-inclusion $\Cal Q \subset^{\Cal F} \Cal P$ is called {\it AFD} (abbreviated from {\it approximately finite dimensional}) 
if given any finite set $F\subset (\Cal P)_1$, any normal state $\varphi$ on $\Cal P$ with $\varphi \circ \Cal F = \varphi$ and any $\varepsilon >0$, 
there exists a finite dimensional $W^*$-inclusion $Q \subset^{E} P$ and a c.sq. embedding of it into $\Cal Q \subset^{\Cal F} \Cal P$
such that for any $x\in F$ there exists $y\in (P)_1$ with $\|x-y\|_{\varphi}< \varepsilon$. 
For an inclusion of II$_1$ factors with finite index $N\subset M$, this amounts to the definition of the AFD  property 
in ([PP88]): for any $F\subset (M)_1$ finite and $\varepsilon > 0$, there exists a commuting square $(Q\subset P)\subset (N\subset M)$ with $Q, P$ 
dimensional such that $\|x-E_{P}(x)\|_2\leq \varepsilon$, for all $x\in F$. 

By Theorem 4.5, amenable subfactors $N\subset M$ are AFD. The hyperfinite II$_1$ subfactors $P_{0\infty}\subset P_{1\infty}$ 
constructed from Markov cells as in Section 2.9 are obviously AFD.  By (3.1.3 in [PP92]), if $N\subset M$ is AFD with the finite dimensional 
approximants $Q\subset P$ so that $(Q\subset P)\subset (N\subset M)$ is non-degenerate, then $N\subset M$ has Tracial representations 
with $(N\subset M)$-hypertrace, so they are weakly amenable. Note that this non-degeneracy condition on the approximating finite dimensional 
subalgebras is automatic if $[M:N]$ is an isolated point in $\Bbb E^2\cap (4, 2+\sqrt{5})$. 

The various properties of representations (smoothness, exactness, Traciality, etc) allow defining several notions of ``weak amenability/injectivity'', where one 
requires existence of a $(N\subset M)$-hypertrace, respectively  of a norm-one projection, from one (or all) representation $(\Cal N\subset^\Cal E \Cal M)$ in some ``special'' 
class. Of particular interest is the following case. 

\vskip.05in 
\noindent
{\bf 5.7. Definition}.  We say that $(N\subset M)$ is {\it ufc-amenable} (resp. {\it ufc-injective}) if the universal exact finite-coupling 
representation $(\Cal N^{u,fc}\subset^{\Cal E^{u,fc}}\Cal M^{u,fc})$ admits an $(N\subset M)$-hypertrace (resp. a norm-one projection with range $N\subset M$). 

\vskip.05in

Note that for an extremal subfactor $N\subset M$ one obviously has 
``amenable $\Rightarrow$ ufc-amenable $\Rightarrow$ weakly amenable''. So in particular, if $N\subset M$ is ufc-amenable, 
then $[M:N]\in \Bbb E^2$. 

 While we will undertake a detailed study of this notion in a follow up to this article, we end this section 
by  stating without proof a result from this forthcoming paper,  
which relates ufc-amenabilty with several interesting structural properties of subfactors.

\proclaim{5.8. Theorem} Let $N\subset M$ be an extremal inclusion of separable $\text{\rm II}_1$ factors. The conditions $(1), (1'), (2), (3)$ below are equivalent 
and they imply condition $(4)$:

\vskip.05in

$(1)$ $N\subset M$ is ufc-amenable: the universal exact fc-representation $(N\subset M)\subset (\Cal N^{u,fc} \subset^{\Cal E^{u,fc}} \Cal M^{u,fc})$ 
admits an $(N\subset M)$-hypertrace. 

$(1')$ $N\subset M$ is ufc-injective:  there exists a norm-one projection of $(\Cal N^{u,fc} \subset^{\Cal E^{u,fc}} \Cal M^{u,fc})$ 
onto $(N\subset M)$.

$(2)$ $\|\Lambda_{N\subset M}^{u,fc}\|^2=[M:N]$ and $N\simeq R \simeq M$. 

$(3)$ $N\subset M$ has the AFDRC-property $($abbreviated for ``AFD by relative commu- tants''$)$: 
given any finite set $F\subset M$ and any $\varepsilon >0$, there exists a subfactor $Q\subset N$ such that $Q'\cap M$ is finite 
dimensional and $F\subset_\varepsilon Q'\cap M$.

$(4)$ $N\subset M$ has the  abc-property $($abbreviated for ``asymptotic bi-centralizer property''$)$: $M=(M'\cap N^\omega)'\cap M^\omega$. 
\endproclaim
\noindent

The proof of some of the implications in the above theorem are quite ellaborate, but let us point out right away that one has 
$(1) \Leftrightarrow (1')$ by Proposition 4.4 and $(1)\Rightarrow (2)$ by Theorem 5.4. Also,  $(3)\Rightarrow (1)$ 
has a proof similar to (3.1.3 of [P92a]). This entails $(3)\Rightarrow (2)$ as well, 
but note that this implication is also a direct consequence of (Theorem 3.3 in [PP88]). In addition, a proof of the implication $(3)\Rightarrow (4)$ 
can be easily completed along the lines of  (proof of Part (6) of Theorem 2.14 in [P91], on page 1676 of that paper).

\heading 6. Further  remarks and open problems \endheading 

\noindent
{\bf 6.1. General questions on} $W^*$-{\bf representations}. Representations of II$_1$ factors seem interesting to study in their own right. 
There is a large number of intriguing problems of ``general'' nature. We mention just a few.  

We have been able to construct only three types of $W^*$-representations for a given subfactor $N\subset M$: the ones in Example 3.1.2, coming from 
graphages $(Q\subset P)\subset (N\subset M)$; the exact representations in Section 3.6; the untamed $W^*$-representations 3.10. 
Representations in this last class tend to be non-smooth,  with infinite coupling constants,  in some sense ``uncontrollable''. 
For instance,  we saw that even if $N\subset M$ has 
finite depth with index $>4$, by Corollary 3.10.2 one can construct $W^*$-representations of $N\subset M$ that have $A_\infty$ inclusion graph.

\vskip.05in
\noindent
{\bf 6.1.1}. Find new constructions of 
Tracial representations with finite couplings. Conceive a method that produces all such representations for a given II$_1$ subfactor. 

\vskip.05in
\noindent
{\bf 6.1.2}. Is Traciality automatic for representations with finite couplings?  

\vskip.05in
\noindent
{\bf 6.1.3}.  Establish whether a representation given by a graphage is necessarily exact. Or at least that it necessarily has ``large'' RC-algebra. 

\vskip.05in
\noindent
{\bf 6.1.4}.  Assume $(N\subset M)\subset (\Cal N\subset^{\Cal E} \Cal M)$ is a  Tracial representation with finite couplings and 
irreducible finite inclusion graph $\Lambda_{\Cal N\subset \Cal M}$. Is this representation necessarily 
arising from a graphage, as in Example 3.1.2? 

\vskip.05in
\noindent
{\bf 6.1.5}. Do there exist representations $(N\subset M)\subset (\Cal N \subset^\Cal E \Cal M)$ 
with RC-factor $M'\cap \Cal N$ of type III (preferably exact) ? Do such representations exist for $M\simeq R$? 

\vskip.05in

Assume $(N\subset M)\subset (\Cal N_P\subset^{\Cal E}  \Cal M_P)$  is an irreducible exact ``right-finite'' representation,  
coming from an irreducible Hilbert-bimodule $_M\Cal H_P$ with $P$ a II$_1$ factor and dim$(\Cal H_P)<\infty$, say dim$(\Cal H_P)=1$.  
Then $_M\Cal H_P=$ $_ML^2P_P$ and the RC-envelope $\Cal P=M'\cap \Cal N_P$ of $P^{op}$ follows of the form $\tilde{T}^{op}$,  
for some  intermediate subfactor $M\subset P\subset \tilde{T} \subset \langle P, e_M\rangle$. But any such  subfactor  comes from an 
intermediate II$_1$ subfactor $M\subset T \subset P$ via basic construction,   
$\tilde{T}=\langle M, e_T\rangle$. 

Thus, by reducing with a finite projection in 
the type II factor $\tilde{P}^{op}=\tilde{T}^{op}$, it follows that $(N\subset M)\subset (\Cal N_P \subset \Cal M_P)$ is stably isomorphic 
to $(N\subset M)\subset (\Cal N_T\subset \Cal M_T)$, a representation that's still right-finite but this time 
the RC-envelope of $T$ is $M'\cap \Cal N_T=T^{op}$. Note that even if one starts with $M\subset P$ with $[P:M]=\infty$, 
during this process we may end up with $M\subset T$ satisfying $[T:M]<\infty$, so a 
sub-representation of $(\Cal N^{u,f}\subset \Cal M^{u,f})$. This justifies the following question:

\vskip.05in
\noindent
{\bf 6.1.6}.   Do there exist examples of irreducible exact right-finite representations $(N\subset M)\subset (\Cal N_P\subset \Cal M_P)$ 
arising from an irreducible embedding $M\subset P$ with $P$ a II$_1$ factor and $[P:M]=\infty$, 
such that $P$ equals its (exacting) RC-envelope ? In other words, while $(\Cal N^{u,f}\subset \Cal M^{u,f})$ is a subrepresentation 
of $(\Cal N^{u,rf}\subset \Cal M^{u,rf})$, are there examples where this inclusion is strict? 

\vskip.05in
\noindent
{\bf 6.1.7}.  Calculate the RC-factor (or exacting factor) for $\Cal N_P\subset \Cal M_P$, for a given irreducible 
embedding $P\subset M^\alpha$, $0< \alpha \leq \infty$. Along these lines, 
one can push the question 6.1.6 even further: Is any irreducible sub-representation 
of $\Cal N^u \subset \Cal M^u$  stably isomorphic to an exact representation with finite couplings? 
(i.e., a subrepresentation of $\Cal N^{u,fc}\subset \Cal M^{u,fc}$). Find concrete examples where 
$(N\subset M)\subset (\Cal N^u\subset \Cal M^u)$ is ``essentially'' equal/non-equal to $(N\subset M)\subset (\Cal N^{u,f}\subset \Cal M^{u,f})$. 

\vskip.05in
\noindent
{\bf 6.1.8}.  Do there exist representations $(N\subset M)\subset (\Cal N \subset^\Cal E \Cal M)$ with trivial RC-algebra, 
$M'\cap \Cal N=\Bbb C$ (preferably with finite couplings)? Do such representations exist for $N\subset M\simeq R$?

\vskip.05in

\noindent
{\bf 6.2. Values of index problems}.  The $W^*$-representation theory for a II$_1$ subfactor $N\subset M$ devised in this paper, 
following up on [P92a], is an analogue of the ``classic'' representation theory of a II$_1$ factor, and in fact it becomes just that when $N=M$. 
But while for a single II$_1$ factor $M$ its representations $M\subset \Cal B(\Cal H)$ are completely classified by the Murray-von Neumann dimension/coupling, 
$\text{\rm dim}(_M\Cal H)$, for a (non-trivial)  irreducible subfactor  $N\subset M$ the ``$W^*$-representation picture'' becomes strikingly complex. 
 It seems to us that it is this framework that's key for investigating rigidity paradigms concerning the values of the index of a given II$_1$ factor $M$, 
 constructed out of specific ``geometric data'', most notably for $M=R$ and for factors with Cartan subalgebras.

Despite the results in ([P90]), showing existence of large families of 
$A_\infty$-subfac-tors of any given index $\lambda^{-1}>4$, 
and more generally the results in [P94]  identifying the abstract objects $\Cal G$ that can occur as higher relative 
 commutants of subfactors (called standard $\lambda$-lattices in [P94]), 
a phenomenon such as $\mycal C(M)=\{4\cos^2(\pi/n) \mid n\geq 3\}\cup [4, \infty)$ seems to only occur 
when the II$_1$ factor $M$ comes from ``random-like'' constructions, such as the ones in ([P90], [P94], [PS01], [GJS11]). 
 In turn, for a II$_1$ factor $M$ with a ``very geometric background'', $\mycal C(M)$ seems more prone to be a subset of $\Bbb E^2$, $\Bbb E^2_0$, or even $\Bbb N$. 
 
 Recall in this respect the huge difference between the two types of existing results where $\mycal C(M)$ could be fully calculated, 
 with the case of the free group factor $M=L(\Bbb F_\infty)$ having $\mycal C(M)$ equal to the whole Jones spectrum $\{4\cos^2(\pi/n) \mid n \geq 3\}\cup[4, \infty)$ 
 by ([PS01]), and the case of free group measure space factors 
 $M= L^\infty X \rtimes \Bbb F_n$ having $\mycal C(M)$ equal to the semigroup of integers $\{1, 2, 3, ...\}$ by ([P01], [OP07], [PV11]).

The case of the hyperfinite II$_1$  factor $M\simeq R$  is particularly puzzling, as $R$  can be constructed  from very geometric  data (finitary, 
more generally amenable, due to [MvN43], [C76]), 
while at the same time it is the playing field for matrix randomness! However, the latter is always ``approximate randomness'', 
in moments. Our belief is that its ``geometric-finitary background'' prevails when it comes to index of subfactors problems.

\vskip.05in
\noindent
{\bf 6.2.1}.  Is $\mycal C(R)$ equal to $\Bbb E^2$?    
This question splits into two types of problems:

\vskip.05in

\noindent
$(a)$.   {\it The restrictions on the index problem}, asking whether the index of any irreducible hyperfinite subfactor is necessarily 
the square norm of a (possibly infinite) bipartite graph, i.e., that $\mycal C(R)\subset \Bbb E^2$. We believe quite strongly 
that this inclusion holds true. But there may be further restrictions for $\mycal C(R)$. 

$W^*$-representation theory 
should be quite useful in approaching  this problem. Since the inclusion  $\mycal C(R)\subset \Bbb E^2$ represents a condition only on the interval $(4, 2+\sqrt{5})$ and any 
irreducible subfactor with index in this interval has $A_\infty$ graph by ([H93]),  proving $\mycal C(R)\subset \Bbb E^2$ 
amounts to showing that any hyperfinite $A_\infty$-subfactor $N\subset R$ 
with index less than $2+\sqrt{5}$ satisfies $[R:N]\in \Bbb E^2$.

\vskip.05in
\noindent 
$(b)$ {\it The  commuting square problem}, asking whether for a given finite connected bipartite graph $\Lambda$, 
there exists a Markov cell $(P_{00}\subset P_{01})\subset (P_{10}\subset P_{11})$  as in Section 2.9, 
with the column $P_{00}\subset P_{10}$ having inclusion graph equal to $\Lambda$. While this would merely show 
$\mycal E(R)\supset \Bbb E_0^2$, since $\mycal E(M)\cap (4, 2+\sqrt{5}) = \mycal C(M) \cap (4, 2+\sqrt{5})$ 
for any II$_1$ factor $M$, it would still imply $\Bbb E_0^2\cap (4, 2+\sqrt{5}) \subset \mycal C(R)$. Ideally, the commuting square problem 
should be solved with control of the higher relative commutants of the resulting subfactor (in particular its irreducibility, see also 6.3.1).  

One should note that there is no known example of an irreducible hyperfinite subfactor with a non-algebraic number as index, 
nor in fact of any number $\mycal C(R)\setminus \Bbb E_0^2$.

\vskip.05in

\noindent
{\bf 6.2.2}.   {\it Conjecture}: Any hyperfinite $A_\infty$-subfactor $N\subset R$ is ufc-amenable, 
and thus any such subfactor satisfies $[N:R]=\|\Lambda_{N\subset R}^{u,fc}\|^2 \in \Bbb E^2$. 
From the above observations, this would imply $\mycal C(R)\subset \Bbb E^2$, thus solving $6.2.1(b)$ above.  

\vskip.05in 
\noindent
{\bf 6.2.3}.  More generally, we believe that if $M$ is any (separable) II$_1$ factor with a Cartan subalgebra, 
then $\mycal C(M)\subset \Bbb E^2$. 

This is of course verified by results in ([P01], [OP07], [PV11]), where for a large class of II$_1$ factors 
with Cartan decomposition one even has $\mycal C(M)\subset \{1, 2, 3, ...\}$. But in these cases the index rigidity is due to the uniqueness 
up to unitary conjugacy of the Cartan subalgebra (what is called $\Cal C_s$-rigidity in [PV11]), a property 
that many group measure space factors, such as $R$, do not have. Nevertheless, 
the presence of a Cartan subalgebra in a II$_1$ factor $M$ seems to make the $W^*$-representation theory of 
its subfactors $N\subset M$ be very ``structured'',  a phenomenon that should entail ``graph-like'' obstructions for $[M:N]$. 

\vskip.05in 
\noindent
{\bf 6.2.4}. Another intriguing question is the calculation of  $\mycal C(M)$ for the II$_1$ factors $M=L\Bbb F_n$ 
associated withe free groups with finitely many generators, $n=2, 3, ...$.  
Since $\mycal C(M)$ is invariant to amplifications of $M$, they are all equal (cf. [V88], [R92], [D92]). 
One would be tempted to believe that, due to its ``pure random nature'', $\mycal C(L\Bbb F_{fin})$ is equal to the entire Jones spectrum 
$\{4\cos^2(\pi/n) \mid n\geq 3\}\cup [4, \infty)$, like in the case $M=L\Bbb F_\infty$.  

But all calculations of $\mycal C(L\Bbb F_\infty)$ 
are based on variations of the construction in ([P90]), which uses amalgamated free product 
involving commuting squares associated with a $\lambda$-lattice (not necessarily standard) $\Cal G$ and some ``initial (semi)finite data'' $Q$. 
When the data $Q$ is $L\Bbb F_n$, $n\geq 1$, this allowed identifying the resulting subfactors $N^\Cal G(Q)\subset M^\Cal G(Q)$ as free group factors, 
using various models in Voiculescu's free probability theory ([R92], [PS01], [GJS10]). This does give $M\simeq L\Bbb F_{fin}$ when $\Cal G$ is  
finite (i.e., a Markov cell), but it always gives $M\simeq L\Bbb F_\infty$ if $\Cal G$ is not finite, notably if $\Cal G$ is  
the TLJ standard $\lambda$-lattice $\Cal G^\lambda$ (equivalently the TLJ standard $\lambda$-cell $\Cal C^\lambda$), with $A_\infty$ graph and index $\lambda^{-1}$. 

Consequently, the only known values $\lambda^{-1}\in \mycal C(L\Bbb F_{fin})$ are square norms of finite bipartite graphs $\Lambda$ for which 
the commuting square problem could be solved !  So exactly the same as the values $\lambda^{-1}$ that have been shown to exist in $\mycal C(R)$. 
This makes quite plausible the prediction $\mycal C(L\Bbb F_{fin})=\mycal C(R)$.   

\vskip.05in 
\noindent
{\bf 6.2.5}. Along the lines of Murray-von Neumann question of characterizing all multiplicative subgroups of $(0, \infty)$ that can be realized as 
fundamental groups of separable II$_1$ factors (where much progress has been done in [PV08]), a natural question is to characterize 
all sub-sets (resp. sub-semigroups) of the Jones spectrum $\{\cos^2(\pi/n) \mid n\geq 3\} \cup [4, \infty)$ that can be realized as 
$\mycal C(M)$ (resp. $\mycal E(M)$) for some separable II$_1$ factor $M$. To approach this, it  
may be useful to revisit the universal construction  in ([P90], [P00]) with the tools of deformation-rigidity theory in hand.  

\vskip.05in 
\noindent
{\bf 6.2.6}.  One obviously has  $\mycal C(Q\overline{\otimes} L\Bbb F_\infty) = \{4\cos^2(\pi/n) \mid n\geq 3\} \cup [4, \infty)$ for any II$_1$ factor $Q$.  
On the other hand, by (Theorem 1.3 in [PS01]) one has $\mycal C(N)\subset \mycal C(N * L\Bbb F_\infty)$ for  any II$_1$ factor $N$. 
Hence, $\mycal C((Q\overline{\otimes} L\Bbb F_\infty) *L\Bbb F_\infty )= \{\cos^2(\pi/n) \mid n\geq 3\} \cup [4, \infty)$
for any II$_1$ factor $Q$. It  would be interesting to find other classes of factors $M$ with $\mycal C(M)$ equal to the entire Jones spectrum. 

\vskip.05in

\noindent
{\bf 6.3. Actions of $\lambda$-lattices on $R$}. Given a II$_1$ factor $M$, we denote by $\mycal G(M)$ the set of all standard $\lambda$-lattices $\Cal G$ 
that can appear as the standard invariant of an extremal subfactor  $N\subset M$, $\Cal G=\Cal G_{N\subset M}$ 
(i.e., which in some sense can ``act'' on $M$). 
 
Note that this set encodes both $\mycal C(M)$ and $\mycal E(M)$, as one has $\mycal E(M)=\{\text{\rm Ind}(\Cal G) \mid \Cal G\in \mycal G(M)\}$ 
(where Ind$(\Cal G)=\lambda^{-1}$ is the {\it index of}  $\Cal G$), and $\mycal C(M)$ is the set of all Ind$(\Cal G)$, $\Cal G\in \mycal G(M)$ 
with $\Gamma_{\Cal G}$ having just one edge from its ``initial'' vertex $*$ (see 2.3). With this in mind, question 6.2.1 can be refined by asking the following:

\vskip.05in 
\noindent
{\bf 6.3.1}. Identify the set $\mycal G(R)$ of all standard $\lambda$-lattices 
that can  act on $R$. Or at least calculate the set $\mycal G_\alpha(R)$ of all $\Cal G\in \mycal G(R)$ 
with Ind$(\Cal G)\leq \alpha$, for some specific number $\alpha>4$, notably for $\alpha=2+\sqrt{5}$, or $\alpha=5$. 

At this point these questions seem extremely 
difficult to answer, with no indication of what the corresponding  sets might be. Note however that if Conjecture 6.2.1 is answered 
in the affirmative, with a solution to the commuting square problem 6.2.1$(b)$ for every finite connected bipartite graph $\Lambda$ 
(as per Section 2.9), then a next step would be to devise a method of constructing Markov cells, 
with given ``vertical'' bipartite graph $\Lambda$ as in 2.9, 
in a way that allows controlling (and computing!)  the standard invariant of the resulting subfactor. From that point on, the 
classification of all standard $\lambda$-lattices (planar algebras) in [JMS15] would help complete the picture at least 
for $\mycal G_5(R)$. 

Along these lines, the following question, complementing 6.2.2 above, is interesting to investigate: 

\vskip.05in

\noindent
{\bf 6.3.2.}  Does there exist  a hyperfinite subfactor $N\subset R$ 
with $[R:N]=\alpha$ and $A_\infty$ graph for any $\alpha\in \mycal C(R)\cap (4, \infty)$? 

\vskip.05in

As noticed in (Remarks $(3), (5)$ of 5.1.5 in [P92a]; see also 4.4 in [93a], 6.2 in [P93b]) 
our classification theorem for  hyperfinite subfactors with amenable 
standard invariant (7.2.1 in [P97a]; cf. also Theorem 4.5 in the present paper)
implies Ocneanu's Theorem [P85] on the uniqueness, up to cocycle conjugacy, of the free actions of a given finitely generated amenable group $\Gamma$ on the hyperfinite II$_1$ factor $R$. 
Jones obtained in [J81] a converse to Ocneanu's result, showing that any countable non-amenable group $\Gamma$ 
admits two actions on $R$ that are not cocycle conjugate. So it is quite natural to predict the following (see Problem 5.4.7 in [P92a]): 

\vskip.05in

\noindent
{\bf 6.3.3.}  {\it Conjecture}. Given any non-amenable $\Cal G \in \mycal G(R)$, there exist  subfactors $Q, P\subset R$ such that $\Cal G_{P\subset R}=\Cal G_{Q\subset R}$ 
but $(P\subset R)\not\simeq (Q\subset R)$. Taking into account ([BNP06]), one can even speculate that given any $\Cal G\in \mycal G(R)$ there exist 
infinitely/uncountably many hypefinite subfactors with $\Cal G$ as standard invariant. 

\vskip.05in 

\noindent
{\bf 6.3.4.} In the spirit of (5.4.3 in [P92a]), we denote by $\mycal C^a(M)$ (resp. $\mycal C^{fd}(M)$) the set of indices of irreducible subfactors 
with amenable (resp. finite depth) graph of the II$_1$ factor $M$. Similarly, we denote $\mycal G^a(M)$ (resp. $\mycal G^{fd}(M)$) the set of 
amenable (resp. finite depth) $\lambda$-lattices that can occur as standard invariants of subfactors of $M$. Note 
that $\mycal G^a(M)\subset \mycal G^a(R)$ for any II$_1$ factor $M$ (because by 7.2.1 in [P97] any amenable $\Cal G$ can be realized 
as the standard invariant of a hyperfinite II$_1$ subfactor; see also 4.4.2 in [P00] and Theorem 2.9 in [P18a]), so in fact we can denote 
$\mycal C^a(R)=\mycal C^a$, $\mycal C^{fd}(R)=\mycal C^{fd}$, $\mycal G^a(R)=\mycal G^a$, 
$\mycal G^{fd}(R)=\mycal G^{fd}$.  

It would be interesting to calculate such sets, or at least obtain some general properties/estimates, 
especially in the case $M=R$. For example, is the set $\mycal C^a$ (resp. $\mycal G^a$) 
countable/uncountable? Can it contain points in $\Bbb E^2\setminus \Bbb E_0^2$ (``limit points'')? 
We refer the reader to (5.4.3 in [P92a]) for a series of questions related to this. 
One should note the early results about $\mycal G^{fd}$ in ([Oc89], [H93]) 
and the more recent  complete description of $\mycal G_5^{fd}$ in ([JMS14]), finalizing a two decades long series of 
impressive results along these lines (see also [AMP15], were the description is pushed to $\mycal G_{5.25}^{fd}$). 

\vskip.05in 

\noindent
{\bf 6.3.5.}  Like for the set $\mycal C(M)$, there have been two types of complete calculations    
of $\mycal G(M)$ for a  II$_1$ factor $M$. On the one hand, $\mycal G(L\Bbb F_\infty)$ was shown in [PS01] to be equal to the set $\mycal G$ of all standard 
$\lambda$-lattices (as constructed in [P94]). Thus, $\mycal G(Q\overline{\otimes} L\Bbb F_\infty)=\mycal G$ and hence, by (Theorem 1.3 in [PS01]),  
$\mycal G ((Q\overline{\otimes}L\Bbb F_\infty)*L\Bbb F_\infty)=\mycal G$ for any II$_1$ factor $Q$. On the other hand, there have been several constructions 
of factors $M$ that have no $\lambda$-symmetries other than the trivial ones, corresponding 
to subfactors that  ``split-off''  $\Bbb M_n(\Bbb C), n\geq 1$ (cf. [V07], [PV21], [CIOS21]). Similar  techniques (deformation-rigidity theory) 
have been used to construct II$_1$ factors $M$ with various prescribed groups $G$ as outer automorphism group, Out$(M)=G$ (see e.g. [IPP05]). 
It would be interesting to obtain results of this type for prescribed ``small'' subsets $\Cal F$ of $\mycal G$.  Like in 6.2.4 above, this should be possible by combining 
the universal construction in [P00] with a way of adding to the building data of a  II$_1$ factor $M$ the 
``right amount'' of both rigid and soft ingredients, to show $\Cal F\subset \mycal G(M)$, then using deformation-rigidity 
to prove that any $\Cal G \in \mycal G(M)$ must lie in $\Cal F$. 

\vskip.05in

\noindent
{\bf 6.4. $W^*$-representations for subfactors of type III}. A type III factor $\Cal P$ is characterized by the property of being ``purely infinite'': 
all its ``parts are same as the whole'', all non-zero projections are equivalent. If $\Cal P$ has separable predual, an alternative characterization 
is that  any two (normal) representations of $\Cal P$ on separable Hilbert spaces (or left Hilbert modules) are unitary conjugate. 

For an inclusion of type III factors $\Cal Q \subset^{\Cal F} \Cal P$, endowed  with an expectation $\Cal F$ of finite index, 
$\text{\rm Ind}(\Cal F)=\lambda^{-1}<\infty$, the expectation $\Cal F$ that one usually considers is  the ``optimal one'' of minimal index ([Hi88]),  
uniquely determined by the condition that the values $\Cal F(q)$ it takes on the minimal projections $q$ in $\Cal Q'\cap \Cal P$ are proportional 
to $[q\Cal Pq:\Cal Qq]^{1/2}$ (see end of Section 2.5). 

Like in the II$_1$ case, a $W^*$-representation for a type III subfactor $\Cal Q\subset^{\Cal F} \Cal P$ 
is defined as a non-degenerate embedding into an atomic 
$W^*$-inclusion $\Cal N\subset^{\Cal E}\Cal M$.  So in this case, all the reps $\Cal P\hookrightarrow \Cal B(\Cal H_j)$ are equivalent. However,   
many of the general considerations in Section 3 work the same, with the obvious adjustments. 

An interesting problem here is to see whether there exist irreducible exact representation $\oplus_{i\in I}\Cal B(\Cal K_i)=\Cal N\subset^{\Cal E}\Cal M=\oplus_{j\in J}\Cal B(\Cal H_j)$ 
of $(\Cal Q\subset^{\Cal F}\Cal P)$ with $\Cal T=\Cal P'\cap \Cal N$ of type II$_1$. More generally, are there irreducible exact representations  
for which there are minimal central projections $q_j\in \Cal M$ 
such that $\Cal Mq_j=\Cal B(\Cal H_j)q_j$  corresponds to an irreducible $\Cal P - \Cal T$ bimodule that does not have finite index, i.e., $[\Cal T':\Cal M]<\infty$?

Viewing bimodules for the single factors $\Cal Q, \Cal P$ as endomorphisms and using the ensuing formalism of superselection sectors in ([L89], [EK97]) may be useful 
for the analysis of $W^*$-representations of type III subfactors.

\head References\endhead

\item{[AP17]} C. Anantharaman, S. Popa: ``An introduction to II$_1$ factors'',
\newline 
www.math.ucla.edu/$\sim$popa/Books/IIun-v13.pdf 
% \url{https://www.math.ucla.edu/~popa/Books/IIun.pdf}

\item{[AMP15]} N. Afzaly, S. Morrison, D. Penneys:
{\it The classification of subfactors with index at most $5\frac {1}{4} $},
arXiv: math.OA/1509.00038. 

\item{[BH94]} D. Bisch, U. Haagerup: {\it Composition of subfactors: new examples of infinite depth subfactors}, Ann. Scient. Ec. Norm. Sup.,
{\bf 29} (1996), 329-383.

\item{[BNP06]} D. Bisch, R. Nicoara, S. Popa:
{\it Continuous families of hyperfinite subfactors with the same
standard invariant}, Intern. Math. Journal, {\bf 18} (2007), 255-267
(math.OA/0604460).

\item{[CIOS21]} I. Chifan, A. Ioana, D. Osin, B. Sun:
{\it Wreath-like product groups and rigidity of their von Neumann algebras},
arXiv: math.OA/2111.04708.

\item{[C75]} A. Connes: {\it On the classification of von Neumann algebras and their automorphisms}, 
Symposia Mathematica, Vol. XX  (Convegno sulle Algebre C$^*$), INDAM, Rome, Academic Press. London 1976, pp. 435-478. 

\item{[C76]} A. Connes: {\it Classification of injective factors}, Ann. of Math., {\bf 104} (1976), 73-115.

\item{[CS75]} A. Connes, E. St\o rmer: {\it Entropy of automorphisms of }  II$_1$
{\it von Neumann algebras}, Acta Math., {\bf 134} (1975), 288-306. 

\item{[CDG82]}  D. Cvetkovic, M. Doob and I. Gutman: 
{\it On graphs whose spectral radius does not exceed $(2+5^{1/2})^{1/2}$}, 
Combinatoria {\bf 14} (1982), 225-239.

\item{[Dy92]} K. Dykema: {\it Interpolated free group factors}, Duke Math J. {\bf 69} (1993), 97-119.

\item{[EK98]} D. E. Evans, Y. Kawahigashi: ``Quantum symmetries on operator algebras'',
Oxford University Press, 1998. 

\item{[G60]} F. Gantmacher: ``The theory of matrices'', Vol 2., Chelsea, New York, 1960. 

\item{[GHJ89]} F. Goodman, P. de la Harpe, V.F.R. Jones: ``Coxeter graphs and towers of algebras'', MSRI Publications {\bf 14}, Springer Verlag 1989. 

\item{[GJS09]} A. Guionnet, V.F.R. Jones, D. Shlyakhtenko: 
{\it A semi-finite algebra associated to a planar algebra}, J. Func. Anal., {\bf 261} (2011), 1345-1360

\item{[H93]} U. Haagerup:
{\it Principal graphs of subfactors in the range $5 < [M:N]< 3 + \sqrt{5}$},
in ``Subfactors'',  Proc. Taneguchi Symposium in Oper. Algebras,
Araki-Kawahigashi-Kosaki Editors, World Scientific 1994, pp 1-38.

\item{[Hi88]} F. Hiai: {\it Minimizing indices of conditional expectations onto a
subfactor}, Publ. RIMS, Kyoto Univ., {\bf 24} (1988), 673-678.

\item{[Ho72]} A. J. Hoffman: {\it On limit points of spectral radii of non-negative symmetric integral matrices}, 
in Lecture Notes in Math. Vol. {\bf 303}, pp 165-172, Springer-Verlag, New-York/Berlin, 1972. 

\item{[IPP05]} A. Ioana, J. Peterson, S. Popa:
{\it Amalgamated Free Products of w-Rigid Factors and Calculation of their
Symmetry Groups},
Acta Math. {\bf 200} (2008), No. 1, 85-153 (math.OA/0505589).

\item{[Iz91]} M Izumi: {\it Application of fusion rules to the classification of subfactors}, Publ. RIMS {\bf 27} (1991), 953-994. 

\item{[J81]} V.F.R. Jones: {\it A converse to Ocneanu theorem}, J. Op. Theory {\bf 10} (1983), 61-64.

\item{[J83]} V. F. R. Jones: {\it Index for subfactors}, Invent. Math., {\bf 72}  (1983), 1-25. 

\item{[J86]} V.F.R. Jones: {\it Subfactors of type} II$_1$ {\it factors and related topics}, Proceedings ICM 1986, Vol. 2, pp 939-947. 

\item{[J99]} V.F.R. Jones: {\it Planar algebras}, arXiv: math.OA/9909027.

\item{[J00]} V.F.R. Jones:  {\it Ten problems}, in ``Mathematics: perspectives and frontieres'', pp. 79-91, AMS 2000, V. Arnold, M. Atiyah, P. Lax, B. Mazur Editors.

\item{[JMS14]} V.F.R. Jones, S. Morrison, N. Snyder: {\it The classification of subfactors of index at most} 5, 
Bull. Amer. Math. Soc. {\bf 51} (2014), 277-327. 

\item{[K86]} H. Kosaki: {\it Extension of Jones theory of index to arbitrary factors}, JFA {\bf 66} (1986), 123-140. 

\item{[L89]} R. Longo: {\it Index of subfactors and statistics of quantum fields}, Comm. Math. Physics {\bf 126} (1989), 217-247. 

\item{[MvN43]} F. Murray, J. von Neumann: {\it Rings of operators} IV, Ann. Math. {\bf 44} (1943), 716-808.

\item{[Oc85]} A. Ocneanu: ``Actions of discrete amenable groups on factors'', Lecture Notes in Math. 
{\bf 1138}, Springer Verlag, Berlin-Heidelberg-New York, 1985. 

\item{[Oc88]} A. Ocneanu: {\it Quantized groups, string algebras and Galois theory for von Neumann algebras},
in ``Operator Algebras and Applications'', London Math. Soc. Lect. Notes Series, Vol. 136,
1988,  pp. 119-172. 

\item{[OP07]} N. Ozawa, S. Popa: {\it On a class of } II$_1$
{\it factors with at most one Cartan subalgebra}, Annals of Math. {\bf 172} (2010), 101-137 (math.OA/0706.3623).

\item{[PP84]} M. Pimsner, S. Popa: {\it Entropy and index for subfactors}, Ann. Sci. Ecole Norm. Sup., {\bf 19} (1986), 57-106.

\item{[PP88]} M. Pimsner,  S. Popa: {\it Finite dimensional approximation for
pairs of algebras and obstructions for the index}, J. Funct. Anal., {\bf 98} (1991), 270-291.

\item{[P81]} S. Popa: {\it Orthogonal pairs of *-subalgebras in
finite von Neumann algebras}, J. Operator Theory, {\bf 9} (1983), 253-268.

\item{[P82]} S. Popa: {\it Maximal injective subalgebras in factors associated with free groups}, Advances in Math., {\bf 50} (1983), 27-48.

\item{[P88]} S. Popa: {\it Relative dimension, towers of projections and the commuting square problem}, Pac. J. Math., {\bf 137} (1989), 181-207. 

\item{[P89]} S. Popa: {\it Classification of subfactors: the reduction to commuting squares}, Invent. Math., {\bf 101} (1990), 19-43.

\item{[P90]} S. Popa: {\it Markov traces on universal Jones algebras and subfactors of finite index}, Invent. Math. {\bf 111} (1993), 375-405 
(IHES preprint No. 43/1990).

\item{[P91]} S. Popa: {\it Classification of actions of discrete amenable groups on amenable subfactors of type } II, 
International Journal of Math {\bf 21} (2010), 1663-1695 (IHES preprint 46/1992).

\item{[P92a]} S. Popa: {\it Classification of amenable subfactors of type} II, Acta Math., {\bf 172} (1994), 163-255.

\item{[P92b]} S. Popa: {\it Free independent sequences in type}  II$_1$ 
{\it factors and related problems}, Asterisque, {\bf 232} (1995), 187-202. 

\item{[P93a]} S. Popa:  {\it Approximate innerness and central freeness for subfactors: A classification result}, in
``Subfactors'' (Proc. Tanegouchi Symposium in Oper. Algebras), Araki-Kawahigashi-Kosaki Editors, World
Scientific 1994, pp 274-293.

\item{[P93b]} S. Popa: {\it Classification of subfactors and their endomorphisms}, CBMS Lecture Notes, {\bf 86} (1995).

\item{[P94]} S. Popa: {\it An axiomatization of the lattice of higher relative commutants of a subfactor}, Invent. Math., {\bf120} (1995), 427-445.

\item{[P97a]} S. Popa: {\it Some properties of the symmetric enveloping algebra
of a subfactor with application to amenability and property T}, Documenta Math., {\bf 4} (1999), 665-744.

\item{[P97b]} S. Popa: {\it The relative Dixmier property for inclusions of  von Neumann algebras}, Ann. Sci. Ec. Norm. Sup., {\bf 32} (1999), 743-767.

\item{[P00]} S. Popa: {\it Universal constructions of subfactors},
J. reine angew. Math. {\bf 543} (2002), 39-81. 

\item{[P01]} S. Popa: {\it On a class of type} II$_1$ {\it factors with
Betti numbers invariants}, Ann. of Math {\bf 163} (2006), 809-899
(math.OA/0209310; MSRI preprint 2001-024).

\item{[P18a]}  S. Popa:  {\it On the vanishing cohomology problem
for cocycle actions of groups on} II$_1$ {\it factors}, 
Ann. Ec. Norm Sup, {\bf 54} (2021), 409-445 (math.OA/180209964).

\item{[P18b]} S. Popa: {\it Coarse decomposition of } II$_1$ {\it factors},
Duke Math. J., {\bf 170} (2021) 
3073-3110 (math.OA/1811.09213).

\item{[PS01]} S. Popa, D. Shlyakhtenko: {\it Universal properties of}
$L(\Bbb F_\infty)$ {\it in subfactor theory}, Acta Mathematica, {\bf 191} (2003), 225-257.

\item{[PV08]} S. Popa, S. Vaes: {\it  Actions of $\Bbb F_\infty$ whose} II$_1$
{\it factors and orbit equivalence  relations have prescribed
fundamental group},  J. Amer. Math. Soc., {\bf 23} (2010), 383-403.

\item{[PV11]} S. Popa, S. Vaes: {\it Unique Cartan decomposition for} II$_1$ 
{\it factors arising from arbitrary actions of free groups}, Acta Math., {\bf 194} (2014), 237-284.

\item{[PV14]} S. Popa, S. Vaes:
{\it Representation theory for subfactors, $\lambda$-lattices and C$^*$-tensor categories}, 
Commun. Math. Phys., {\bf 340}  (2015), 1239-1280 
\newline
(math.OA/1412.2732).

\item{[PV21]} S. Popa, S. Vaes: $W^*$-{\it rigidity paradigms for embeddings of} II$_1$ {\it factors},
% \newline 
arXiv: math.OA/2102.01664.

\item{[PSt70]} R. Powers,  E. St\o rmer: {\it Free states of the canonical  anticommutation relations}, Comm. Math. Phys., 16 (1970), 1-33. 

\item{[R92]} F. Radulescu: {\it Random matrices, amalgamated free products and subfactors of the von Neumann algebra of a free group, of noninteger 
index}, Invent. Math., {\bf 115} (1994), 347-389.

\item{[S90]} J. Schou: {\it Commuting squares and index for subfactors},
\newline
arXiv: math.OA/1304.5907.

\item{[T72]} M. Takesaki: {\it Conditional expectations in von Neumann algebras},
JFA, {\bf 9} (1972), 306-321. 

\item{[T03]} M. Takesaki: ``Theory of Operator Algebras, II'', Encyclopedia of Math. Sciences, Vol. {\bf 125}, Springer-Verlag, Berlin, 2003. 

\item{[To71]} M. Tomiyama:  ``Tensor products and projections of norm one in von Neumann algebras'', Lecture Notes, Copenhagen University, 1971. 

\item{[V06]} S. Vaes: {\it Factors of type} II$_1$ {\it without non-trivial finite index subfactors}, 
Trans. Amer. Math. Soc. {\bf 361} (2009), 2587-2606.

\item{[Vo88]} D. Voiculescu: {\it Circular and semicircular systems and free product factors},
Prog. in Math. {\bf 92}, Birkh\"auser, Boston, 1990, pp. 45-60.

\item{[W88]} H. Wenzl: {\it Hecke algebras of type $A_n$ and subfactors}, Invent. Math. {\bf 92} (1988), 349-383. 

\enddocument